\documentclass[a4paper,12pt]{article} 
\usepackage{amssymb}
\usepackage{amscd}
\usepackage{amsmath} 
\usepackage{graphicx}
\usepackage{latexsym} 
\usepackage{enumerate}

\usepackage{theorem}
\theoremstyle{plain}
\theorembodyfont{\itshape}
\newtheorem{theorem}{Theorem}[section]
\newtheorem{proposition}[theorem]{Proposition}
\newtheorem{lemma}[theorem]{Lemma}
\newtheorem{corollary}[theorem]{Corollary}
\theorembodyfont{\rmfamily}
\newtheorem{definition}[theorem]{Definition}
\newtheorem{remark}[theorem]{Remark}
\newtheorem{example}[theorem]{Example}
\title{\bf Rational Surface Automorphisms \\
with Positive Entropy } 
\author{Takato Uehara \\ \\
Department of Information Engineering, Niigata University \\
8050, Ikarashi 2-no-cho, Nishi-ku, Niigata, 950-2181,  
Japan\thanks{E-mail address: {\tt uehara@ie.niigata-u.ac.jp}}} 
\begin{document}
\setlength{\lineskiplimit}{-7pt}
\setlength{\baselineskip}{15.6pt}
\maketitle
\begin{abstract} 
The aim of this paper is to construct many examples of rational surface 
automorphisms with positive entropy by means of the concept of orbit data. 
The concept enables us to introduce some mild and verifiable condition, and 
to show that if an orbit data satisfies the condition, 
then there exists an automorphism realizing the orbit data. 
Applying this result, we describe the set of entropy values of 
the rational surface automorphisms in terms of Weyl groups. 
\end{abstract} 
\section{Introduction} \label{sec:intro}
In this paper, we consider automorphisms on compact complex surfaces with 
positive entropy. According to a result of S. Cantat \cite{C}, 
a surface admitting an automorphism with positive entropy must be 
either a K3 surface, an Enriques surface, a complex torus or a rational 
surface. For rational surfaces, rather few examples had been known 
(see \cite{C}, Section 2). 
However, some rational surface automorphisms with invariant 
anticanonical curves have been constructed recently. 
Bedford and Kim \cite{BK1, BK2} found some examples of 
automorphisms by studying an explicit family of quadratic birational maps 
on $\mathbb{P}^2$, and then McMullen \cite{M} gave a synthetic construction 
of many examples. More recently, Diller \cite{D} sought automorphisms 
from quadratic maps that preserve a cubic curve 
by using the group law for the cubic curve. 
We stress the point that these automorphisms can be all obtained from 
quadratic birational maps. 
The aim of this paper is to construct yet more examples of rational surface 
automorphisms with positive entropy from general birational maps on 
$\mathbb{P}^2$ preserving a cuspidal cubic curve. 
\par 
Let $F : X \to X$ be an automorphism on a rational surface $X$. 
From results of Gromov and Yomdin \cite{G,Y}, the {\sl topological entropy} 
$h_{\mathrm{top}}(F)$ of $F$ is calculated as 
$h_{\mathrm{top}}(F)= \log \lambda(F^*)$, where $\lambda(F^*)$ is the 
spectral radius of the action 
$F^* : H^2(X;\mathbb{Z}) \to H^2(X;\mathbb{Z})$ 
on the cohomology group. Therefore, when handling the topological entropy of 
a map, we need to discuss its action on the cohomology group, which can be 
described as an element of a Weyl group acting on a Lorentz lattice. 
The {\sl Lorentz lattice $\mathbb{Z}^{1,N}$} is 
the lattice with the Lorentz inner product given by 
\[
\mathbb{Z}^{1,N}= \bigoplus_{i=0}^{N} \mathbb{Z} \cdot e_{i}, \quad \quad 
(e_i,e_j) = 
\left\{
\begin{array}{ll}
1 & (i=j=0) \\[2mm]
-1 & (i=j=1,\dots, N) \\[2mm]
0 & (i \neq j). 
\end{array}
\right. 
\]
For $N \ge 3$, the {\sl Weyl group $W_N \subset 
O(\mathbb{Z}^{1,N})$} is the group generated by $(\rho_i)_{i=0}^{N-1}$, where 
$\rho_{i} : \mathbb{Z}^{1,N} \to \mathbb{Z}^{1,N}$ 
is a reflection defined by 
\[
\rho_i(x)=x + (x,\alpha_i) \cdot \alpha_i, \quad \quad 
\alpha_i:= 
\left\{
\begin{array}{ll}
e_0-e_1-e_2-e_3 & (i=0) \\[2mm]
e_i-e_{i+1} & (i=1,\dots, N-1). 
\end{array}
\right. 
\]
We call the $W_N$-translate $\Phi_N:= \bigcup_{i=0}^{N-1} W_N \cdot \alpha_i$ 
of the elements $(\alpha_0, \dots,\alpha_{N-1})$ the {\sl root system } 
of $W_N$, and each element of $\Phi_N$ a {\sl root}. 
On the other hand, if $\lambda(F^*)>1$, then there is a blowup 
$\pi : X \to \mathbb{P}^2$ of $N$ points $(p_1, \dots, p_N)$ with $N \ge 10$ 
(see \cite{N1}), which gives an 
expression of the cohomology group : 
$H^2(X;\mathbb{Z})=\mathbb{Z} [H] \oplus \mathbb{Z} [E_1] 
\oplus \cdots \oplus \mathbb{Z} [E_N]$, 
where $H$ is the total transform of a line in $\mathbb{P}^2$, 
and $E_i$ is the exceptional divisor over $p_i$. 
Moreover, there is a natural marking isomorphism 
\begin{equation} \label{eqn:marking}
\phi_{\pi} : \mathbb{Z}^{1,N} \to H^2(X;\mathbb{Z}), \quad 
\phi_{\pi}(e_0) = [H], \quad \phi_{\pi}(e_i)=[E_i] ~~~ (i=1, \dots, N). 
\end{equation}
It is known (see \cite{N2}) that there is a unique element 
$w \in W_N$ such that the following diagram commutes: 
\begin{equation} \label{eqn:diag}
~~~~~~~~~~~~~~~~~~~~~~~~~~~
\begin{CD}
\mathbb{Z}^{1,N} @> w >> \mathbb{Z}^{1,N} \\
@V \phi_{\pi} VV @VV \phi_{\pi} V \\
H^2(X;\mathbb{Z}) @> F^* >> H^2(X;\mathbb{Z}). 
\end{CD}
~~~~~~~~~~~~~~~~~~~~~~~~~~~
\end{equation}
Then $w$ is said to be {\sl realized} by $(\pi,F)$ (see also \cite{M}). 
A question at this stage is whether a given element $w \in W_N$ 
is realized by some pair $(\pi,F)$. 
McMullen \cite{M} states that if $w$ 
has spectral radius $\lambda(w) > 1$ and no periodic roots, 
that is, 
\begin{equation} \label{eqn:McCondi}
w^{k}(\alpha) \neq \alpha \qquad 
(\alpha \in \Phi_N, \quad k \ge 1), 
\end{equation}
then $w$ is realized by a pair $(\pi,F)$. 
However, since the root system $\Phi_N$ is an infinite set when $N \ge 10$, 
it is rather difficult to see whether $w$ has no periodic roots. 
Indeed, he shows condition (\ref{eqn:McCondi}) only 
for the so called Coxeter element. 
One of our interest is to introduce 
a more verifiable condition and to construct realizations of much more 
Weyl group elements. 
\par 
Another interest is to find the entropy values of 
rational surface automorphisms. 
In general, the topological entropy of any automorphism 
$F : X \to X$ is expressed as $h_{\mathrm{top}}(F)= \log \lambda(w)$ 
for some $w \in W_N$ (see also Proposition \ref{pro:expent}). 
Namely, we have 
\[
\mathbb{E}:= \{h_{\mathrm{top}}(F) \, | \, F : X \to X 
\text{ is a rational surface automorphism} \} \subset \log \Lambda :=
\{\log \lambda \, | \, \lambda \in \Lambda\}, 
\]
where $\Lambda$ is given by 
\begin{equation} \label{eqn:PV}
\Lambda:=\{ \lambda(w) \ge 1 \, | \, w \in W_N, \, N \ge 3  \}.  
\end{equation}
The entropy values of automorphisms having been found so far 
seem to be contained in a very thin subset of $\log \Lambda$. 
On the other hand, 
one of our main results bridges the gap between two sets 
$\mathbb{E}$ and $\log \Lambda$, which is stated as follows. 
\begin{theorem} \label{thm:main0}
The logarithm of any value $\lambda \in \Lambda$ is realized 
by the entropy of some rational surface automorphism $F: X \to X$. 
In particular, we have
\[
\mathbb{E}=\log \Lambda. 
\]
\end{theorem}
We will show this theorem by introducing the concept of orbit data. 
\par
Now let us consider an $n$-tuple $\overline{f}=(f_1, \dots,f_n)$ of quadratic 
birational maps $f_\ell : \mathbb{P}_{\ell-1}^2 \to \mathbb{P}_\ell^2$ 
with each $\mathbb{P}_\ell^2$ being a copy of $\mathbb{P}^2$ and 
$\mathbb{P}^2_0=\mathbb{P}^2_n$. 
Note that the inverse of any quadratic map is also a quadratic map and a quadratic 
map has three points of indeterminacy. So the indeterminacy sets of $f_\ell^{\pm 1}$ 
can be denoted by 
$I(f_\ell)=\{ p_{\ell,1}^{+},p_{\ell,2}^{+},p_{\ell,3}^{+} \} 
\subset \mathbb{P}_{\ell-1}^2$ 
and 
$I(f_\ell^{-1})=\{ p_{\ell,1}^{-},p_{\ell,2}^{-},p_{\ell,3}^{-} \} 
\subset \mathbb{P}_{\ell}^2$
with a suitable matching of the indices $(\ell, j)$ between forward and backward 
indeterminacies to be specified later (see Section \ref{sec:const}). 
Then we assume that the orbit of each backward indeterminacy point reaches some 
forward one. More precisely, with the notation 
$p_{\iota}^{\pm}=p_{i,j}^{\pm}$ for 
$\iota \in \mathcal{K}(n) 
:=\{ \iota=(i,j) \, | \, i=1,2, \dots, n, \, j=1,2,3 \}$, 
suppose that there is a permutation $\sigma$ of $\mathcal{K}(n)$ and 
a function $\mu : \mathcal{K}(n) \to \mathbb{Z}_{\ge 0}$ such that 
the following condition holds for any $\iota \in \mathcal{K}(n)$: 
\begin{equation} \label{eqn:orbit}
p_{\iota}^{m} \neq p_{\iota'}^+ \quad 
(0 \le m < \mu(\iota), \, \iota' \in \mathcal{K}(n)), 
\qquad p_{\iota}^{\mu(\iota)}=
p_{\sigma(\iota)}^+, 
\end{equation}
where $p_{\iota}^{m}$ is defined inductively by 
\begin{equation} \label{eqn:traj}
p_{\iota}^{0}:=p_{\iota}^{-} \in \mathbb{P}_{i}^2, \qquad 
p_{\iota}^{m}:=f_{\ell} (p_{\iota}^{m-1}) \in \mathbb{P}_{\ell}^2 \quad
(\ell \equiv i+m~\text{mod}~n). 
\end{equation}
Then, by blowing up the orbit segments $p_{\iota}^{-}=p_{\iota}^{0}, 
p_{\iota}^{1}, \dots , p_{\iota}^{\mu(\iota)}=
p_{\sigma(\iota)}^+$ for $\iota \in \mathcal{K}(n)$, 
we can cancel all indeterminacy points of 
$(f_{\ell})$ and $(f_{\ell}^{-1})$. That is, if 
$\pi_{\ell} : X_{\ell} \to \mathbb{P}_{\ell}^2$ is a blowup of points 
$p_{\iota}^{m}$ 
with $0 \le m \le \mu(\iota)$ and $i+m \equiv \ell~\text{mod}~n$, 
then the birational maps 
$f_{\ell} : \mathbb{P}_{\ell-1}^2 \to \mathbb{P}_{\ell}^2$ 
lift to biholomorphisms 
$F_{\ell} : X_{\ell-1} \to X_{\ell}$, whose composition gives 
an automorphism $F := F_n \circ \cdots \circ F_1 : X_0 \to X_n=X_0$.   
Now, we denote by $\kappa(\iota)$ the number of points 
among $p_{\iota}^0, p_{\iota}^1,\dots, 
p_{\iota}^{\mu(\iota)}$ lying on $\mathbb{P}_n^2$  
or, in other words, 
$\kappa(\iota)=(\mu(\iota)+i-i_1+1)/n$ with 
$\sigma(\iota)=(i_1,j_1)$. It is easy to see that 
$\kappa(\iota) \ge 1$ provided $i_1 \le i$. 
This observation leads us to the following definition. 
\begin{definition} \label{def:data}
An {\sl orbit data} is a triplet $\tau=(n,\sigma,\kappa)$ consisting of 
\begin{itemize}
\item a positive integer $n$, 
\item a permutation $\sigma$ of 
$\mathcal{K}(n)$, 
and 
\item a function $\kappa : \mathcal{K}(n) \to \mathbb{Z}_{\ge 0}$ such that 
$\mu(\iota)= \kappa(\iota) \cdot n + i_1 - i-1 \ge 0$. 
\end{itemize}
\end{definition}
\begin{definition} \label{def:real'}
An $n$-tuple $\overline{f}=(f_1, \dots,f_n)$ of quadratic 
birational maps $f_{\ell}$ is 
called a {\sl realization} of an orbit data $\tau$ if 
condition (\ref{eqn:orbit}) holds for any 
$\iota \in \mathcal{K}(n)$. 
\end{definition}
A question here is whether a given orbit data $\tau$ admits some 
realization $\overline{f}$. 
\par
To answer this, we consider a class of birational maps preserving 
a cuspidal cubic $C$ on $\mathbb{P}^2$. 
Let $\mathcal{Q}(C)$ be the set of quadratic birational maps 
$f : \mathbb{P}^2 \to \mathbb{P}^2$ satisfying 
$f(C)=C$ and $I(f) \subset C^*$, where $C^*$ is the smooth locus of $C$. 
The smooth locus $C^*$ is isomorphic to $\mathbb{C}$ and is preserved by 
any map $f \in \mathcal{Q}(C)$. Thus, the restriction $f|_{C^*}$ is 
an automorphism expressed as 
\begin{equation*}
f|_{C^*} : \mathbb{C} \to \mathbb{C}, 
\quad t \mapsto \delta(f) \cdot t + c(f) 
\end{equation*}
for some $\delta(f) \in \mathbb{C}^\times$ and $c(f) \in \mathbb{C}$. 
For an $n$-tuple 
$\overline{f}=(f_1,\dots,f_n) \in \mathcal{Q}(C)^n$, 
the {\sl determinant} of $\overline{f}$ is defined 
by $\delta(\overline{f}):=\prod_{\ell=1}^n \delta(f_\ell)$. 
\par
Moreover, we take advantage of introducing the concept of orbit data 
to state a more verifiable condition than (\ref{eqn:McCondi}) 
in terms of a finite subset $\Gamma(\tau)$ of $\Phi_N$, that is,  
\begin{equation} \label{eqn:condi}
w_{\tau}^{k}(\alpha) \neq \alpha \qquad 
(\alpha \in \Gamma(\tau), \quad k \ge 1), 
\end{equation}
where $w_{\tau}$ is an element of $W_N$ with 
$N:=\sum_{\iota \in \mathcal{K}(n)} \kappa(\iota)$.
The orbit data $\tau$ canonically determines 
$\Gamma(\tau)$ and $w_{\tau}$, whose definitions will be given later 
(see Definitions \ref{def:latiso} and \ref{def:root}). 
It will be also seen later that any element $w \in W_N$ 
is expressed as $w=w_{\tau}$ 
for some orbit data $\tau$ (see Proposition \ref{pro:iden}). 
Thus, once an orbit data $\tau$ with $w=w_{\tau}$ is fixed, 
the finiteness of $\Gamma(\tau)$ enables us to check easily that 
$w$ satisfies condition (\ref{eqn:condi}). 
It should be noted that an expression for $w$ in terms of 
factorization into a product of the generators $(\rho_i)_{i=0}^{N-1}$ yields 
the orbit data $\tau=(n,\sigma,\kappa)$, 
and the number of $\rho_0$ in the expression is the length $n$, where 
$\rho_0$ corresponds to the standard cremona transformation.  
As an expression for a given element $w \in W_N$ is not unique, 
neither is an orbit data $\tau$ satisfying $w=w_{\tau}$.   
\par
Condition (\ref{eqn:condi}) is referred to as the {\sl 
realizability condition}, for reasons that become clear in the following 
theorem. 
\begin{theorem} \label{thm:main1}
Let $\tau=(n,\sigma,\kappa)$ be an orbit data with $\lambda (w_{\tau}) > 1$. Then, 
$\tau$ satisfies the realizability condition (\ref{eqn:condi}) 
if and only if there is a realization 
$\overline{f}_{\tau}=(f_1,\dots,f_n) \in \mathcal{Q}(C)^n$ of $\tau$ 
such that $\delta(\overline{f}_{\tau})=\lambda(w_{\tau})$. 
The realization $\overline{f}_{\tau} \in \mathcal{Q}(C)^n$ of $\tau$ with 
$\delta(\overline{f}_{\tau})=\lambda(w_{\tau})$ is uniquely determined. 
Moreover, $\tau$ determines 
a blowup $\pi_{\tau} : X_{\tau} \to \mathbb{P}^2$ of 
$N$ points 
on $C^*$ in a canonical way, which lifts 
$f_{\tau}:=f_n \circ \cdots \circ f_1$ to an automorphism 
$F_{\tau} : X_{\tau} \to X_{\tau}$: 
\begin{equation*} 
\begin{CD}
X_{\tau} @> F_{\tau} >> X_{\tau} \\
@V \pi_{\tau} VV @VV \pi_{\tau} V \\
\mathbb{P}^2 @> f_{\tau} >> ~\mathbb{P}^2. 
\end{CD}
\end{equation*}
Finally, $(\pi_{\tau},F_{\tau})$ realizes $w_{\tau}$ and 
$F_{\tau}$ has positive entropy 
$h_{\mathrm{top}}(F_{\tau})= \log \lambda(w_{\tau}) > 0$. 
\end{theorem}
As seen in Theorem \ref{thm:main3}, almost all orbit data satisfy 
the realizability condition (\ref{eqn:condi}). 
Furthermore, 
even if an orbit data $\tau$ does not satisfy the realizability condition 
(\ref{eqn:condi}), another orbit data $\check{\tau}$ with the same spectral 
radius, called the {\sl sibling} of $\tau$, 
does satisfy the condition. 
\begin{theorem} \label{thm:main2}
For any orbit data $\tau$ with $\lambda(w_{\tau}) > 1$, 
there is an orbit data $\check{\tau}$ satisfying 
$\lambda(w_{\tau})=\lambda(w_{\check{\tau}}) > 1$ and the realizability 
condition (\ref{eqn:condi}). 
In particular, $\check{\tau}$ is realized by 
$\overline{f}_{\check{\tau}}$. 
\end{theorem}
Theorem \ref{thm:main0} is a consequence of Theorem \ref{thm:main2} and 
the fact that any element $w \in W_N$ is expressed as $w=w_{\tau}$ 
for some $\tau$. 
\par
Finally, we give a sufficient condition for (\ref{eqn:condi}), which 
enables us to see clearly that almost all orbit data are realized, and to 
obtain an estimate for the entropy. 
\begin{theorem} \label{thm:main3}
Assume that an orbit data $\tau=(n,\sigma,\kappa)$ satisfies 
\begin{enumerate} 
\item $n \ge 2$, 
\item $\kappa(\iota) \ge 3$ for any 
$\iota \in \mathcal{K}(n)$, and 
\item if $\iota \neq \iota'$ satisfy 
$i_m=i'_m$ and 
$\kappa(\sigma^m(\iota))=\kappa(\sigma^m(\iota'))$ 
for any $m \ge 0$, then 
$\iota' \neq \sigma^m(\iota)$ for any $m \ge 0$, 
where $\sigma^m(\iota)=(i_m,j_m)$ and $\sigma^m(\iota')=(i_m',j_m')$. 
\end{enumerate}
Then the orbit data $\tau$ satisfies $2^n-1 < \lambda(w_{\tau}) < 2^n$ 
and the realizability condition (\ref{eqn:condi}). 
In particular, $F_{\tau}$ has positive entropy 
$\log (2^{n}-1) < h_{\mathrm{top}}(F_{\tau}) < \log 2^n$. 
\end{theorem}
Diller \cite{D} constructs, by studying single quadratic 
maps preserving $C$, automorphisms with positive entropy 
realizing orbit data $\hat{\tau}=(1,\hat{\sigma},\hat{\kappa})$. 
As is seen in Example \ref{ex:auto}, there is 
an orbit data $\tau$ such that $F_{\tau}$ is not topologically conjugate 
to the iterates of $F_{\hat{\tau}}$ that Diller constructs for any 
$\hat{\tau}=(1,\hat{\sigma},\hat{\kappa})$. 
Moreover, the element $w_{\tau}$ determined by 
this orbit data $\tau$ admits periodic roots and thus does not satisfy 
condition (\ref{eqn:McCondi}). 
\par
Any quadratic map in $\mathcal{Q}(C)$ is determined completely, 
up to a linear conjugation, by the configuration of the three indeterminacy points, 
which lie on the smooth locus $C^* \cong \mathbb{C}$, 
and the map is degenerate to a linear one when the indeterminacy points are collinear 
(see Lemma \ref{lem:quad} and Remark \ref{rem:explquad}).  
Hence for an orbit data $\tau=(n,\sigma,\kappa)$, 
the $3n$ conditions $p_{\iota}^{\mu(\iota)}=p_{\sigma(\iota)}^+$ 
determine $3n$ indeterminacy points $\{ p_{\iota}^+\}_{\iota \in \mathcal{K}(n)}$ and 
an $n$-tuple $\overline{f}=(f_1,\dots,f_n)$. 
Our investigations on the existence of a realization are divided into two steps. 
The first step is to check that $\tau$ admits a tentative realization, namely, 
each map $f_i$ is indeed quadratic (see also Definition \ref{def:tenta}).  
Proposition \ref{pro:tenta} states that a tentative realization $\overline{f}$ of $\tau$ exists 
if and only if $w_{\tau}$ has no periodic roots in a subset $\Gamma_{1}(\tau)$ of $\Phi$ 
with $n$ elements. 
When $\tau$ does not satisfy the condition, another orbit data 
$\check{\tau}=(\check{n}, \check{\sigma},\check{\kappa})$ with $\check{n} < n$ 
satisfies $\lambda(w_{\check{\tau}})=\lambda(w_{\tau})$ and admits 
a tentative realization. 
The second step is to check that the tentative realization $\overline{f}$ is indeed 
a realization of $\tau$, that is, check whether an orbit 
$\{ p_{\iota}^m \}_{m \ge 0}$ satisfies $p_{\iota}^m \neq p_{\iota'}^+$ 
before reaching $p_{\sigma(\iota)}^+$. 
Proposition \ref{pro:ind} shows that 
$\overline{f}$ is a realization of $\tau$ if and only if $w_{\tau}$ has 
no periodic roots in 
$\Gamma_2(\tau)=\Gamma(\tau) \setminus \Gamma_1(\tau)$. 
Even if $\tau$ does not satisfy the condition, 
$\overline{f}$ is a realization of another orbit data $\check{\tau}$. 
Note that  the configuration of the orbits 
$\{p_{\iota}^0, \dots, p_{\iota}^{\mu(\iota)}\}_{\iota \in \mathcal{K}(n)}$, 
which are blown up to yield an automorphism, 
is closely related to the eigenvalue problem of $w_{\tau}$ (see Proposition  \ref{pro:eigen}). 
Then the sibling $\check{\tau}$ of $\tau$ determines essentially 
the same configuration of the blown up points as $\tau$.  
On the other hand, under the assumptions in Theorem \ref{thm:main3}, 
Proposition \ref{pro:det} gives an estimate for the spectral radius 
$\lambda(w_{\tau})$ and shows the absence of periodic roots in 
$\Gamma_1(\tau)$, and then Proposition \ref{pro:cri} guarantees 
the absence of periodic roots in $\Gamma_2(\tau)$, 
which proves Theorem \ref{thm:main3}. 
\par
This article is organized as follows. 
After a preliminary study in Section \ref{sec:pre}, 
Section \ref{sec:const} is devoted to developing a 
method for constructing a rational surface automorphism 
from a realization of $\tau$. 
In Section \ref{sec:tenta}, we discuss the existence of a 
tentative realization of $\tau$, and 
in Section \ref{sec:real}, we investigate whether 
it is indeed a realization and prove Theorems \ref{thm:main0} and 
\ref{thm:main1}--\ref{thm:main3}. 
Finally, Propositions \ref{pro:det} and \ref{pro:cri} are 
proved in Section \ref{sec:proof}. 
\begin{remark} \label{rem:notation}
In the sequel, we will often use the following notations for a given orbit data $\tau=(n, \sigma,\kappa)$. 
\begin{enumerate}
\item For $\iota=(i,j) \in \mathcal{K}(n)$, put $\iota_m=(i_m,j_m)=\sigma^{m}(\iota)$. 
We will also use a similar notation for another orbit data 
(e.g. for $\iota'=(i',j') \in \mathcal{K}(n)$, put $\iota_m'=(i_m',j_m')=\sigma^{m}(\iota')$). 
\item For each $\iota \in \mathcal{K}(n)$, put $\bar{\sigma}(\iota):= \sigma^k(\iota)$, where 
$k \ge 0$ is determined by the relations $\kappa(\sigma^\ell (\iota))=0$ for $0 \le \ell < k$, and $\kappa(\sigma^k(\iota)) \ge 1$. 
\end{enumerate}
\end{remark}
\section{Preliminary} \label{sec:pre} 
In this section, we give a preliminary study of birational maps between 
complex surfaces, mainly  in order to clarify the meaning of the equality in (\ref{eqn:orbit}), 
since $p_{\iota}^m$ may be an infinitely near point on $\mathbb{P}^2$. 
To this end, from a birational map $f : X \to X$ on a compact surface $X$, 
we build up a new surface map 
$f_\sim : X_\sim \to X_\sim$. 
Although the surface $X_\sim$, which may be regarded as the set of all proper and 
infinitely near points on $X$, becomes noncompact and rather larger than $X$, 
it gives certain nice properties to the map $f_\sim$ as is mentioned below. 
These properties are used in our arguments of this article. 
\par
Now let $X$ be a smooth irreducible projective surface, and 
consider a pair $(x,\pi : \hat{X} \to X)$ of a point $x \in \hat{X}$ and 
a proper modification $\pi : \hat{X} \to X$. 
Two pairs $(x_1,\pi_1)$ and $(x_2,\pi_2)$ with $\pi_i : X_i \to X$ 
are said to be {\sl equivalent}, denoted by $(x_1,\pi_1) \sim (x_2,\pi_2)$, 
if $\pi_1^{-1} \circ \pi_2 : X_2 \to X_1$ 
is locally biholomorphic at $x_2$ and 
$x_1=\pi_1^{-1} \circ \pi_2(x_2)$. 
Let $X_\sim$ be the set of all equivalence classes of $(x,\pi)$. 
The equivalence class of $(x,\pi)$ is denoted by $[x,\pi]$. 
Then, $x \in X$ can be identified with $[x,\mathrm{id}_X] \in X_\sim$, 
which is said to be {\sl proper}. 
Moreover $X_\sim$ is equipped with the topology generated by 
\[
\{ U_\sim \, | \, [x,\pi : \hat{X} \to X] \in X_\sim 
\text{ and } U \subset \hat{X} \text{ is an open neighbourhood of } x \}, 
\]
where $U_\sim:=\{ [y,\pi] \in X_\sim \, | \, y \in U \}$. 
We notice that $U_\sim$ is identified with $U$, via $y \mapsto [y,\pi]$, 
which gives $X_\sim$ the complex structure induced from that on $X$, 
namely, $X_\sim$ becomes a (noncompact) complex surface. 
\par
Let $f : X \to Y$ be a birational map with its inverse $f^{-1} : Y \to X$. 
Moreover assume that $Y \subset \mathbb{P}^N$ and $f(X)$ is contained in 
no hyperplane of $\mathbb{P}^N$, so that $f^* H$ is a curve on $X$ 
for any hyperplane $H \subset \mathbb{P}^N$. 
Put 
\[
I(f):= \{ x \in X \, | \, 
x \in f^* H \text{ for any hyperplane } H \subset \mathbb{P}^N  \}. 
\]
We call $I(f)$ the {\sl indeterminacy set} of $f$, on which $f$ is not defined. 
\begin{remark} \label{rem:biholo}
Let $f : X \to Y$ be a birational map and $x \in X$ be a point. 
Then $f$ is locally holomorphic at $x$ if and only if $x \notin I(f)$. 
Moreover, $f$ is locally biholomorphic at $x$ if and only if 
$x \notin I(f)$ and $f(x) \notin I(f^{-1})$. 
\end{remark}
\par
For a curve $C \subset X$, let $\nu_{[x,\pi]}(C)$ be the multiplicity 
$\nu_{x}(\pi^{-1}(C))$ of $\pi^{-1}(C)$ at $x$ and put 
\[ 
C_\sim :=\{ [x,\pi] \in X_\sim \, | \, 
\nu_{[x,\pi]}(C) \ge 1 \} \subset X_\sim, 
\]
where $\pi^{-1}(C):=\overline{\pi^{-1}(C \setminus I(f))}$ 
is the strict transform of $C$. 
We notice that the definition of multiplicity is well-defined 
by virtue of the following lemma. 
\begin{lemma} \label{lem:independece}
If $(x_1,\pi_1)$ and $(x_2,\pi_2)$ are equivalent, 
then we have $\nu_{x_1}(\pi_1^{-1}(C))=\nu_{x_2}(\pi_2^{-1}(C))$. 
\end{lemma}
{\it Proof}. 
For $i=1,2$, assume that $U_i \subset X_i$ are open subsets containing $x_i$ 
such that $\pi_2^{-1} \circ \pi_1 : U_1 \to U_2$ is a biholomorphism 
with $x_2=\pi_2^{-1} \circ \pi_1(x_1)$. 
Since $U:=\pi_1(U_1) =\pi_2 (\pi_2^{-1} \circ \pi_1(U_1))=\pi_2(U_2) \subset X$, 
one may assume $I(\pi_i^{-1}) \cap U  \subset \{ x \}$ by taking sufficiently 
small $U_i$, where $x:= \pi_1(x_1)=\pi_2(x_2) \in X$. 
Note that $\nu_{x_i}(\pi_i^{-1}(C))$ equals the multiplicity of 
$\overline{\pi_i^{-1}(C \cap U \setminus \{ x \})} \subset U_i$ at $x_i$, and that 
$\overline{\pi_2^{-1}(C \cap U \setminus \{ x \})}= 
\overline{\pi_2^{-1} \circ \pi_1(\pi_1^{-1}(C \cap U \setminus \{ x \}))}=
\pi_2^{-1} \circ \pi_1(\overline{\pi_1^{-1}(C \cap U \setminus \{ x \})})$. 
Therefore, we have $\nu_{x_1}(\pi_1^{-1}(C))=\nu_{x_2}(\pi_2^{-1}(C))$ 
as $\pi_2^{-1} \circ \pi_1 : U_1 \to U_2$ is a biholomorphism. 
\hfill $\Box$ \par\medskip
\par
In what follows, two proper modifications 
$\pi_1 : X_1 \to X$ and $\pi_2 : X_2 \to X$ are identified if 
$\pi_1^{-1} \circ \pi_2 : X_2 \to X_1$ is an isomorphism. 
Under the identification, a blowup $\pi_x : X_x \to X$ of a point $x \in X$ 
is uniquely determined. 
Then, we have the following proposition. 
\begin{proposition}[\cite{A, B}] \label{prop:UP}
Let $f: X \to Y$ be a birational morphism of surfaces, 
and assume that $y \in I(f^{-1})$. 
Then, $f$ factorizes as 
\[
f : X  \overset{\hat{f}}{\longrightarrow} Y_{y}
\overset{\pi_y}{\longrightarrow} Y, 
\]
where $\pi_y : Y_y \to Y$ is the blowup of the point $y$, and 
$\hat{f} : X \to Y_y$ is a birational morphism. 
Moreover there exists a sequence of blowups $\pi_i : Y_i \to Y_{i-1}$ 
of points with $Y_0=Y$ and $Y_m=X$ such that $f= \pi_1 \circ \cdots \circ \pi_m$. 
\end{proposition}
\begin{proposition} \label{prop:iden}
Let $\pi_x : X_x \to X$ be the blowup of a point $x \in X$. 
Then there exists an isomorphism 
$(X_x)_\sim \to X_\sim \setminus \{ [x, \mathrm{id}_X] \}$, 
given by $[z,\pi] \mapsto [z,\pi_x \circ \pi]$. 
\end{proposition}
{\it Proof}. 
First we notice that $[z,\pi_x \circ \pi] \neq [x, \mathrm{id}_X]$ 
since $x \in I((\pi_x \circ \pi)^{-1})$. 
Now we construct the inverse of 
$(X_x)_\sim \to X_\sim \setminus \{ [x, \mathrm{id}_X] \}$. 
Take an element $[z,\hat{\pi}] \in X_\sim \setminus \{ [x, \mathrm{id}_X] \}$. 
Then we may assume that $x \in I(\hat{\pi}^{-1})$. 
Indeed, if $x= \hat{\pi}(z)$, then $\hat{\pi}$ is not locally biholomorphic 
at $z$ as $[z,\hat{\pi}] \neq [x,\mathrm{id}_X]$, 
which means that $x \in I(\hat{\pi}^{-1})$. 
On the other hand, if $x \neq \hat{\pi}(z)$ and $x \notin I(\hat{\pi}^{-1})$, 
that is, $\hat{\pi}^{-1}$ is locally biholomorphic at $x$, then put 
$\bar{\pi}:= \hat{\pi} \circ \pi_y$ with $y:= \hat{\pi}^{-1}(x)$. 
Since $z \neq y$, $\bar{z}:= \pi_y^{-1}(z)$ satisfies 
$[\bar{z},\bar{\pi}] = [z,\hat{\pi}]$ and $x \in I(\bar{\pi}^{-1})$. 
Hence we can assume that $x \in I(\hat{\pi}^{-1})$. 
Proposition \ref{prop:UP} says that $\hat{\pi}$ factorizes as 
$\hat{\pi}=\pi_x \circ \pi$ for some proper modification $\pi$, 
which yields inverse of 
$(X_x)_\sim \to X_\sim \setminus \{ [x, \mathrm{id}_X] \}$. 
Therefore the proposition is established. 
\hfill $\Box$ \par\medskip
For two pairs $[x_i,\pi_i : X_i \to X]$, we put 
$[x_1,\pi_1] < [x_2,\pi_2]$ if $\pi_1^{-1} \circ \pi_2 : X_2 \to X_1$ 
is locally {\sl holomorphic} (not necessarily locally biholomorphic) 
at $x_2$ and $x_1=\pi_1^{-1} \circ \pi_2(x_2)$. 
Note that the definition is independent of the choice of 
$(x_i,\pi_i) \in [x_i,\pi_i]$. 
Proposition \ref{prop:UP} shows that $[x_1,\pi_1]$ is proper on $X_\sim$ if and only if 
there is no point $[x,\pi] \neq [x_1,\pi_1]$ with $[x,\pi] < [x_1,\pi_1]$. 
We write $[x_1,\pi_1] \approx [x_2,\pi_2]$ if either $[x_1,\pi_1] < [x_2,\pi_2]$ 
or $[x_1,\pi_1] > [x_2,\pi_2]$. 
\begin{definition} \label{def:cluster}
A finite subset $K \subset X_\sim$ is called a {\sl cluster} 
if $[x_1,\pi_1] \in K$ 
for any $[x_1,\pi_1] < [x_2,\pi_2]$ with $[x_2,\pi_2] \in K$. 
\end{definition}
\begin{proposition} \label{prop:multi}
Let $C$ be a curve on $X$, and $C_1$, $C_2$ be curves having no irreducible 
component in common. 
Then we have the following. 
\begin{enumerate}
\item $\nu_{[x_1,\pi_1]}(C) \ge \nu_{[x_2,\pi_2]}(C)$ for 
$[x_1,\pi_1] < [x_2,\pi_2]$. 
\item $\displaystyle (C_1,C_2)=\sum_{[x,\pi] \in (C_1)_\sim \cap (C_2)_\sim}
\nu_{[x,\pi]}(C_1) \cdot \nu_{[x,\pi]}(C_2)$, 
where $(\cdot, \cdot)$ is the intersection form on $X$. 
\item $(C_1)_\sim \cap (C_2)_\sim$ is a cluster. 
\end{enumerate}
\end{proposition}
{\it Proof}. 
Assertion (1) is a consequence of Proposition \ref{prop:UP} and the fact that 
$\nu_{z_1}(C) \ge \nu_{z_2}(\pi_x^{-1}(C))$ for any blowup $\pi_x : X_x \to X$ 
of a point $x$ and for any points $z_2 \in \hat{X}$ with $z_1 =\pi(z_2) \in X$. 
In order to prove assertion (2), we use the fact that 
$(\pi_x^{-1}(C_1),\pi_x^{-1}(C_2))= (C_1,C_2)- \nu_{x}(C_1) \cdot \nu_{x}(C_2)$ 
for the blowup $\pi_x : X_x \to X$ of a point $x$. 
If $(C_1,C_2)=0$, then it follows that $\nu_x(C_1) \cdot \nu_x(C_2)=0$ for any $x \in X$ 
and thus $\nu_{[x,\pi]}(C_1) \cdot \nu_{[x,\pi]}(C_2)=0$ for any $[x,\pi] \in X_\sim$, 
which yields $(C_1)_\sim \cap (C_2)_\sim=\emptyset$. 
Hence assertion (2) holds when $(C_1,C_2)=0$. 
On the other hand, if $(C_1,C_2)>0$, then there is a point $x \in X$ such that 
$\nu_{[x,\mathrm{id}]}(C_1) \cdot \nu_{[x,\mathrm{id}]}(C_2)=\nu_x(C_1) \cdot \nu_x(C_2)>0$ and thus 
$(\pi_x^{-1}(C_1),\pi_x^{-1}(C_2))=
 (C_1,C_2)- \nu_{[x,\mathrm{id}]}(C_1) \cdot \nu_{[x,\mathrm{id}]}(C_2) < (C_1,C_2)$. 
Hence by replacing $X$ with $X_x$,  namely, $X_\sim$ with $(X_x)_\sim = X_\sim \setminus \{ [x,\mathrm{id}] \} $ by 
Proposition \ref{prop:iden}, and $C_i$ with $\pi_x^{-1}(C_i)$, we can repeat this argument finitely many times 
to yield assertion (2). 
Finally we notice that $(C_1)_\sim \cap (C_2)_\sim$ is a finite set by assertion (2). 
Assertion (1) shows that $(C_1)_\sim \cap (C_2)_\sim$ becomes a cluster, which yields assertion (3). 
\hfill $\Box$ \par\medskip
For $[z_0,\pi_0] \in X_\sim$, an element 
$[z_m,\pi_m] \in X_\sim$ satisfying 
$\pi_m=\pi_0 \circ \pi_{z_0} \circ \cdots \circ \pi_{z_{m-1}}$ and 
$z_0=\pi_{z_0} \circ \cdots \circ \pi_{z_{i-1}} (z_i)$ 
for $1 \le i \le m$ is called a {\sl point in the $m$-th infinitesimal 
neighbourhood of $[z_0,\pi_0]$} or a {\sl point infinite near to $[z_0,\pi_0]$}. 
A point in the $0$-th infinitesimal neighbourhood of $[z_0,\pi_0]$ is interpreted 
as $[z_0,\pi_0]$ itself. 
If $[z_0,\pi_0]$ is proper, $[z_m,\pi_m]$ is called an {\sl $m$-th infinitely 
near point on $X$} or an {\sl $m$-th point} for short . 
Note that if $[x_1, \pi_1] < [x_2, \pi_2]$, then $[x_2,\pi_2]$ is a point 
in the $m$-th infinitesimal neighbourhood of $[x_1,\pi_1]$ for some 
$m \in \mathbb{Z}_{\ge 0}$ by Proposition \ref{prop:UP}. 
\par
Next we construct a proper modification $\pi_K : X_K \to X$ 
from a cluster $K \subset X_{\sim}$. 
Put $K=\{ [x_0,\pi_0], [x_1,\pi_1], \dots [x_{m-1},\pi_{m-1}] \}$ 
so that if $[x_i,\pi_i] < [x_j,\pi_j]$ for $i \neq j$ then 
$[x_i,\pi_i] \neq [x_j,\pi_j]$ and $i<j$. 
We also put $Y_0:=X$ and $\nu_0:= \text{id}_{Y_0} : Y_0 \to Y_0$. 
For $k \in \{ 0,\dots,m-1\}$, let $y_k$ be a point of $Y_k$
inductively given by the relation 
\begin{equation} \label{eqn:indcpt}
(y_k, \nu_0 \circ \nu_1 \circ \cdots \circ \nu_k) \in [x_k , \pi_k], 
\end{equation}
and let $\nu_{k+1} : Y_{k+1} \to Y_k$ be the blowup of $y_k \in Y_k$. 
It should be noted that the point $y_k$ is determined uniquely. 
Indeed, when $k=0$, $\pi_0$ is locally biholomorphic at $x_0$ as $[x_0,\pi_0]$ is proper, 
and hence $y_0$ is given by $y_0=\nu_0^{-1} \circ \pi_0(x_0)=\pi_0(x_0)$. 
Moreover, under the assumption that $y_0, \dots, y_{k-1}$ are already given by (\ref{eqn:indcpt}), 
the point $[x_k,\pi_k] \in X_\sim \setminus \{ [x_0,\pi_0], \dots [x_{k-1},\pi_{k-1}]  \} \cong (Y_k)_\sim$ 
is proper on $(Y_k)_\sim$ and $y_k$ is also determined uniquely in a similar argument. 
\par
The proper modification $\pi_K : X_K \to X$ is given by the composition 
\[
\pi_K : X_K = Y_m \overset{\nu_m}{\longrightarrow} Y_{m-1} 
\overset{\nu_{m-1}}{\longrightarrow} \cdots \overset{\nu_2}{\longrightarrow} 
Y_1 \overset{\nu_1}{\longrightarrow} Y_0 =X, 
\]
called {\sl the blowup of the cluster $K$}. 
The definition of $\pi_K$ is independent of the choice of 
ordering in the cluster $K$. 
Moreover, the total transform 
$E_k=(\nu_{k+1} \circ \cdots \circ \nu_{m})^*(E_k')$ of the 
exceptional curve $E_k'$ of $\nu_k$ is called the 
{\sl exceptional divisor of $\pi_K$ over the point $[x_{k-1}, \pi_{k-1}]$}. 
\par 
Now let $X$ and $Y$ be projective surfaces and $f : X \to Y$ be a birational map. 
Assume that $Y \subset \mathbb{P}^N$ and $f(X)$ is contained in no hyperplane 
of $\mathbb{P}^N$, so that $f^* H$ is a curve on $X$ for any hyperplane 
$H \subset \mathbb{P}^N$. 
Put 
\[
I(f)_\sim:= \{ [x, \pi] \in X_\sim \, | \, 
[x,\pi] \in (f^* H)_\sim \text{ for any hyperplanes } 
H \subset \mathbb{P}^N  \}. 
\]
It follows from Proposition \ref{prop:multi} that the set $I(f)_\sim$ 
is a cluster. 
Moreover, for proper modifications $\pi : \hat{X} \to X$ and 
$\nu : \hat{Y} \to Y$, 
put $f_{\nu,\pi} := \nu^{-1} \circ f \circ \pi : \hat{X} \to \hat{Y}$. 
\begin{proposition}[\cite{A, B}] \label{prop:U2}
For any birational map $f : X \to Y$ of surfaces, put 
$\pi_0 := \pi_{I(f)_\sim} : X_0 := X_{I(f)_\sim} \to X$ and 
$\nu_0 := \pi_{I(f^{-1})_\sim} : Y_0 := Y_{I(f^{-1})_\sim} \to Y$. 
Then, the map $f_{\nu_0,\pi_0} : X_0 \to Y_0$ is a biholomorphism. 
\end{proposition}
\begin{proposition} \label{prop:newmap}
For any $[x, \pi : \hat{X} \to X] \in X_\sim \setminus I(f)_\sim$, 
there is a unique element 
$[y, \nu : \hat{Y} \to Y] \in Y_\sim \setminus I(f^{-1})_\sim$ 
such that $f_{\nu,\pi}$ is locally biholomorphic at $x \in \hat{X}$ and 
$y=f_{\nu,\pi}(x) \in \hat{Y}$ for some $($and any$)$ $(x, \pi) \in [x, \pi]$ 
and $(y, \nu) \in [y, \nu]$. 
\end{proposition}
{\it Proof}. 
First we prove the existence of 
$[y, \nu] \in Y_\sim \setminus I(f^{-1})_\sim$. 
Let $f_0 := f_{\nu_0,\pi_0} : X_0 \to Y_0$ be the biholomorphism given in 
Proposition \ref{prop:U2}, 
and let $\hat{\pi} : \hat{X} \to X_0$ be a proper modification with 
$(x_0,\pi_0 \circ \hat{\pi}) \in [x,\pi] \in X_\sim \setminus I(f)_\sim \cong (X_0)_\sim$ 
(see Proposition \ref{prop:iden}). 
Since $f_0 \circ \hat{\pi} : \hat{X} \to Y_0$ is a birational morphism, 
Proposition \ref{prop:UP} shows that there is a proper modification  
$\hat{\nu} : \hat{Y} \to Y_0$ such that 
$\hat{\nu}^{-1} \circ f_0 \circ \hat{\pi}  : \hat{X} \to \hat{Y}$ is a biholomorphism. 
As $f_{\nu,\pi_0 \circ \hat{\pi}}=\hat{\nu}^{-1} \circ f_0 \circ \hat{\pi}$ with $\nu:=\nu_0 \circ \hat{\nu}$, 
we can see that 
$[y,\nu] \in (Y_0) \cong Y_\sim \setminus I(f^{-1})_\sim$ for $y := f_{\nu,\pi_0 \circ \hat{\pi}}(x)$ 
is a desired element. 
\par 
Next we prove the uniqueness of $[y,\nu]$. 
For $i=1,2$, assume that there are $(y_i, \nu_i : \hat{Y}_{i} \to Y)$ 
such that $f_{\nu_i,\pi} : \hat{X} \to \hat{Y}_i$ are locally biholomorphic 
at $x$ and $y_i=f_{\nu_i,\pi}(x)$. 
Then $\nu_2^{-1} \circ \nu_1 = f_{\nu_2,\pi} \circ f_{\nu_1,\pi}^{-1} : 
\hat{Y}_1 \to \hat{Y}_2$ is locally biholomorphic at $y_1$ and 
$y_2=f_{\nu_2,\pi} \circ f_{\nu_1,\pi}^{-1}(y_1) = \nu_2^{-1} \circ \nu_1(y_1)$, 
which means that $(y_1, \nu_1) \sim (y_2, \nu_2)$. 
\par 
Finally, we assume that 
$(x_1,\pi_1 : \hat{X}_1 \to X) \sim (x_2,\pi_2 : \hat{X}_2 \to X)$ and 
$(y_1,\nu_1 : \hat{Y}_1 \to Y) \sim (y_2,\nu_2 : \hat{Y}_2 \to Y)$, and that 
$f_{\nu_1,\pi_1} : \hat{X}_1 \to \hat{Y}_1$ is locally biholomorphic 
at $x_1$ with $y_1=f_{\nu_1,\pi_1}(x_1)$. 
Then 
$f_{\nu_2,\pi_2} =\nu_2^{-1} \circ f \circ \pi_2 = (\nu_2^{-1} 
\circ \nu_1) \circ f_{\nu_1,\pi_1} \circ (\pi_2^{-1} \circ \pi_1)^{-1}$ 
is locally biholomorphic with $y_2=f_{\nu_2,\pi_2}(x_2)$. 
Therefore, the relation $[x,\pi] \mapsto [y,\nu]$ is well-defined, 
and the proposition is established. 
\hfill $\Box$ \par\medskip
By virtue of Proposition \ref{prop:newmap}, we can define a map 
\[
f_\sim : X_\sim \setminus I(f)_\sim 
\to Y_\sim \setminus I(f^{-1})_\sim, \quad 
f_\sim([x, \pi])=[y, \nu], 
\]
which is in fact a biholomorphism with the inverse $(f^{-1})_\sim$. 
\begin{remark} \label{rem:ed}
Let $I_X=\{[x_1,\pi_1], \dots, [x_m,\pi_m] \} \subset X_\sim$ and $I_Y= \{ [y_1,\nu_1], \dots, [y_m,\nu_m] \} \subset Y_\sim$ 
be clusters satisfying $I(f)_\sim \subset I_X$, $I(f^{-1})_\sim \subset I_Y$ and $f_\sim([x_k,\pi_k])=[y_k,\nu_k]$ 
for $[x_k,\pi_k] \notin I(f)_\sim$, 
and let $\widetilde{X} \to X$ and $\widetilde{Y} \to Y$ be the blowups of $I_X$ and $I_Y$ respectively. 
Then the blowups lift $f : X \to Y$ to a biholomorphism $\widetilde{f} : \widetilde{X} \to \widetilde{Y}$, 
and $\widetilde{f}$ sends $E_k^X$ to $E_k^Y$, where $E_k^X$ and $E_k^Y$ are 
the exceptional divisors over $[x_k,\pi_k] \notin I(f)_\sim$ and $[y_k,\nu_k] \notin I(f^{-1})_\sim$ respectively. 
\end{remark}
Hereafter, a surface $X_\sim$ is denoted simply by $X$, 
a point $[x,\pi]$ by $x$, 
a map $f_\sim$ by $f$, 
a curve $C_\sim$ by $C$, and a cluster $I(f)_\sim$ 
by $I(f)$ whenever no confusion arises. 
\begin{remark} \label{rem:propmap}
Under the above notations, we have $f(x_1) < f(x_2) \in Y$ 
for any $x_1 < x_2 \in X \setminus I(f)$. 
Hence, for an $m$-th point $x \in X \setminus I(f)$ on $X$, 
it is seen that $f(x)$ is an $(m - m_+ + m_-)$-th point on $Y$, where   
\[
\begin{array}{l}
m_+:=\max \{ \ell \ge 0 \, | \, \text{there is an } (\ell-1) \text{-th point } 
x_0 \in I(f) \text{ with } x \approx x_0 \}, \\[2mm]
m_-:=\max \{ \ell \ge 0 \, | \, \text{there is an } (\ell-1)\text{-th point }
y_0 \in I(f^{-1}) \text{ with } f(x) \approx y_0 \}. 
\end{array}
\]
In particular, for a proper point $x \in X \setminus I(f)$, 
the image $f(x)$ is also proper on $Y$ if and only if 
there is no (proper) point $y \in I(f^{-1})$ with $f(x) \approx y$, 
where $f(x) \approx y$ is equivalent to $f(x) > y$ as $I(f^{-1})$ is a cluster. 
\end{remark}
\section{Construction of Rational Surface Automorphisms} \label{sec:const} 
In this section, we develop a method for constructing a rational surface 
automorphism from a composition 
$f=f_n \circ \cdots \circ f_1 : \mathbb{P}^2 \to \mathbb{P}^2$ 
of quadratic birational maps 
$f_i : \mathbb{P}^2 \to \mathbb{P}^2$ and an orbit data $\tau$. 
If $\tau$ is 
compatible with the maps $\overline{f}=(f_1,\dots,f_n)$, $f$ lifts to 
an automorphism $F: X \to X$ through a blowup 
$\pi : X \to \mathbb{P}^2$. 
Moreover, we calculate the action 
$F^* : H^2(X;\mathbb{Z}) \to H^2(X;\mathbb{Z})$ of the automorphism $F$, 
and also calculate a Weyl group element $w_{\tau}$ realized by 
the pair $(\pi,F)$. 
\par
First we consider a smooth rational surface $X$, that is, a surface 
birationally equivalent to 
$\mathbb{P}^2$, and an automorphism $F : X \to X$ of $X$. By theorems of 
Gromov and Yomdin, the {\sl topological entropy of $F$} is given by 
$h_{\mathrm{top}}(F)= \log \lambda(F^*)$, where $\lambda(F^*)$ is the 
spectral radius of the action $F^* : H^2(X;\mathbb{Z}) \to H^2(X;\mathbb{Z})$ 
on the cohomology group. In this paper, we are interested in the case where 
$F : X \to X$ has positive entropy $h_{\mathrm{top}}(F)>0$ or, in other words, 
$\lambda(F^*)>1$. Then, the surface $X$ is characterized as follows 
(see \cite{H2, N1}). 
\begin{proposition} \label{pro:H1} 
If $X$ admits an automorphism $F : X \to X$ with $\lambda(F^*) > 1$, 
then there is a birational morphism $\pi : X \to \mathbb{P}^2$. 
\end{proposition}
It is known that any birational morphism $\pi : X \to \mathbb{P}^2$ is 
expressed as $\pi=\pi_K$ for some cluster 
$K=\{x_1,\dots,x_N \}$, where $\pi_K$ is the blowup of $K$. 
Then $\pi : X \to \mathbb{P}^2$ gives an expression of 
the cohomology group: 
$H^2(X;\mathbb{Z}) \cong \mathrm{Pic}(X) =\mathbb{Z} [H] \oplus \mathbb{Z} [E_1] 
\oplus \cdots \oplus \mathbb{Z} [E_N]$, 
where $H$ is the total transform of a line in $\mathbb{P}^2$, and $E_i$ is 
the exceptional divisor over the point $x_i$. The intersection form on 
$H^2(X;\mathbb{Z})$ is given by 
\[
\left\{
\begin{array}{ll}
([H],[H])=1 & ~ \\[2mm]
([E_i],[E_j])=-\delta_{i,j} & (i,j=1,\dots, N) \\[2mm]
([H],[E_i])=0 & (i=1,\dots, N), 
\end{array}
\right. 
\]
where $\delta_{i,j}$ is the Kronecker delta. 
Therefore $H^2(X;\mathbb{Z})$ is isometric to the Lorentz lattice 
$\mathbb{Z}^{1,N}$ via the marking 
$\phi_{\pi} : \mathbb{Z}^{1,N} \to H^2(X;\mathbb{Z})$
given in (\ref{eqn:marking}). 
The marking $\phi_{\pi}$ is isometric and determined uniquely by 
$\pi : X \to \mathbb{P}^2$ 
in the sense that if $\phi_{\pi}$ and $\phi_{\pi}'$ are 
markings determined by $\pi$, then there is an element 
$\wp \in \langle \rho_1,\dots,\rho_{N-1} \rangle$, 
acting by a permutation on the 
basis elements $(e_1,\dots,e_N)$, such that 
$\phi_{\pi}= \phi_{\pi}' \circ \wp$. 
Moreover the Weyl group $W_N$ plays an important role in our discussion 
as is mentioned in the following proposition (see \cite{DO, H1, N2}). 
\begin{proposition} \label{pro:H2} 
For any birational morphism  $\pi : X \to \mathbb{P}^2$ and any automorphism 
$F: X \to X$, there is a unique element $w \in W_N$ such that diagram 
(\ref{eqn:diag}) commutes. 
\end{proposition}
Thus, a pair $(\pi,F)$ determines $w$ uniquely, up to 
conjugacy by an element of $\langle \rho_1,\dots,\rho_{N-1} \rangle$. 
In this case, the element $w$ is said to be {\sl realized} by $(\pi,F)$,  
and the entropy of $F$ is expressed as 
$h_{\mathrm{top}}(F)= \log \lambda (w)$. 
Summing up these discussions, we have the following proposition. 
\begin{proposition} \label{pro:expent} 
The entropy of any automorphism $F: X \to X$ on a rational surface $X$ is 
given by $h_{\mathrm{top}}(F)= \log \lambda$ for some $\lambda \in \Lambda$, 
where $\Lambda$ is given in (\ref{eqn:PV}).  
\end{proposition}
Indeed, when $F : X \to X$ satisfies $\lambda(F^*)=1$, the entropy of 
$F$ is expressed as $h_{\mathrm{top}}(F)= \log \lambda(e)$ with the unit 
element $e \in W_N$. 
\begin{remark} \label{rem:zeroent}
If $\pi : X \to \mathbb{P}^2$ is a blowup of $N$ points with $N \le 9$, and 
$F : X \to X$ is an automorphism, then it follows that $h_{\mathrm{top}}(F)=0$ 
(see e.g. \cite{M}). 
\end{remark}
\par
Next we recall some properties of quadratic maps. 
Let $f : \mathbb{P}^2 \to \mathbb{P}^2$ be a quadratic birational map on 
$\mathbb{P}^2$. It is known that $f$ can be expressed as 
$f= l_- \circ g \circ l_+^{-1}$, 
where $l_+, l_- : \mathbb{P}^2 \to \mathbb{P}^2$ are linear 
transformations, and $g : \mathbb{P}^2 \to \mathbb{P}^2$ with 
$I(g^{\pm 1})=\{p_1,p_2,p_3\}$ is a 
simple quadratic birational map given in exactly one of the following 
three cases: 
\begin{description}
\item[Case 1] $g=g_1 : \mathbb{P}^2 \ni 
[x:y:z] \mapsto  [yz:zx:xy] \in \mathbb{P}^2$, \quad and \quad 
$\left\{ \begin{array}{l}
p_{1} = [1:0:0] \\
p_{2} = [0:1:0] \\
p_{3} = [0:0:1],  
\end{array} \right. $
\item[Case 2] $g=g_2 : \mathbb{P}^2 \ni [x:y:z] \mapsto 
[xz:yz:x^2] \in \mathbb{P}^2$, \quad and \quad 
$\left\{ \begin{array}{l}
p_{1} = [0:1:0] \\
p_{2} = [0:0:1] \\
p_{3} > p_1,  
\end{array} \right. $ 
\item[Case 3] $g=g_3 : \mathbb{P}^2 \ni [x:y:z] \mapsto 
[x^2:xy:y^2+xz] \in \mathbb{P}^2$, \quad and \quad
$\left\{ \begin{array}{l}
p_{1} = [0:0:1] \\
p_{2} > p_1 \\
p_{3} > p_2.  
\end{array} \right. $ 
\end{description}
Let $\pi : X \to \mathbb{P}^2$ 
be the blowup of the cluster $\{ p_{1}, p_{2}, p_{3} \}$, 
and let $H$ be the total transform of a line in $\mathbb{P}^2$, 
$L_{1}$, $L_{2}$, $L_{3}$ 
be the strict transforms of the lines ${x=0}$, ${y=0}$, ${z=0}$, respectively, 
and $E_{i}$ be the exceptional divisor over the point 
$p_{i}$ for $i=1,2,3$. Then $L_{i}$ is linearly equivalent to 
$H - E_{j} - E_{k}$ for $\{i,j,k\}=\{1,2,3\}$. 
The birational map $g$ lifts to an automorphism 
$\widetilde{g} : X \to X$, which sends irreducible 
rational curves as follows: 
\begin{description}
\item[Case 1] $\widetilde{g}=\widetilde{g}_1 
: ~E_{i} \to L_{i} \quad (i \in \{1,2,3\})$, 
\\[-2mm]
\item[Case 2] 
$\widetilde{g}=\widetilde{g}_2 : \left\{ 
\begin{array}{rcl}
E_{1}-E_{3} & \to & E_{1}-E_{3} \\
E_{i} & \to & L_{i} \quad (i\in \{2,3\}), 
\end{array} \right. $ 
\\[-2mm]
\item[Case 3] 
$\widetilde{g}=\widetilde{g}_3 : \left\{ 
\begin{array}{rcl}
E_{1}-E_{2} & \to & E_{1}-E_{2} \\
E_{2}-E_{3} & \to & E_{2}-E_{3} \\
E_{3} & \to & L_{3}, 
\end{array} \right. $ 
\end{description}
\begin{figure}[!t] 
\begin{center}
\unitlength 0.1in
\begin{picture}( 55.6100, 24.1300)(  2.0000,-26.7300)
%
\special{pn 13}%
\special{pa 200 2538}%
\special{pa 2422 2538}%
\special{fp}%
%
\special{pn 13}%
\special{pa 200 2674}%
\special{pa 1396 1712}%
\special{fp}%
%
\special{pn 13}%
\special{pa 2422 2674}%
\special{pa 1224 1712}%
\special{fp}%
%
\special{pn 20}%
\special{sh 1.000}%
\special{ar 1310 1780 36 30  0.0000000 6.2831853}%
%
\special{pn 20}%
\special{sh 1.000}%
\special{ar 2252 2538 34 28  0.0000000 6.2831853}%
%
\special{pn 20}%
\special{sh 1.000}%
\special{ar 370 2538 36 28  0.0000000 6.2831853}%
%
\special{pn 13}%
\special{pa 2592 2124}%
\special{pa 3276 2124}%
\special{fp}%
\special{sh 1}%
\special{pa 3276 2124}%
\special{pa 3210 2104}%
\special{pa 3224 2124}%
\special{pa 3210 2144}%
\special{pa 3276 2124}%
\special{fp}%
%
\special{pn 8}%
\special{pa 1312 1438}%
\special{pa 1312 1712}%
\special{fp}%
\special{sh 1}%
\special{pa 1312 1712}%
\special{pa 1332 1646}%
\special{pa 1312 1660}%
\special{pa 1292 1646}%
\special{pa 1312 1712}%
\special{fp}%
%
\special{pn 8}%
\special{pa 4644 1438}%
\special{pa 4644 1712}%
\special{fp}%
\special{sh 1}%
\special{pa 4644 1712}%
\special{pa 4664 1646}%
\special{pa 4644 1660}%
\special{pa 4624 1646}%
\special{pa 4644 1712}%
\special{fp}%
%
\special{pn 13}%
\special{pa 714 1316}%
\special{pa 1910 1316}%
\special{fp}%
%
\special{pn 13}%
\special{pa 884 1396}%
\special{pa 200 752}%
\special{da 0.070}%
%
\special{pn 13}%
\special{pa 1738 1396}%
\special{pa 2422 752}%
\special{da 0.070}%
%
\special{pn 13}%
\special{pa 200 914}%
\special{pa 884 270}%
\special{fp}%
%
\special{pn 13}%
\special{pa 2422 914}%
\special{pa 1738 270}%
\special{fp}%
%
\special{pn 13}%
\special{pa 714 352}%
\special{pa 1910 352}%
\special{da 0.070}%
%
\special{pn 13}%
\special{pa 4046 1316}%
\special{pa 5240 1316}%
\special{da 0.070}%
\put(29.0100,-21.0400){\makebox(0,0)[lb]{$g_1$}}%
%
\special{pn 13}%
\special{pa 2592 832}%
\special{pa 3276 832}%
\special{fp}%
\special{sh 1}%
\special{pa 3276 832}%
\special{pa 3210 812}%
\special{pa 3224 832}%
\special{pa 3210 852}%
\special{pa 3276 832}%
\special{fp}%
\put(29.0100,-8.0700){\makebox(0,0)[lb]{$\widetilde{g}_1$}}%
\put(14.3900,-18.5500){\makebox(0,0)[lb]{$p_{1}$}}%
\put(3.6000,-27.3000){\makebox(0,0)[lb]{$p_{2}$}}%
\put(20.7000,-27.4000){\makebox(0,0)[lb]{$p_{3}$}}%
\put(19.5000,-4.3000){\makebox(0,0)[lb]{$E_{1}$}}%
\put(18.0000,-11.0000){\makebox(0,0)[lb]{$E_{3}$}}%
\put(2.5000,-12.2000){\makebox(0,0)[lb]{$E_{2}$}}%
\put(2.0000,-22.5000){\makebox(0,0)[lb]{$\{z=0\}$}}%
\put(19.1000,-22.5000){\makebox(0,0)[lb]{$\{y=0\}$}}%
\put(11.0000,-27.4000){\makebox(0,0)[lb]{$\{x=0\}$}}%
%
\special{pn 13}%
\special{pa 3540 2538}%
\special{pa 5762 2538}%
\special{fp}%
%
\special{pn 13}%
\special{pa 3540 2674}%
\special{pa 4736 1712}%
\special{fp}%
%
\special{pn 13}%
\special{pa 5762 2674}%
\special{pa 4564 1712}%
\special{fp}%
%
\special{pn 20}%
\special{sh 1.000}%
\special{ar 4650 1780 36 30  0.0000000 6.2831853}%
%
\special{pn 20}%
\special{sh 1.000}%
\special{ar 5592 2538 34 28  0.0000000 6.2831853}%
%
\special{pn 20}%
\special{sh 1.000}%
\special{ar 3710 2538 36 28  0.0000000 6.2831853}%
%
\special{pn 8}%
\special{pa 4652 1438}%
\special{pa 4652 1712}%
\special{fp}%
\special{sh 1}%
\special{pa 4652 1712}%
\special{pa 4672 1646}%
\special{pa 4652 1660}%
\special{pa 4632 1646}%
\special{pa 4652 1712}%
\special{fp}%
%
\special{pn 13}%
\special{pa 4054 1316}%
\special{pa 5250 1316}%
\special{fp}%
%
\special{pn 13}%
\special{pa 5078 1396}%
\special{pa 5762 752}%
\special{da 0.070}%
%
\special{pn 13}%
\special{pa 3540 914}%
\special{pa 4224 270}%
\special{fp}%
%
\special{pn 13}%
\special{pa 5762 914}%
\special{pa 5078 270}%
\special{fp}%
%
\special{pn 13}%
\special{pa 4054 352}%
\special{pa 5250 352}%
\special{da 0.070}%
\put(47.7900,-18.5500){\makebox(0,0)[lb]{$p_{1}$}}%
\put(37.0000,-27.3000){\makebox(0,0)[lb]{$p_{2}$}}%
\put(54.1000,-27.4000){\makebox(0,0)[lb]{$p_{3}$}}%
\put(52.9000,-4.3000){\makebox(0,0)[lb]{$E_{1}$}}%
\put(51.4000,-11.0000){\makebox(0,0)[lb]{$E_{3}$}}%
\put(35.9000,-12.2000){\makebox(0,0)[lb]{$E_{2}$}}%
\put(35.4000,-22.5000){\makebox(0,0)[lb]{$\{z=0\}$}}%
\put(52.5000,-22.5000){\makebox(0,0)[lb]{$\{y=0\}$}}%
\put(44.4000,-27.4000){\makebox(0,0)[lb]{$\{x=0\}$}}%
%
\special{pn 13}%
\special{pa 4210 1390}%
\special{pa 3526 746}%
\special{da 0.070}%
%
\special{pn 8}%
\special{pa 3940 610}%
\special{pa 3968 626}%
\special{pa 3994 642}%
\special{pa 4020 658}%
\special{pa 4048 672}%
\special{pa 4074 688}%
\special{pa 4100 704}%
\special{pa 4126 720}%
\special{pa 4154 736}%
\special{pa 4180 754}%
\special{pa 4206 770}%
\special{pa 4234 786}%
\special{pa 4260 804}%
\special{pa 4286 822}%
\special{pa 4314 840}%
\special{pa 4340 858}%
\special{pa 4366 876}%
\special{pa 4392 896}%
\special{pa 4420 916}%
\special{pa 4446 936}%
\special{pa 4472 956}%
\special{pa 4498 978}%
\special{pa 4524 998}%
\special{pa 4552 1022}%
\special{pa 4578 1044}%
\special{pa 4604 1068}%
\special{pa 4628 1092}%
\special{pa 4650 1118}%
\special{pa 4670 1144}%
\special{pa 4686 1172}%
\special{pa 4698 1200}%
\special{pa 4708 1230}%
\special{pa 4710 1262}%
\special{pa 4710 1270}%
\special{sp}%
%
\special{pn 8}%
\special{pa 2240 350}%
\special{pa 2272 350}%
\special{pa 2306 350}%
\special{pa 2338 350}%
\special{pa 2370 350}%
\special{pa 2402 348}%
\special{pa 2434 348}%
\special{pa 2466 348}%
\special{pa 2498 348}%
\special{pa 2530 348}%
\special{pa 2564 348}%
\special{pa 2596 350}%
\special{pa 2628 350}%
\special{pa 2660 350}%
\special{pa 2692 350}%
\special{pa 2724 352}%
\special{pa 2756 352}%
\special{pa 2788 354}%
\special{pa 2820 356}%
\special{pa 2852 358}%
\special{pa 2884 360}%
\special{pa 2916 362}%
\special{pa 2948 364}%
\special{pa 2980 366}%
\special{pa 3012 370}%
\special{pa 3042 372}%
\special{pa 3074 376}%
\special{pa 3106 380}%
\special{pa 3138 384}%
\special{pa 3170 388}%
\special{pa 3200 392}%
\special{pa 3232 398}%
\special{pa 3264 404}%
\special{pa 3294 410}%
\special{pa 3326 416}%
\special{pa 3358 422}%
\special{pa 3388 430}%
\special{pa 3420 436}%
\special{pa 3450 444}%
\special{pa 3482 452}%
\special{pa 3512 460}%
\special{pa 3544 468}%
\special{pa 3574 478}%
\special{pa 3606 486}%
\special{pa 3636 494}%
\special{pa 3666 504}%
\special{pa 3698 514}%
\special{pa 3728 522}%
\special{pa 3758 532}%
\special{pa 3790 542}%
\special{pa 3820 550}%
\special{pa 3820 550}%
\special{sp}%
%
\special{pn 8}%
\special{pa 2110 1090}%
\special{pa 2144 1096}%
\special{pa 2176 1102}%
\special{pa 2210 1106}%
\special{pa 2242 1112}%
\special{pa 2276 1116}%
\special{pa 2308 1122}%
\special{pa 2342 1126}%
\special{pa 2374 1130}%
\special{pa 2408 1134}%
\special{pa 2440 1138}%
\special{pa 2472 1142}%
\special{pa 2506 1146}%
\special{pa 2538 1150}%
\special{pa 2570 1152}%
\special{pa 2604 1154}%
\special{pa 2636 1156}%
\special{pa 2668 1158}%
\special{pa 2700 1158}%
\special{pa 2732 1160}%
\special{pa 2764 1160}%
\special{pa 2796 1158}%
\special{pa 2828 1158}%
\special{pa 2858 1156}%
\special{pa 2890 1154}%
\special{pa 2922 1150}%
\special{pa 2952 1146}%
\special{pa 2984 1142}%
\special{pa 3014 1136}%
\special{pa 3046 1130}%
\special{pa 3076 1124}%
\special{pa 3106 1116}%
\special{pa 3136 1108}%
\special{pa 3166 1098}%
\special{pa 3196 1088}%
\special{pa 3226 1078}%
\special{pa 3254 1066}%
\special{pa 3284 1054}%
\special{pa 3314 1040}%
\special{pa 3342 1026}%
\special{pa 3372 1012}%
\special{pa 3400 998}%
\special{pa 3428 984}%
\special{pa 3458 968}%
\special{pa 3486 954}%
\special{pa 3510 940}%
\special{sp}%
%
\special{pn 8}%
\special{pa 600 1060}%
\special{pa 630 1046}%
\special{pa 658 1030}%
\special{pa 686 1016}%
\special{pa 716 1002}%
\special{pa 744 986}%
\special{pa 774 972}%
\special{pa 802 958}%
\special{pa 830 944}%
\special{pa 860 930}%
\special{pa 888 916}%
\special{pa 918 902}%
\special{pa 946 888}%
\special{pa 976 874}%
\special{pa 1004 860}%
\special{pa 1034 846}%
\special{pa 1064 834}%
\special{pa 1092 820}%
\special{pa 1122 808}%
\special{pa 1152 796}%
\special{pa 1180 784}%
\special{pa 1210 772}%
\special{pa 1240 760}%
\special{pa 1270 748}%
\special{pa 1300 738}%
\special{pa 1330 728}%
\special{pa 1360 718}%
\special{pa 1390 708}%
\special{pa 1420 698}%
\special{pa 1452 688}%
\special{pa 1482 680}%
\special{pa 1512 672}%
\special{pa 1544 664}%
\special{pa 1574 656}%
\special{pa 1606 650}%
\special{pa 1636 642}%
\special{pa 1668 636}%
\special{pa 1700 630}%
\special{pa 1730 624}%
\special{pa 1762 618}%
\special{pa 1794 612}%
\special{pa 1826 606}%
\special{pa 1858 600}%
\special{pa 1888 596}%
\special{pa 1920 590}%
\special{pa 1952 586}%
\special{pa 1980 580}%
\special{sp}%
%
\special{pn 13}%
\special{pa 1980 580}%
\special{pa 1980 580}%
\special{pa 1980 580}%
\special{pa 1980 580}%
\special{pa 1980 580}%
\special{pa 1980 580}%
\special{pa 1980 580}%
\special{pa 1980 580}%
\special{pa 1980 580}%
\special{pa 1980 580}%
\special{pa 1980 580}%
\special{pa 1980 580}%
\special{pa 1980 580}%
\special{pa 1980 580}%
\special{pa 1980 580}%
\special{pa 1980 580}%
\special{pa 1980 580}%
\special{pa 1980 580}%
\special{pa 1980 580}%
\special{pa 1980 580}%
\special{pa 1980 580}%
\special{pa 1980 580}%
\special{pa 1980 580}%
\special{pa 1980 580}%
\special{pa 1980 580}%
\special{pa 1980 580}%
\special{pa 1980 580}%
\special{pa 1980 580}%
\special{pa 1980 580}%
\special{pa 1980 580}%
\special{pa 1980 580}%
\special{pa 1980 580}%
\special{pa 1980 580}%
\special{pa 1980 580}%
\special{pa 1980 580}%
\special{pa 1980 580}%
\special{pa 1980 580}%
\special{pa 1980 580}%
\special{pa 1980 580}%
\special{pa 1980 580}%
\special{pa 1980 580}%
\special{pa 1980 580}%
\special{pa 1980 580}%
\special{pa 1980 580}%
\special{pa 1980 580}%
\special{pa 1980 580}%
\special{pa 1980 580}%
\special{pa 1980 580}%
\special{pa 1980 580}%
\special{pa 1980 580}%
\special{pa 1980 580}%
\special{pa 1980 580}%
\special{pa 1980 580}%
\special{sp}%
%
\special{pn 8}%
\special{pa 3980 580}%
\special{pa 5350 580}%
\special{fp}%
\special{sh 1}%
\special{pa 5350 580}%
\special{pa 5284 560}%
\special{pa 5298 580}%
\special{pa 5284 600}%
\special{pa 5350 580}%
\special{fp}%
%
\special{pn 8}%
\special{pa 3440 980}%
\special{pa 3500 930}%
\special{fp}%
\special{sh 1}%
\special{pa 3500 930}%
\special{pa 3436 958}%
\special{pa 3460 964}%
\special{pa 3462 988}%
\special{pa 3500 930}%
\special{fp}%
%
\special{pn 8}%
\special{pa 2130 580}%
\special{pa 3790 580}%
\special{fp}%
%
\special{pn 8}%
\special{pa 4700 1210}%
\special{pa 4700 1260}%
\special{fp}%
\special{sh 1}%
\special{pa 4700 1260}%
\special{pa 4720 1194}%
\special{pa 4700 1208}%
\special{pa 4680 1194}%
\special{pa 4700 1260}%
\special{fp}%
\end{picture}%
\end{center}
\caption{Geometry of $g_1 : \mathbb{P}^2 \to \mathbb{P}^2$} 
\label{fig:quadratic1} 
\begin{center}
\unitlength 0.1in
\begin{picture}( 62.4100, 25.1000)(  2.0000,-27.1000)
%
\special{pn 13}%
\special{pa 324 2134}%
\special{pa 2772 2134}%
\special{fp}%
%
\special{pn 20}%
\special{sh 1.000}%
\special{ar 2584 2134 40 30  0.0000000 6.2831853}%
%
\special{pn 20}%
\special{sh 1.000}%
\special{ar 510 2134 40 30  0.0000000 6.2831853}%
%
\special{pn 13}%
\special{pa 2960 2132}%
\special{pa 3712 2132}%
\special{fp}%
\special{sh 1}%
\special{pa 3712 2132}%
\special{pa 3646 2112}%
\special{pa 3660 2132}%
\special{pa 3646 2152}%
\special{pa 3712 2132}%
\special{fp}%
%
\special{pn 8}%
\special{pa 1546 1564}%
\special{pa 1546 1848}%
\special{fp}%
\special{sh 1}%
\special{pa 1546 1848}%
\special{pa 1566 1782}%
\special{pa 1546 1796}%
\special{pa 1526 1782}%
\special{pa 1546 1848}%
\special{fp}%
\put(32.9800,-21.1100){\makebox(0,0)[lb]{$g_2$}}%
%
\special{pn 13}%
\special{pa 2960 938}%
\special{pa 3712 938}%
\special{fp}%
\special{sh 1}%
\special{pa 3712 938}%
\special{pa 3646 918}%
\special{pa 3660 938}%
\special{pa 3646 958}%
\special{pa 3712 938}%
\special{fp}%
\put(32.9800,-9.1100){\makebox(0,0)[lb]{$\widetilde{g}_2$}}%
\put(5.6200,-23.8000){\makebox(0,0)[lb]{$p_{1}<p_{3}$}}%
\put(25.1100,-23.8000){\makebox(0,0)[lb]{$p_{2}$}}%
%
\special{pn 13}%
\special{pa 510 2700}%
\special{pa 510 1848}%
\special{fp}%
%
\special{pn 13}%
\special{pa 324 1424}%
\special{pa 2772 1424}%
\special{fp}%
%
\special{pn 13}%
\special{pa 794 1492}%
\special{pa 228 926}%
\special{da 0.070}%
%
\special{pn 13}%
\special{pa 228 1066}%
\special{pa 794 500}%
\special{da 0.070}%
%
\special{pn 13}%
\special{pa 698 712}%
\special{pa 698 200}%
\special{fp}%
%
\special{pn 13}%
\special{pa 2584 1564}%
\special{pa 2584 996}%
\special{da 0.070}%
\put(2.1000,-7.9000){\makebox(0,0)[lb]{$E_{3}$}}%
\put(24.7800,-9.8100){\makebox(0,0)[lb]{$E_{2}$}}%
%
\special{pn 8}%
\special{pa 766 500}%
\special{pa 200 1066}%
\special{fp}%
\put(5.4200,-12.7000){\makebox(0,0)[lb]{$E_{1}-E_{3}$}}%
\put(13.7000,-23.6000){\makebox(0,0)[lb]{$\{x=0\}$}}%
\put(3.1000,-18.5000){\makebox(0,0)[lb]{$\{z=0\}$}}%
%
\special{pn 13}%
\special{pa 3994 2144}%
\special{pa 6442 2144}%
\special{fp}%
%
\special{pn 20}%
\special{sh 1.000}%
\special{ar 6254 2144 40 30  0.0000000 6.2831853}%
%
\special{pn 20}%
\special{sh 1.000}%
\special{ar 4180 2144 40 30  0.0000000 6.2831853}%
%
\special{pn 8}%
\special{pa 5216 1574}%
\special{pa 5216 1858}%
\special{fp}%
\special{sh 1}%
\special{pa 5216 1858}%
\special{pa 5236 1792}%
\special{pa 5216 1806}%
\special{pa 5196 1792}%
\special{pa 5216 1858}%
\special{fp}%
\put(42.3200,-23.9000){\makebox(0,0)[lb]{$p_{1}<p_{3}$}}%
\put(61.8100,-23.9000){\makebox(0,0)[lb]{$p_{2}$}}%
%
\special{pn 13}%
\special{pa 4180 2710}%
\special{pa 4180 1858}%
\special{fp}%
%
\special{pn 13}%
\special{pa 3994 1434}%
\special{pa 6442 1434}%
\special{fp}%
%
\special{pn 13}%
\special{pa 4464 1502}%
\special{pa 3898 936}%
\special{da 0.070}%
%
\special{pn 13}%
\special{pa 3898 1076}%
\special{pa 4464 510}%
\special{da 0.070}%
%
\special{pn 13}%
\special{pa 4368 722}%
\special{pa 4368 210}%
\special{fp}%
%
\special{pn 13}%
\special{pa 6254 1574}%
\special{pa 6254 1006}%
\special{da 0.070}%
\put(38.8000,-8.0000){\makebox(0,0)[lb]{$E_{3}$}}%
%
\special{pn 8}%
\special{pa 4436 510}%
\special{pa 3870 1076}%
\special{fp}%
\put(42.1200,-12.8000){\makebox(0,0)[lb]{$E_{1}-E_{3}$}}%
\put(50.4000,-23.7000){\makebox(0,0)[lb]{$\{x=0\}$}}%
\put(39.8000,-18.6000){\makebox(0,0)[lb]{$\{z=0\}$}}%
%
\special{pn 8}%
\special{pa 1400 1200}%
\special{pa 2490 1200}%
\special{fp}%
%
\special{pn 8}%
\special{pa 2680 1200}%
\special{pa 4120 1200}%
\special{fp}%
\special{sh 1}%
\special{pa 4120 1200}%
\special{pa 4054 1180}%
\special{pa 4068 1200}%
\special{pa 4054 1220}%
\special{pa 4120 1200}%
\special{fp}%
%
\special{pn 8}%
\special{pa 2610 790}%
\special{pa 2640 780}%
\special{pa 2670 768}%
\special{pa 2700 758}%
\special{pa 2732 746}%
\special{pa 2762 734}%
\special{pa 2792 724}%
\special{pa 2822 712}%
\special{pa 2852 702}%
\special{pa 2882 692}%
\special{pa 2912 680}%
\special{pa 2942 670}%
\special{pa 2972 660}%
\special{pa 3004 650}%
\special{pa 3034 640}%
\special{pa 3064 630}%
\special{pa 3094 620}%
\special{pa 3124 610}%
\special{pa 3156 600}%
\special{pa 3186 590}%
\special{pa 3216 582}%
\special{pa 3248 574}%
\special{pa 3278 564}%
\special{pa 3308 556}%
\special{pa 3340 548}%
\special{pa 3370 540}%
\special{pa 3402 532}%
\special{pa 3432 526}%
\special{pa 3464 518}%
\special{pa 3494 512}%
\special{pa 3526 506}%
\special{pa 3558 498}%
\special{pa 3588 494}%
\special{pa 3620 488}%
\special{pa 3652 482}%
\special{pa 3684 476}%
\special{pa 3714 472}%
\special{pa 3746 466}%
\special{pa 3778 462}%
\special{pa 3810 458}%
\special{pa 3842 452}%
\special{pa 3874 448}%
\special{pa 3904 444}%
\special{pa 3936 440}%
\special{pa 3968 436}%
\special{pa 4000 434}%
\special{pa 4032 430}%
\special{pa 4064 426}%
\special{pa 4096 422}%
\special{pa 4128 420}%
\special{pa 4160 416}%
\special{pa 4192 412}%
\special{pa 4224 410}%
\special{pa 4256 406}%
\special{pa 4288 402}%
\special{pa 4310 400}%
\special{sp}%
%
\special{pn 8}%
\special{pa 540 860}%
\special{pa 572 862}%
\special{pa 604 864}%
\special{pa 636 868}%
\special{pa 668 870}%
\special{pa 700 872}%
\special{pa 732 874}%
\special{pa 764 876}%
\special{pa 796 878}%
\special{pa 828 882}%
\special{pa 860 884}%
\special{pa 892 886}%
\special{pa 924 888}%
\special{pa 956 892}%
\special{pa 988 894}%
\special{pa 1020 896}%
\special{pa 1052 900}%
\special{pa 1084 902}%
\special{pa 1116 906}%
\special{pa 1148 908}%
\special{pa 1180 912}%
\special{pa 1210 916}%
\special{pa 1242 918}%
\special{pa 1274 922}%
\special{pa 1306 926}%
\special{pa 1338 930}%
\special{pa 1370 934}%
\special{pa 1402 938}%
\special{pa 1434 942}%
\special{pa 1464 946}%
\special{pa 1496 950}%
\special{pa 1528 956}%
\special{pa 1560 960}%
\special{pa 1592 966}%
\special{pa 1622 970}%
\special{pa 1654 976}%
\special{pa 1686 982}%
\special{pa 1718 986}%
\special{pa 1748 992}%
\special{pa 1780 998}%
\special{pa 1812 1004}%
\special{pa 1842 1010}%
\special{pa 1874 1016}%
\special{pa 1906 1024}%
\special{pa 1936 1030}%
\special{pa 1968 1036}%
\special{pa 2000 1042}%
\special{pa 2030 1050}%
\special{pa 2062 1056}%
\special{pa 2094 1064}%
\special{pa 2124 1070}%
\special{pa 2156 1078}%
\special{pa 2188 1084}%
\special{pa 2218 1092}%
\special{pa 2250 1098}%
\special{pa 2280 1106}%
\special{pa 2312 1114}%
\special{pa 2344 1120}%
\special{pa 2374 1128}%
\special{pa 2406 1136}%
\special{pa 2436 1142}%
\special{pa 2468 1150}%
\special{pa 2500 1158}%
\special{pa 2510 1160}%
\special{sp}%
%
\special{pn 8}%
\special{pa 4230 410}%
\special{pa 4310 400}%
\special{fp}%
\special{sh 1}%
\special{pa 4310 400}%
\special{pa 4242 388}%
\special{pa 4258 408}%
\special{pa 4246 428}%
\special{pa 4310 400}%
\special{fp}%
%
\special{pn 8}%
\special{pa 2660 1240}%
\special{pa 2692 1250}%
\special{pa 2722 1258}%
\special{pa 2754 1266}%
\special{pa 2784 1274}%
\special{pa 2816 1282}%
\special{pa 2846 1290}%
\special{pa 2878 1300}%
\special{pa 2908 1308}%
\special{pa 2940 1314}%
\special{pa 2970 1322}%
\special{pa 3002 1330}%
\special{pa 3034 1338}%
\special{pa 3064 1344}%
\special{pa 3096 1352}%
\special{pa 3126 1358}%
\special{pa 3158 1364}%
\special{pa 3190 1370}%
\special{pa 3220 1376}%
\special{pa 3252 1382}%
\special{pa 3284 1386}%
\special{pa 3316 1392}%
\special{pa 3346 1396}%
\special{pa 3378 1400}%
\special{pa 3410 1404}%
\special{pa 3442 1406}%
\special{pa 3474 1410}%
\special{pa 3506 1412}%
\special{pa 3538 1414}%
\special{pa 3570 1416}%
\special{pa 3602 1418}%
\special{pa 3634 1420}%
\special{pa 3666 1422}%
\special{pa 3698 1424}%
\special{pa 3730 1424}%
\special{pa 3762 1426}%
\special{pa 3794 1426}%
\special{pa 3826 1428}%
\special{pa 3858 1428}%
\special{pa 3890 1430}%
\special{pa 3922 1430}%
\special{pa 3940 1430}%
\special{sp}%
%
\special{pn 8}%
\special{pa 3880 1430}%
\special{pa 3930 1430}%
\special{fp}%
\special{sh 1}%
\special{pa 3930 1430}%
\special{pa 3864 1410}%
\special{pa 3878 1430}%
\special{pa 3864 1450}%
\special{pa 3930 1430}%
\special{fp}%
\put(61.4800,-9.9100){\makebox(0,0)[lb]{$E_{2}$}}%
\end{picture}%
\end{center}
\caption{Geometry of $g_2 : \mathbb{P}^2 \to \mathbb{P}^2$} 
\label{fig:quadratic2}
\begin{center}
\unitlength 0.1in
\begin{picture}( 65.0300, 20.1000)(  2.0000,-22.3000)
%
\special{pn 13}%
\special{pa 200 2158}%
\special{pa 2814 2158}%
\special{fp}%
%
\special{pn 20}%
\special{sh 1.000}%
\special{ar 1508 2158 42 32  0.0000000 6.2831853}%
%
\special{pn 13}%
\special{pa 3014 2158}%
\special{pa 3818 2158}%
\special{fp}%
\special{sh 1}%
\special{pa 3818 2158}%
\special{pa 3752 2138}%
\special{pa 3766 2158}%
\special{pa 3752 2178}%
\special{pa 3818 2158}%
\special{fp}%
%
\special{pn 8}%
\special{pa 1508 1568}%
\special{pa 1508 1866}%
\special{fp}%
\special{sh 1}%
\special{pa 1508 1866}%
\special{pa 1528 1798}%
\special{pa 1508 1812}%
\special{pa 1488 1798}%
\special{pa 1508 1866}%
\special{fp}%
\put(33.7600,-21.3600){\makebox(0,0)[lb]{$g_3$}}%
%
\special{pn 13}%
\special{pa 3014 926}%
\special{pa 3818 926}%
\special{fp}%
\special{sh 1}%
\special{pa 3818 926}%
\special{pa 3752 906}%
\special{pa 3766 926}%
\special{pa 3752 946}%
\special{pa 3818 926}%
\special{fp}%
\put(33.7600,-8.9600){\makebox(0,0)[lb]{$\widetilde{g}_3$}}%
\put(15.0000,-24.0000){\makebox(0,0)[lb]{$p_{1}<p_{2}<p_{3}$}}%
%
\special{pn 13}%
\special{pa 200 1424}%
\special{pa 2814 1424}%
\special{fp}%
\put(22.4100,-6.0400){\makebox(0,0)[lb]{$E_{1}-E_{2}$}}%
\put(14.2600,-3.9200){\makebox(0,0)[lb]{$E_{2}-E_{3}$}}%
\put(22.5100,-10.3500){\makebox(0,0)[lb]{$E_{3}$}}%
%
\special{pn 13}%
\special{pa 1508 1500}%
\special{pa 1508 398}%
\special{da 0.070}%
%
\special{pn 13}%
\special{pa 804 544}%
\special{pa 2210 544}%
\special{da 0.070}%
%
\special{pn 8}%
\special{pa 804 984}%
\special{pa 2210 984}%
\special{fp}%
%
\special{pn 8}%
\special{pa 1538 1500}%
\special{pa 1538 398}%
\special{fp}%
%
\special{pn 13}%
\special{pa 804 962}%
\special{pa 2210 962}%
\special{da 0.070}%
%
\special{pn 13}%
\special{pa 804 1004}%
\special{pa 2210 1004}%
\special{dt 0.045}%
\put(3.6000,-23.7000){\makebox(0,0)[lb]{$\{x=0\}$}}%
%
\special{pn 13}%
\special{pa 4090 2156}%
\special{pa 6704 2156}%
\special{fp}%
%
\special{pn 20}%
\special{sh 1.000}%
\special{ar 5398 2156 42 32  0.0000000 6.2831853}%
%
\special{pn 8}%
\special{pa 5398 1566}%
\special{pa 5398 1864}%
\special{fp}%
\special{sh 1}%
\special{pa 5398 1864}%
\special{pa 5418 1796}%
\special{pa 5398 1810}%
\special{pa 5378 1796}%
\special{pa 5398 1864}%
\special{fp}%
\put(53.9000,-23.9800){\makebox(0,0)[lb]{$p_{1}<p_{2}<p_{3}$}}%
%
\special{pn 13}%
\special{pa 4090 1422}%
\special{pa 6704 1422}%
\special{fp}%
\put(61.3100,-6.0200){\makebox(0,0)[lb]{$E_{1}-E_{2}$}}%
\put(53.1600,-3.9000){\makebox(0,0)[lb]{$E_{2}-E_{3}$}}%
\put(61.4100,-10.3300){\makebox(0,0)[lb]{$E_{3}$}}%
%
\special{pn 13}%
\special{pa 5398 1498}%
\special{pa 5398 396}%
\special{da 0.070}%
%
\special{pn 13}%
\special{pa 4694 542}%
\special{pa 6100 542}%
\special{da 0.070}%
%
\special{pn 8}%
\special{pa 4694 982}%
\special{pa 6100 982}%
\special{fp}%
%
\special{pn 8}%
\special{pa 5428 1498}%
\special{pa 5428 396}%
\special{fp}%
%
\special{pn 13}%
\special{pa 4694 960}%
\special{pa 6100 960}%
\special{da 0.070}%
%
\special{pn 13}%
\special{pa 4694 1002}%
\special{pa 6100 1002}%
\special{dt 0.045}%
\put(42.5000,-23.6800){\makebox(0,0)[lb]{$\{x=0\}$}}%
%
\special{pn 8}%
\special{pa 1600 1200}%
\special{pa 5310 1200}%
\special{fp}%
\special{sh 1}%
\special{pa 5310 1200}%
\special{pa 5244 1180}%
\special{pa 5258 1200}%
\special{pa 5244 1220}%
\special{pa 5310 1200}%
\special{fp}%
%
\special{pn 8}%
\special{pa 3140 530}%
\special{pa 4610 530}%
\special{fp}%
\special{sh 1}%
\special{pa 4610 530}%
\special{pa 4544 510}%
\special{pa 4558 530}%
\special{pa 4544 550}%
\special{pa 4610 530}%
\special{fp}%
%
\special{pn 8}%
\special{pa 3510 1260}%
\special{pa 3540 1274}%
\special{pa 3570 1286}%
\special{pa 3598 1298}%
\special{pa 3628 1310}%
\special{pa 3658 1322}%
\special{pa 3688 1334}%
\special{pa 3718 1346}%
\special{pa 3748 1356}%
\special{pa 3778 1366}%
\special{pa 3810 1376}%
\special{pa 3840 1384}%
\special{pa 3872 1392}%
\special{pa 3904 1398}%
\special{pa 3934 1404}%
\special{pa 3966 1408}%
\special{pa 3998 1410}%
\special{pa 4030 1410}%
\special{pa 4030 1410}%
\special{sp}%
%
\special{pn 8}%
\special{pa 2590 960}%
\special{pa 2622 964}%
\special{pa 2654 968}%
\special{pa 2686 972}%
\special{pa 2718 976}%
\special{pa 2750 980}%
\special{pa 2780 986}%
\special{pa 2812 992}%
\special{pa 2844 998}%
\special{pa 2874 1006}%
\special{pa 2906 1014}%
\special{pa 2936 1022}%
\special{pa 2966 1030}%
\special{pa 2998 1040}%
\special{pa 3028 1050}%
\special{pa 3058 1060}%
\special{pa 3088 1070}%
\special{pa 3118 1080}%
\special{pa 3150 1092}%
\special{pa 3180 1102}%
\special{pa 3210 1114}%
\special{pa 3240 1126}%
\special{pa 3270 1138}%
\special{pa 3300 1148}%
\special{pa 3330 1160}%
\special{pa 3330 1160}%
\special{sp}%
%
\special{pn 8}%
\special{pa 3960 1410}%
\special{pa 4040 1410}%
\special{fp}%
\special{sh 1}%
\special{pa 4040 1410}%
\special{pa 3974 1390}%
\special{pa 3988 1410}%
\special{pa 3974 1430}%
\special{pa 4040 1410}%
\special{fp}%
\end{picture}%
\end{center}
\caption{Geometry of $g_3 : \mathbb{P}^2 \to \mathbb{P}^2$} 
\label{fig:quadratic3}
\end{figure} 
(see also Figures \ref{fig:quadratic1}-\ref{fig:quadratic3}). 
Note that $g$ sends a generic line to a conic passing through 
the three points $p_{1},p_{2},p_{3}$ in either case. Therefore, the 
action 
$\widetilde{g}^* : H^2(X;\mathbb{Z}) 
\to H^2(X;\mathbb{Z})$ 
on the cohomology group 
$H^2(X;\mathbb{Z}) = 
\mathbb{Z} [H] \oplus \mathbb{Z} [E_{1}] \oplus 
\mathbb{Z} [E_{2}] \oplus \mathbb{Z} [E_{3}]$ 
is given by 
\begin{equation} \label{eqn:relelem}
\widetilde{g}^* : \left\{
\begin{array}{lll}
[H] & \mapsto 2 [H] - \sum_{i=1}^3 [E_{i}] & ~ \\[2mm]
[E_{i}] & \mapsto [H] - [E_{j}] - [E_{k}] \, & 
(\{i,j,k\}=\{1,2,3\}). 
\end{array}
\right. 
\end{equation}
\par
For a general quadratic birational map 
$f=l_- \circ g \circ l_+^{-1} : \mathbb{P}^2 \to \mathbb{P}^2$, put 
\begin{equation} \label{eqn:label}
p_i^{\pm}= l_{\pm}(p_{1}), \quad 
p_j^{\pm}= l_{\pm}(p_{2}), \quad 
p_k^{\pm}= l_{\pm}(p_{3}) 
\end{equation}
with  $\{i,j,k\}=\{1,2,3\}$. 
Then the indeterminacy points of $f^{\pm 1}$ are expressed as 
\begin{equation} \label{eqn:indquad}
I(f^{\pm 1})=\{p_1^{\pm}, p_2^{\pm}, p_3^{\pm} \}. 
\end{equation}
Moreover, 
formula (\ref{eqn:relelem}) leads to the following lemma, which is stated 
in a general situation (see also Remark \ref{rem:ed}). 
\begin{lemma} \label{lem:QA}
In the above notations, assume that 
$I_{\pm}:=\{ p_1^{\pm},\dots,p_N^{\pm} \}$ are clusters satisfying 
$f(p_{i}^+)=p_{i}^-$ for any $i=4,\dots,N$. 
Let $\pi^{\pm} : X^{\pm} \to \mathbb{P}^2$ be the blowups of 
$I_{\pm}$, and let 
$H^{\pm} \subset X^{\pm}$ be the total transforms of lines in 
$\mathbb{P}^2$ under $\pi^{\pm}$, 
and $E_i^{\pm} \subset X^{\pm}$ be the exceptional divisors over 
the points $p_i^{\pm}$. Then, the quadratic birational map 
$f : \mathbb{P}^2 \to \mathbb{P}^2$ 
lifts to an isomorphism $\widetilde{f} : X^+ \to X^-$ and 
its cohomological action 
$\widetilde{f}^* : H^2(X^-;\mathbb{Z}) \to H^2(X^+;\mathbb{Z})$ is given by 
\begin{equation} \label{eqn:relquad}
\widetilde{f}^* : \left\{
\begin{array}{lll}
[H^-] & \mapsto 2 [H^+] - \sum_{i=1}^3 [E_i^+] & ~ \\[2mm]
[E_i^-] & \mapsto [H^+] - [E_{j}^+] - [E_{k}^+] \, & 
(\{i,j,k\}=\{1,2,3\}) \\[2mm]
[E_{\ell}^-] & \mapsto [E_{\ell}^+] \, & 
(\ell=4,\dots,N). 
\end{array}
\right. 
\end{equation}
\end{lemma}
\begin{remark} \label{rem:labeling}
From here on, we assume that a quadratic birational map 
$f=l_- \circ g \circ l_+^{-1} : \mathbb{P}^2 \to \mathbb{P}^2$ 
lifts to 
$\widetilde{f}: X^+ \to X^-$ whose cohomological action is given as in 
(\ref{eqn:relquad}). 
Then the points $p_{i}^{\pm}$ given in (\ref{eqn:indquad}) 
are expressed as (\ref{eqn:label}) for 
some $\{i,j,k\}=\{1,2,3\}$, and hence 
the indices of the forward indeterminacies 
are determined uniquely by those of the backward indeterminacies 
and vice versa. 
In particular, it follows that $p_i^+ < p_j^+$ if and only if 
$p_i^- < p_j^-$. 
\end{remark}
\par
Next we turn our attention to a method for constructing rational surface 
automorphisms in a general context. 
Let $Y_1, \dots,Y_n$ be smooth rational surfaces, and 
$\overline{f}:=(f_1,\dots,f_n)$ be an $n$-tuple of birational 
maps $f_{\ell} : Y_{\ell-1} \to Y_{\ell}$ with $Y_0:=Y_n$. 
Let  
$I(f_\ell)=\{ p_{\ell,1}^{+},\dots,p_{\ell,\eta_+(\ell)}^{+} \} 
\subset Y_{\ell-1}$ and 
$I(f_\ell^{-1})=\{ p_{\ell,1}^{-},\dots,p_{\ell,\eta_-(\ell)}^{-} 
\} \subset Y_{\ell}$ be the clusters of indeterminacy points, and 
let $\mathcal{K}_{\pm} := 
\{ \iota=(i,j)\,|\, i=1,\dots,n, \, j =1,\dots, \eta_{\pm}(i) \}$ 
be the sets of indices. 
Then it turns out that the cardinalities of the sets 
$\mathcal{K}_{\pm}$ are the same, that is, 
$\sum_{\ell=1}^n \eta_+(\ell)=\sum_{\ell=1}^n \eta_-(\ell)$, 
since $Y_n=Y_0$. 
Moreover, for $m \ge 0$ and 
$\iota =(i,j) \in \mathcal{K}_-$, 
we inductively put 
\[
p_{\iota}^0 := p_{\iota}^- \in Y_{i}, \qquad \qquad 
p_{\iota}^m := f_{\ell} (p_{\iota}^{m-1}) \in Y_{\ell} \quad 
(\ell \equiv i +m~(\text{mod}~n)). 
\]
Note that a point $p_{\iota}^m$ is well-defined if 
$p_{\iota}^{m-1} \notin I(f_{\ell})$. 
Moreover, let us introduce 
a {\sl generalized orbit data $\tau=(n,\sigma,\kappa)$} 
consisting of the integer $n \ge 1$, a bijection 
$\sigma : \mathcal{K}_- \to \mathcal{K}_+$ 
and a function $\kappa : \mathcal{K}_- \to \mathbb{Z}_{\ge 0}$ 
such that $\kappa(\iota) \ge 1$ provided $i_1 \le i$, or in other words, a function 
$\kappa$ satisfying $\mu(\iota) \ge 0$ for any $\iota \in \mathcal{K}_-$, 
where $\sigma(\iota)=\iota_1=(i_1,j_1)$ and $\mu : \mathcal{K}_- \to \mathbb{Z}_{\ge 0}$ is 
given by 
\begin{equation} \label{eqn:mu}
\mu(\iota)=\kappa(\iota) \cdot n + i_1- i -1 
=\theta_{i,i_1-1}(\kappa(\iota)) 
\end{equation}
with 
\begin{equation} \label{eqn:vartheta}
\theta_{i,i'}(k) := 
k \cdot n + i'-i.  \\[2mm]
\end{equation}
\begin{definition} \label{def:real}
Let $\overline{f}$ be an $n$-tuple of birational maps and 
$\tau=(n,\sigma,\kappa)$ be a generalized orbit data. 
Then $\overline{f}$ is called a {\sl realization} of $\tau$ if 
the following condition holds for any 
$\iota \in \mathcal{K}_-$: 
\begin{equation} \label{eqn:orbit1}
p_{\iota}^{m} \neq p_{\iota'}^+ \quad 
(0 \le m < \mu(\iota), \, \iota' \in 
\mathcal{K}_+), \qquad 
p_{\iota}^{\mu(\iota)}=p_{\sigma(\iota)}^+. 
\end{equation}
\end{definition}
It should be noted that in condition (\ref{eqn:orbit1}), two points $p_{\iota}^m$ and $p_{\iota'}^+$ may satisfy 
$p_{\iota}^m \approx p_{\iota'}^+$. 
From a realization $\overline{f}$ of $\tau$, we will 
construct an automorphism. 
So let us give the following lemma. 
\begin{lemma} \label{lem:selpro}
There is an element $\iota^o \in \mathcal{K}_-$ such that 
$p_{\iota^o}^{m}$ is proper for any $0 \le m \le \mu(\iota^o)$. 
\end{lemma}
{\it Proof}. 
Take an element $\iota^o=(i^o,j^o) \in \mathcal{K}_-$ 
such that $p_{\iota^o}^-$ is proper and 
\[
\mu(\iota^o) = \min \{ \mu(\iota) \, | \, 
\iota \in \mathcal{K}_- \text{ and } 
p_{\iota}^- \text{ is proper} \}. 
\]
Then, we claim that $p_{\iota^o}^{m}$ is proper for any 
$0 < m \le \mu(\iota^o)$.
Indeed, assume the contrary that $p_{\iota^o}^{m-1}$ is proper 
but $p_{\iota^o}^{m}$ is not proper for some $0 < m \le \mu(\iota^o)$. 
Then it follows from Remark \ref{rem:propmap} 
that there is a proper point $p_{\iota}^-$ such that 
$p_{\iota}^- < p_{\iota^o}^{m}$. 
The minimality of $\mu(\iota^o)$ yields 
$p_{\iota}^{k} \notin I(f_{\ell})$ for any 
$0 \le k < \mu(\iota^o) - m$, 
and so $p_{\iota}^{\mu(\iota^o) - m} < p_{\iota^o}^{\mu(\iota^o)}$ 
by Remark \ref{rem:propmap}. 
Since 
$p_{\iota^o}^{\mu(\iota^o)}=p_{\sigma(\iota^o)}^+ \in I(f_{i_1^o})$ 
and $I(f_{i_1^o})$ is a cluster, 
$p_{\iota}^{\mu(\iota^o) - m}$ 
is also an element of $I(f_{i_1^o})$ and 
thus is equal to $p_{\sigma(\iota)}^+$. This means that 
$\mu(\iota)=\mu(\iota^o)-m < \mu(\iota^o)$, 
which contradicts the assumption that 
$\mu(\iota^o)$ is minimal. 
Thus, $p_{\iota^o}^{m}$ is proper for any 
$0 \le m \le \mu(\iota^o)$. 
\hfill $\Box$ \par\medskip 
\par
\begin{figure}[t]
\begin{center}
\unitlength 0.1in
\begin{picture}( 48.5000, 25.7000)(  3.5000,-26.0000)
%
\special{pn 20}%
\special{sh 1.000}%
\special{ar 1200 2000 46 46  0.0000000 6.2831853}%
%
\special{pn 20}%
\special{sh 1.000}%
\special{ar 2000 2000 46 46  0.0000000 6.2831853}%
%
\special{pn 20}%
\special{sh 1.000}%
\special{ar 4400 2000 46 46  0.0000000 6.2831853}%
%
\special{pn 13}%
\special{pa 400 2600}%
\special{pa 396 2568}%
\special{pa 392 2536}%
\special{pa 388 2504}%
\special{pa 384 2472}%
\special{pa 382 2442}%
\special{pa 378 2410}%
\special{pa 374 2378}%
\special{pa 370 2346}%
\special{pa 368 2314}%
\special{pa 364 2282}%
\special{pa 362 2250}%
\special{pa 360 2218}%
\special{pa 356 2186}%
\special{pa 354 2154}%
\special{pa 354 2122}%
\special{pa 352 2090}%
\special{pa 352 2058}%
\special{pa 350 2026}%
\special{pa 350 1994}%
\special{pa 350 1962}%
\special{pa 352 1930}%
\special{pa 352 1898}%
\special{pa 354 1868}%
\special{pa 356 1836}%
\special{pa 358 1804}%
\special{pa 360 1772}%
\special{pa 362 1740}%
\special{pa 366 1708}%
\special{pa 368 1676}%
\special{pa 372 1644}%
\special{pa 376 1612}%
\special{pa 378 1580}%
\special{pa 382 1548}%
\special{pa 386 1516}%
\special{pa 390 1484}%
\special{pa 394 1452}%
\special{pa 398 1420}%
\special{pa 400 1400}%
\special{sp}%
%
\special{pn 13}%
\special{pa 5200 2600}%
\special{pa 5196 2568}%
\special{pa 5192 2536}%
\special{pa 5188 2504}%
\special{pa 5184 2472}%
\special{pa 5182 2442}%
\special{pa 5178 2410}%
\special{pa 5174 2378}%
\special{pa 5170 2346}%
\special{pa 5168 2314}%
\special{pa 5164 2282}%
\special{pa 5162 2250}%
\special{pa 5160 2218}%
\special{pa 5156 2186}%
\special{pa 5154 2154}%
\special{pa 5154 2122}%
\special{pa 5152 2090}%
\special{pa 5152 2058}%
\special{pa 5150 2026}%
\special{pa 5150 1994}%
\special{pa 5150 1962}%
\special{pa 5152 1930}%
\special{pa 5152 1898}%
\special{pa 5154 1868}%
\special{pa 5156 1836}%
\special{pa 5158 1804}%
\special{pa 5160 1772}%
\special{pa 5162 1740}%
\special{pa 5166 1708}%
\special{pa 5168 1676}%
\special{pa 5172 1644}%
\special{pa 5176 1612}%
\special{pa 5178 1580}%
\special{pa 5182 1548}%
\special{pa 5186 1516}%
\special{pa 5190 1484}%
\special{pa 5194 1452}%
\special{pa 5198 1420}%
\special{pa 5200 1400}%
\special{sp}%
%
\special{pn 8}%
\special{pa 1310 1958}%
\special{pa 1342 1948}%
\special{pa 1374 1940}%
\special{pa 1404 1932}%
\special{pa 1436 1924}%
\special{pa 1468 1918}%
\special{pa 1498 1912}%
\special{pa 1530 1906}%
\special{pa 1562 1902}%
\special{pa 1592 1900}%
\special{pa 1624 1900}%
\special{pa 1656 1902}%
\special{pa 1688 1904}%
\special{pa 1720 1908}%
\special{pa 1752 1912}%
\special{pa 1784 1918}%
\special{pa 1814 1924}%
\special{pa 1846 1930}%
\special{pa 1878 1938}%
\special{pa 1882 1938}%
\special{sp}%
%
\special{pn 8}%
\special{pa 1868 1938}%
\special{pa 1900 1960}%
\special{fp}%
\special{sh 1}%
\special{pa 1900 1960}%
\special{pa 1856 1906}%
\special{pa 1856 1930}%
\special{pa 1834 1940}%
\special{pa 1900 1960}%
\special{fp}%
%
\special{pn 13}%
\special{pa 2860 1970}%
\special{pa 3580 1980}%
\special{da 0.070}%
%
\special{pn 8}%
\special{pa 2090 1958}%
\special{pa 2122 1948}%
\special{pa 2154 1940}%
\special{pa 2184 1932}%
\special{pa 2216 1924}%
\special{pa 2248 1918}%
\special{pa 2278 1912}%
\special{pa 2310 1906}%
\special{pa 2342 1902}%
\special{pa 2372 1900}%
\special{pa 2404 1900}%
\special{pa 2436 1902}%
\special{pa 2468 1904}%
\special{pa 2500 1908}%
\special{pa 2532 1912}%
\special{pa 2564 1918}%
\special{pa 2594 1924}%
\special{pa 2626 1930}%
\special{pa 2658 1938}%
\special{pa 2662 1938}%
\special{sp}%
%
\special{pn 8}%
\special{pa 2648 1938}%
\special{pa 2680 1960}%
\special{fp}%
\special{sh 1}%
\special{pa 2680 1960}%
\special{pa 2636 1906}%
\special{pa 2636 1930}%
\special{pa 2614 1940}%
\special{pa 2680 1960}%
\special{fp}%
%
\special{pn 8}%
\special{pa 3730 1948}%
\special{pa 3762 1938}%
\special{pa 3794 1930}%
\special{pa 3824 1922}%
\special{pa 3856 1914}%
\special{pa 3888 1908}%
\special{pa 3918 1902}%
\special{pa 3950 1896}%
\special{pa 3982 1892}%
\special{pa 4012 1890}%
\special{pa 4044 1890}%
\special{pa 4076 1892}%
\special{pa 4108 1894}%
\special{pa 4140 1898}%
\special{pa 4172 1902}%
\special{pa 4204 1908}%
\special{pa 4234 1914}%
\special{pa 4266 1920}%
\special{pa 4298 1928}%
\special{pa 4302 1928}%
\special{sp}%
%
\special{pn 8}%
\special{pa 4288 1928}%
\special{pa 4320 1950}%
\special{fp}%
\special{sh 1}%
\special{pa 4320 1950}%
\special{pa 4276 1896}%
\special{pa 4276 1920}%
\special{pa 4254 1930}%
\special{pa 4320 1950}%
\special{fp}%
%
\special{pn 8}%
\special{pa 4470 1948}%
\special{pa 4502 1938}%
\special{pa 4534 1930}%
\special{pa 4564 1922}%
\special{pa 4596 1914}%
\special{pa 4628 1908}%
\special{pa 4658 1902}%
\special{pa 4690 1896}%
\special{pa 4722 1892}%
\special{pa 4752 1890}%
\special{pa 4784 1890}%
\special{pa 4816 1892}%
\special{pa 4848 1894}%
\special{pa 4880 1898}%
\special{pa 4912 1902}%
\special{pa 4944 1908}%
\special{pa 4974 1914}%
\special{pa 5006 1920}%
\special{pa 5038 1928}%
\special{pa 5042 1928}%
\special{sp}%
%
\special{pn 8}%
\special{pa 5028 1928}%
\special{pa 5060 1950}%
\special{fp}%
\special{sh 1}%
\special{pa 5060 1950}%
\special{pa 5016 1896}%
\special{pa 5016 1920}%
\special{pa 4994 1930}%
\special{pa 5060 1950}%
\special{fp}%
%
\special{pn 8}%
\special{pa 510 1948}%
\special{pa 542 1938}%
\special{pa 574 1930}%
\special{pa 604 1922}%
\special{pa 636 1914}%
\special{pa 668 1908}%
\special{pa 698 1902}%
\special{pa 730 1896}%
\special{pa 762 1892}%
\special{pa 792 1890}%
\special{pa 824 1890}%
\special{pa 856 1892}%
\special{pa 888 1894}%
\special{pa 920 1898}%
\special{pa 952 1902}%
\special{pa 984 1908}%
\special{pa 1014 1914}%
\special{pa 1046 1920}%
\special{pa 1078 1928}%
\special{pa 1082 1928}%
\special{sp}%
%
\special{pn 8}%
\special{pa 1068 1928}%
\special{pa 1100 1950}%
\special{fp}%
\special{sh 1}%
\special{pa 1100 1950}%
\special{pa 1056 1896}%
\special{pa 1056 1920}%
\special{pa 1034 1930}%
\special{pa 1100 1950}%
\special{fp}%
\put(7.5000,-18.5000){\makebox(0,0)[lb]{$f_{i}$}}%
\put(15.1000,-18.6000){\makebox(0,0)[lb]{$f_{i+1}$}}%
\put(22.9000,-18.6000){\makebox(0,0)[lb]{$f_{i+2}$}}%
\put(38.5000,-18.5000){\makebox(0,0)[lb]{$f_{i_1-1}$}}%
\put(46.5000,-18.5000){\makebox(0,0)[lb]{$f_{i_1}$}}%
\put(13.0000,-22.0000){\makebox(0,0)[lb]{$p_{\iota}^0$}}%
\put(21.0000,-22.0000){\makebox(0,0)[lb]{$p_{\iota}^1$}}%
\put(40.0000,-22.4000){\makebox(0,0)[lb]{$p_{\iota}^{\mu(\iota)}=p_{\sigma(\iota)}^+$}}%
%
\special{pn 20}%
\special{sh 1.000}%
\special{ar 1200 700 46 46  0.0000000 6.2831853}%
%
\special{pn 20}%
\special{sh 1.000}%
\special{ar 2000 700 46 46  0.0000000 6.2831853}%
%
\special{pn 20}%
\special{sh 1.000}%
\special{ar 4400 700 46 46  0.0000000 6.2831853}%
%
\special{pn 13}%
\special{pa 400 1300}%
\special{pa 396 1268}%
\special{pa 392 1236}%
\special{pa 388 1204}%
\special{pa 384 1172}%
\special{pa 382 1142}%
\special{pa 378 1110}%
\special{pa 374 1078}%
\special{pa 370 1046}%
\special{pa 368 1014}%
\special{pa 364 982}%
\special{pa 362 950}%
\special{pa 360 918}%
\special{pa 356 886}%
\special{pa 354 854}%
\special{pa 354 822}%
\special{pa 352 790}%
\special{pa 352 758}%
\special{pa 350 726}%
\special{pa 350 694}%
\special{pa 350 662}%
\special{pa 352 630}%
\special{pa 352 598}%
\special{pa 354 568}%
\special{pa 356 536}%
\special{pa 358 504}%
\special{pa 360 472}%
\special{pa 362 440}%
\special{pa 366 408}%
\special{pa 368 376}%
\special{pa 372 344}%
\special{pa 376 312}%
\special{pa 378 280}%
\special{pa 382 248}%
\special{pa 386 216}%
\special{pa 390 184}%
\special{pa 394 152}%
\special{pa 398 120}%
\special{pa 400 100}%
\special{sp}%
%
\special{pn 13}%
\special{pa 5200 1300}%
\special{pa 5196 1268}%
\special{pa 5192 1236}%
\special{pa 5188 1204}%
\special{pa 5184 1172}%
\special{pa 5182 1142}%
\special{pa 5178 1110}%
\special{pa 5174 1078}%
\special{pa 5170 1046}%
\special{pa 5168 1014}%
\special{pa 5164 982}%
\special{pa 5162 950}%
\special{pa 5160 918}%
\special{pa 5156 886}%
\special{pa 5154 854}%
\special{pa 5154 822}%
\special{pa 5152 790}%
\special{pa 5152 758}%
\special{pa 5150 726}%
\special{pa 5150 694}%
\special{pa 5150 662}%
\special{pa 5152 630}%
\special{pa 5152 598}%
\special{pa 5154 568}%
\special{pa 5156 536}%
\special{pa 5158 504}%
\special{pa 5160 472}%
\special{pa 5162 440}%
\special{pa 5166 408}%
\special{pa 5168 376}%
\special{pa 5172 344}%
\special{pa 5176 312}%
\special{pa 5178 280}%
\special{pa 5182 248}%
\special{pa 5186 216}%
\special{pa 5190 184}%
\special{pa 5194 152}%
\special{pa 5198 120}%
\special{pa 5200 100}%
\special{sp}%
%
\special{pn 8}%
\special{pa 1310 658}%
\special{pa 1342 648}%
\special{pa 1374 640}%
\special{pa 1404 632}%
\special{pa 1436 624}%
\special{pa 1468 618}%
\special{pa 1498 612}%
\special{pa 1530 606}%
\special{pa 1562 602}%
\special{pa 1592 600}%
\special{pa 1624 600}%
\special{pa 1656 602}%
\special{pa 1688 604}%
\special{pa 1720 608}%
\special{pa 1752 612}%
\special{pa 1784 618}%
\special{pa 1814 624}%
\special{pa 1846 630}%
\special{pa 1878 638}%
\special{pa 1882 638}%
\special{sp}%
%
\special{pn 8}%
\special{pa 1868 638}%
\special{pa 1900 660}%
\special{fp}%
\special{sh 1}%
\special{pa 1900 660}%
\special{pa 1856 606}%
\special{pa 1856 630}%
\special{pa 1834 640}%
\special{pa 1900 660}%
\special{fp}%
%
\special{pn 13}%
\special{pa 2860 670}%
\special{pa 3580 680}%
\special{da 0.070}%
%
\special{pn 8}%
\special{pa 2090 658}%
\special{pa 2122 648}%
\special{pa 2154 640}%
\special{pa 2184 632}%
\special{pa 2216 624}%
\special{pa 2248 618}%
\special{pa 2278 612}%
\special{pa 2310 606}%
\special{pa 2342 602}%
\special{pa 2372 600}%
\special{pa 2404 600}%
\special{pa 2436 602}%
\special{pa 2468 604}%
\special{pa 2500 608}%
\special{pa 2532 612}%
\special{pa 2564 618}%
\special{pa 2594 624}%
\special{pa 2626 630}%
\special{pa 2658 638}%
\special{pa 2662 638}%
\special{sp}%
%
\special{pn 8}%
\special{pa 2648 638}%
\special{pa 2680 660}%
\special{fp}%
\special{sh 1}%
\special{pa 2680 660}%
\special{pa 2636 606}%
\special{pa 2636 630}%
\special{pa 2614 640}%
\special{pa 2680 660}%
\special{fp}%
%
\special{pn 8}%
\special{pa 3730 648}%
\special{pa 3762 638}%
\special{pa 3794 630}%
\special{pa 3824 622}%
\special{pa 3856 614}%
\special{pa 3888 608}%
\special{pa 3918 602}%
\special{pa 3950 596}%
\special{pa 3982 592}%
\special{pa 4012 590}%
\special{pa 4044 590}%
\special{pa 4076 592}%
\special{pa 4108 594}%
\special{pa 4140 598}%
\special{pa 4172 602}%
\special{pa 4204 608}%
\special{pa 4234 614}%
\special{pa 4266 620}%
\special{pa 4298 628}%
\special{pa 4302 628}%
\special{sp}%
%
\special{pn 8}%
\special{pa 4288 628}%
\special{pa 4320 650}%
\special{fp}%
\special{sh 1}%
\special{pa 4320 650}%
\special{pa 4276 596}%
\special{pa 4276 620}%
\special{pa 4254 630}%
\special{pa 4320 650}%
\special{fp}%
%
\special{pn 8}%
\special{pa 4470 648}%
\special{pa 4502 638}%
\special{pa 4534 630}%
\special{pa 4564 622}%
\special{pa 4596 614}%
\special{pa 4628 608}%
\special{pa 4658 602}%
\special{pa 4690 596}%
\special{pa 4722 592}%
\special{pa 4752 590}%
\special{pa 4784 590}%
\special{pa 4816 592}%
\special{pa 4848 594}%
\special{pa 4880 598}%
\special{pa 4912 602}%
\special{pa 4944 608}%
\special{pa 4974 614}%
\special{pa 5006 620}%
\special{pa 5038 628}%
\special{pa 5042 628}%
\special{sp}%
%
\special{pn 8}%
\special{pa 5028 628}%
\special{pa 5060 650}%
\special{fp}%
\special{sh 1}%
\special{pa 5060 650}%
\special{pa 5016 596}%
\special{pa 5016 620}%
\special{pa 4994 630}%
\special{pa 5060 650}%
\special{fp}%
%
\special{pn 8}%
\special{pa 510 648}%
\special{pa 542 638}%
\special{pa 574 630}%
\special{pa 604 622}%
\special{pa 636 614}%
\special{pa 668 608}%
\special{pa 698 602}%
\special{pa 730 596}%
\special{pa 762 592}%
\special{pa 792 590}%
\special{pa 824 590}%
\special{pa 856 592}%
\special{pa 888 594}%
\special{pa 920 598}%
\special{pa 952 602}%
\special{pa 984 608}%
\special{pa 1014 614}%
\special{pa 1046 620}%
\special{pa 1078 628}%
\special{pa 1082 628}%
\special{sp}%
%
\special{pn 8}%
\special{pa 1068 628}%
\special{pa 1100 650}%
\special{fp}%
\special{sh 1}%
\special{pa 1100 650}%
\special{pa 1056 596}%
\special{pa 1056 620}%
\special{pa 1034 630}%
\special{pa 1100 650}%
\special{fp}%
\put(7.5000,-5.5000){\makebox(0,0)[lb]{$f_{i}'$}}%
\put(15.1000,-5.6000){\makebox(0,0)[lb]{$f_{i+1}'$}}%
\put(22.7000,-5.5000){\makebox(0,0)[lb]{$f_{i+2}'$}}%
\put(38.5000,-5.5000){\makebox(0,0)[lb]{$f_{i_1-1}'$}}%
\put(46.5000,-5.5000){\makebox(0,0)[lb]{$f_{i_1}'$}}%
%
\special{pn 20}%
\special{pa 2600 1000}%
\special{pa 2600 1600}%
\special{fp}%
\special{sh 1}%
\special{pa 2600 1600}%
\special{pa 2620 1534}%
\special{pa 2600 1548}%
\special{pa 2580 1534}%
\special{pa 2600 1600}%
\special{fp}%
%
\special{pn 20}%
\special{pa 4450 1300}%
\special{pa 4446 1268}%
\special{pa 4442 1236}%
\special{pa 4438 1204}%
\special{pa 4434 1172}%
\special{pa 4432 1142}%
\special{pa 4428 1110}%
\special{pa 4424 1078}%
\special{pa 4420 1046}%
\special{pa 4418 1014}%
\special{pa 4414 982}%
\special{pa 4412 950}%
\special{pa 4410 918}%
\special{pa 4406 886}%
\special{pa 4404 854}%
\special{pa 4404 822}%
\special{pa 4402 790}%
\special{pa 4402 758}%
\special{pa 4400 726}%
\special{pa 4400 694}%
\special{pa 4400 662}%
\special{pa 4402 630}%
\special{pa 4402 598}%
\special{pa 4404 568}%
\special{pa 4406 536}%
\special{pa 4408 504}%
\special{pa 4410 472}%
\special{pa 4412 440}%
\special{pa 4416 408}%
\special{pa 4418 376}%
\special{pa 4422 344}%
\special{pa 4426 312}%
\special{pa 4428 280}%
\special{pa 4432 248}%
\special{pa 4436 216}%
\special{pa 4440 184}%
\special{pa 4444 152}%
\special{pa 4448 120}%
\special{pa 4450 100}%
\special{sp}%
%
\special{pn 20}%
\special{pa 2050 1300}%
\special{pa 2046 1268}%
\special{pa 2042 1236}%
\special{pa 2038 1204}%
\special{pa 2034 1172}%
\special{pa 2032 1142}%
\special{pa 2028 1110}%
\special{pa 2024 1078}%
\special{pa 2020 1046}%
\special{pa 2018 1014}%
\special{pa 2014 982}%
\special{pa 2012 950}%
\special{pa 2010 918}%
\special{pa 2006 886}%
\special{pa 2004 854}%
\special{pa 2004 822}%
\special{pa 2002 790}%
\special{pa 2002 758}%
\special{pa 2000 726}%
\special{pa 2000 694}%
\special{pa 2000 662}%
\special{pa 2002 630}%
\special{pa 2002 598}%
\special{pa 2004 568}%
\special{pa 2006 536}%
\special{pa 2008 504}%
\special{pa 2010 472}%
\special{pa 2012 440}%
\special{pa 2016 408}%
\special{pa 2018 376}%
\special{pa 2022 344}%
\special{pa 2026 312}%
\special{pa 2028 280}%
\special{pa 2032 248}%
\special{pa 2036 216}%
\special{pa 2040 184}%
\special{pa 2044 152}%
\special{pa 2048 120}%
\special{pa 2050 100}%
\special{sp}%
%
\special{pn 20}%
\special{pa 1250 1300}%
\special{pa 1246 1268}%
\special{pa 1242 1236}%
\special{pa 1238 1204}%
\special{pa 1234 1172}%
\special{pa 1232 1142}%
\special{pa 1228 1110}%
\special{pa 1224 1078}%
\special{pa 1220 1046}%
\special{pa 1218 1014}%
\special{pa 1214 982}%
\special{pa 1212 950}%
\special{pa 1210 918}%
\special{pa 1206 886}%
\special{pa 1204 854}%
\special{pa 1204 822}%
\special{pa 1202 790}%
\special{pa 1202 758}%
\special{pa 1200 726}%
\special{pa 1200 694}%
\special{pa 1200 662}%
\special{pa 1202 630}%
\special{pa 1202 598}%
\special{pa 1204 568}%
\special{pa 1206 536}%
\special{pa 1208 504}%
\special{pa 1210 472}%
\special{pa 1212 440}%
\special{pa 1216 408}%
\special{pa 1218 376}%
\special{pa 1222 344}%
\special{pa 1226 312}%
\special{pa 1228 280}%
\special{pa 1232 248}%
\special{pa 1236 216}%
\special{pa 1240 184}%
\special{pa 1244 152}%
\special{pa 1248 120}%
\special{pa 1250 100}%
\special{sp}%
\put(13.0000,-2.0000){\makebox(0,0)[lb]{$E_{p_{\iota}^0}$}}%
\put(21.0000,-2.0000){\makebox(0,0)[lb]{$E_{p_{\iota}^1}$}}%
\put(45.0000,-2.0000){\makebox(0,0)[lb]{$E_{p_{\iota}^{\mu(\iota)}}$}}%
\put(26.6000,-13.4000){\makebox(0,0)[lb]{blowup}}%
\end{picture}%
\end{center}
\caption{Blowup of indeterminacy points} 
\label{fig:blowup}
\end{figure}
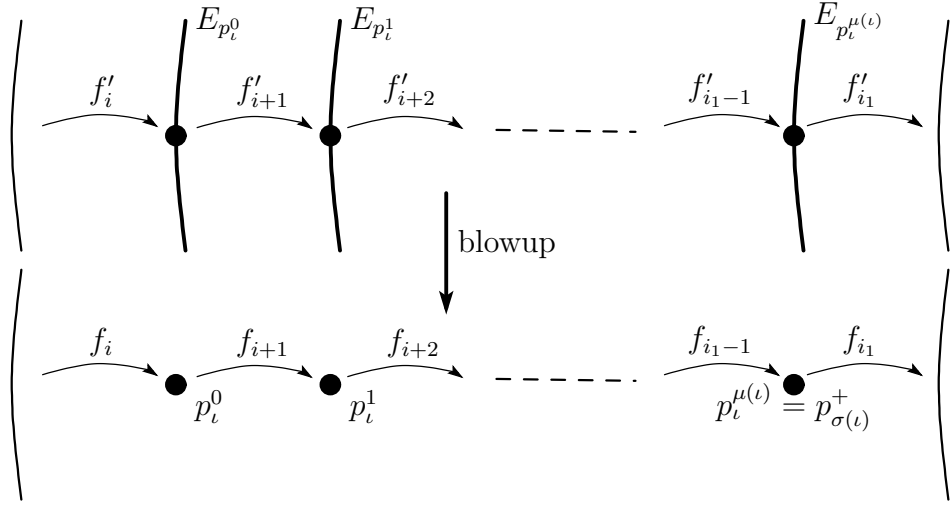

For $\iota^o=(i^o,j^o) \in \mathcal{K}_-$ given as in Lemma \ref{lem:selpro}, 
let $Y_{\ell}' \to Y_{\ell}$ be the blowups of distinct proper points 
$\{ p_{\iota^o}^m \, | \, 0 \le m \le \mu(\iota^o), \, 
i^o +m \equiv \ell ~(\text{mod}~n) \}$. 
These blowups lift $f_{\ell} : Y_{\ell-1} \to Y_{\ell}$ to 
$f_{\ell}' : Y_{\ell-1}' \to Y_{\ell}'$ (see Figure \ref{fig:blowup}). 
In this case, one has 
\begin{equation} \label{eqn:elim}
I(f_{\ell}')=
\left\{
\begin{array}{ll}
I(f_{i_1^o}) \setminus \{ p_{\iota_1^o}^+ \} & 
(\ell=i_1^o) \\[2mm]
I(f_{\ell}) & (\ell \neq i_1^o), 
\end{array}
\right. \quad
I((f_{\ell}')^{-1})=
\left\{
\begin{array}{ll}
I(f_{i^o}^{-1}) \setminus \{ p_{\iota^o}^- \} & (\ell = i^o) \\[2mm]
I(f_{\ell}^{-1}) & (\ell \neq i^o). 
\end{array}
\right.  
\end{equation}
We notice that $\overline{f}'=(f_1',\dots,f_n')$ also satisfies 
condition (\ref{eqn:orbit1}) for any 
$\iota \in \mathcal{K}_-' := \mathcal{K}_- \setminus \{ \iota^o \}$, 
which means that $\overline{f}'$ is a realization of 
$\tau'=(n,\sigma|_{\mathcal{K}_-'},
\mu|_{\mathcal{K}_-'})$. 
One can therefore repeat the above argument 
by replacing $\overline{f}$ with $\overline{f}'$, $\mathcal{K}_-$ with 
$\mathcal{K}_-'$ and 
$\tau$ with $\tau'$. In the end, from (\ref{eqn:elim}), 
the resulting map becomes a 
biholomorphism. 
Namely, let $\pi_{\ell} : X_{\ell} \to Y_{\ell}$ be the 
composition of the above blowups, that is, the blowup of the cluster 
$I_{\ell}:= \{ p_{\iota}^m \, | \, \iota \in \mathcal{K}_-, 
\, 0 \le m \le \mu(\iota), \, i +m \equiv \ell ~(\text{mod}~n) \}$. 
Then the blowups $\pi_{\ell}$ lift $f_{\ell} : Y_{\ell-1} \to Y_{\ell}$ to 
biholomorphisms $F_{\ell}: X_{\ell-1} \to X_{\ell}$: 
\[
\begin{CD}
X_{\ell-1} @> F_{\ell} >> X_{\ell} \\
@V \pi_{\ell-1} VV @VV \pi_{\ell} V \\
Y_{\ell-1} @> f_{\ell} >> Y_{\ell}, 
\end{CD}
\]
and $\pi_{\tau}:=\pi_{n} : X_{\tau} \to Y$ also lifts 
$f:=f_n \circ \cdots \circ f_1 : Y \to Y$ to the automorphism 
$F_{\tau}:=F_n \circ \cdots \circ F_1 : X_{\tau} \to X_{\tau}$, 
where $Y:=Y_0=Y_n$ and $X_{\tau}:=X_0=X_n$. 
\par
We now restrict our attention to the case 
where each component of $\overline{f}=(f_1,\dots,f_n)$ is a 
quadratic birational map with $Y_{\ell}=\mathbb{P}_{\ell}^2$ and 
$\tau$ is an original orbit data with 
$\mathcal{K}_+=\mathcal{K}_-=\mathcal{K}(n)$. 
Note that $f_{\ell}$ satisfies the assumption in Lemma \ref{lem:QA} 
with $I_+=I_{\ell-1}$ and $I_-=I_{\ell}$. 
Therefore, the cohomological action 
of the biholomorphism $F_{\ell} : X_{\ell-1} \to X_{\ell}$ 
is expressed as in the form of (\ref{eqn:relquad}), 
and that of the automorphism $F_{\tau} : X_{\tau} \to X_{\tau}$ can be 
calculated in terms of the composition $F_{\tau}^*=F_1^* \circ \cdots \circ F_n^*$. 
The Weyl group element $w_{\tau} \in W_N$ realized by 
$(\pi_{\tau}, F_{\tau})$ is given as follows. 
\par
Let $H \subset X_{\tau}$ be the total transform of a line in 
$\mathbb{P}^2$ under $\pi_{\tau}$ and 
$E_{\iota}^{k} \subset X_{\tau}$ be the exceptional divisor over 
$p_{\iota}^{m}$, where $m=\theta_{i,0}(k)$. 
Then the cohomology group of $X_{\tau}$ is expressed as 
$H^2(X_{\tau};\mathbb{Z})= \mathbb{Z} [H] \oplus 
\bigr( \oplus_{\iota \in \mathcal{K}(n)} 
\oplus_{k=1}^{\kappa(\iota)} \mathbb{Z} [E_{\iota}^{k}]\bigl)$. 
Now, we consider the lattice 
\[
\mathbb{Z}^{\tau}:= \mathbb{Z} e_0 \oplus \bigr( \oplus_{\iota \in 
\mathcal{K}(n)} \oplus_{k=1}^{\kappa(\iota)} \mathbb{Z} 
e_{\iota}^{k} \bigl) \cong \mathbb{Z}^{1,N} \qquad 
(N=\sum_{\iota \in \mathcal{K}(n)} \kappa(\iota) ), 
\]
with the inner product given by 
\[
\left\{
\begin{array}{ll}
(e_0,e_0)=1 & \\[2mm]
(e_{\iota}^{k},e_{\iota}^{k})=-1 \qquad & 
(\iota \in \mathcal{K}(n), \quad 1 \le k \le 
\kappa(\iota)) \\[2mm]
(e_0,e_{\iota}^{k})=
(e_{\iota}^{k},e_{\iota'}^{k'})=0 \qquad & 
((\iota,k) \neq (\iota',k')). \\
\end{array}
\right. 
\]
Then an isomorphism $\phi_{\pi_{\tau}} : \mathbb{Z}^{1,N} 
\to H^2(X_{\tau}:\mathbb{Z})$ 
defined by $\phi_{\pi_{\tau}}(e_0)=[H]$ and 
$\phi_{\pi_{\tau}}(e_{\iota}^k)=[E_{\iota}^{k}]$ is the marking 
corresponding to $\pi_{\tau}$. 
Moreover, for each $\iota \in \mathcal{K}(n)$, put 
$\bar{\sigma}(\iota):=\sigma^{k}(\iota)$, 
where $k \ge 0$ is determined by the relations 
$\kappa(\sigma^{\ell}(\iota))=0$ for $0 \le \ell <k$, 
and $\kappa(\sigma^{k}(\iota)) \ge 1$. 
Then an automorphism $r_{\tau} : \mathbb{Z}^{\tau} \to \mathbb{Z}^{\tau}$ 
is defined by 
\[
r_{\tau} : \left\{
\begin{array}{lll}
e_0 & \mapsto e_0 & ~ \\[2mm]
e_{\bar{\sigma}(\iota_1)}^1 & \mapsto 
e_{\iota}^{\kappa(\iota)} \, & (\kappa(\iota) \ge 1) \\[2mm]
e_{\iota}^k & \mapsto e_{\iota}^{k-1} & 
(2 \le k \le \kappa(\iota)), 
\end{array}
\right. 
\]
Note that the map $\iota \mapsto \bar{\sigma}(\iota_1)=\bar{\sigma}(\sigma(\iota))$ becomes a permutation of 
$\{ \iota \in \mathcal{K}(n) \, | \, 
\kappa(\iota) \ge 1\}$, and so 
$e_{\bar{\sigma}(\iota_1)}^1$ is well-defined. 
The automorphism $r_{\tau}$ is an element of the subgroup 
$\langle \rho_1,\dots,\rho_{N-1} \rangle \subset W_N$ generated by 
$\rho_1,\dots,\rho_{N-1}$. 
On the other hand, for $1 \le m \le n$, 
an automorphism $q_m : \mathbb{Z}^{\tau} \to \mathbb{Z}^{\tau}$ is defined by 
\begin{equation} \label{eqn:q}
q_m : \left\{
\begin{array}{lll}
e_0 & \mapsto 2 e_0 - \sum_{\ell=1}^3 e_{\bar{\sigma}(m,\ell)}^1 & ~ \\[2mm]
e_{\bar{\sigma}(m,i)}^1 & \mapsto e_0 
- e_{\bar{\sigma}(m,j)}^1 - e_{\bar{\sigma}(m,k)}^1 \, & 
(\{i,j,k\}=\{1,2,3\}) \\[2mm]
e_{\iota}^k & \mapsto e_{\iota}^k & (\text{otherwise}). 
\end{array}
\right. 
\end{equation}
The automorphism $q_m$ is conjugate to $\rho_0$ under the action of 
$\langle \rho_1,\dots,\rho_{N-1} \rangle$. 
We notice that 
if $i \neq j \in \{1,2,3\}$ then 
$\bar{\sigma}(m,i) \neq \bar{\sigma}(m,j)$. 
Indeed, assume the contrary that $\bar{\sigma}(m,i) = \bar{\sigma}(m,j)$. 
Let $k_i \ge 0$ be the integer determined by the relations 
$\sigma^{k_i}(m,i)=\bar{\sigma}(m,i)$ and 
$\kappa(\sigma^{\ell}(\iota))=0$ for $0 \le \ell <k_i$. 
One may assume that $k_i > k_j$ and thus 
$(m,j)=\sigma^k(m,i)$ with $k=k_i-k_j$. 
As $(m_{\ell},i_{\ell}):=\sigma^{\ell}(m,i)$ satisfies 
$\kappa(m_{\ell},i_{\ell})=0$ for $0 \le \ell \le k-1 \le k_i-1$, 
we have $m=m_0 < m_1 < \cdots < m_{k}=m$, which is a contradiction. 
\par
Now we define the lattice automorphism $w_{\tau} : \mathbb{Z}^{1,N} \to 
\mathbb{Z}^{1,N}$. 
\begin{definition} \label{def:latiso}
For an orbit data $\tau$, we define the lattice automorphism 
$w_{\tau} : \mathbb{Z}^{\tau} \to \mathbb{Z}^{\tau}$ by
\begin{equation*}
~~~~~~~~~~~~~~~~~~~~~~~~~~~
w_{\tau} := r_{\tau} \circ q_1 \circ 
\cdots \circ q_n : \mathbb{Z}^{\tau} \to \mathbb{Z}^{\tau}.   
~~~~~~~~~~~~~~~~~~~~~~~~~~~
\end{equation*}
We sometimes write $w_{\tau} : \mathbb{Z}^{1,N} \to \mathbb{Z}^{1,N}$. 
\end{definition}
Indeed, it will be seen that $w_{\tau} \in W_N$ 
is realized by $(\pi_{\tau},F_{\tau})$, that is, 
$\phi_{\pi_{\tau}} \circ w_{\tau}= F_{\tau}^* \circ \phi_{\pi_{\tau}} 
: \mathbb{Z}^{1,N} \to H^2(X_{\tau};\mathbb{Z})$. 
Summing up these discussions, we have the following proposition. 
\begin{proposition} \label{pro:auto}
Assume that $\overline{f}$ is a realization of $\tau$. 
Then the blowup 
$\pi_{\tau} : X_{\tau} \to \mathbb{P}^2$ of 
$N=\sum_{\iota \in \mathcal{K}(n)} \kappa(\iota)$ points 
$\{ p_{\iota}^m \, | \, \iota =(i,j) \in \mathcal{K}(n), \, 
m=\theta_{i,0}(k), 1 \le k \le \kappa(\iota) \}$
lifts $f=f_n \circ \cdots \circ f_1$ to the automorphism 
$F_{\tau} : X_{\tau} \to X_{\tau}$. 
Moreover, $(\pi_{\tau},F_{\tau})$ realizes $w_{\tau}$ and $F_{\tau}$ 
has positive entropy 
$h_{\mathrm{top}}(F_{\tau})= \log \lambda(w_{\tau}) > 0$. 
\end{proposition}
{\it Proof}. 
We will only show that $(\pi_{\tau},F_{\tau})$ realizes 
the Weyl group element $w_{\tau}$ given in Definition \ref{def:latiso}. 
For the blowup $\pi_{\ell} : X_{\ell} \to \mathbb{P}_{\ell}^2$, 
let $H_{\ell} \subset X_{\ell}$ be 
the total transform of a line in $\mathbb{P}_{\ell}^2$ and, for $k \ge 1$, 
$E_{(i,j),\ell}^k \subset X_{\ell}$ be the exceptional divisor over 
the point $p_{i,j}^{m}$ with 
\[
m=
\left\{
\begin{array}{ll}
\theta_{i,\ell}(k-1)  \quad & (i \le \ell) \\[2mm]
\theta_{i,\ell}(k)  \quad & (i > \ell).  
\end{array}
\right. 
\] 
The cohomology group of $X_{\ell}$ admits an expression 
$H^2(X_{\ell};\mathbb{Z})= \mathbb{Z} [H_{\ell}] \oplus 
\bigr( \oplus_{\iota \in \mathcal{K}(n)} 
\oplus_{k=1}^{\kappa(\iota,\ell)} \mathbb{Z} 
[E_{\iota,\ell}^{k}]\bigl)$, 
where $\kappa(\iota,\ell)$ is the number of points among 
$p_{\iota}^0, p_{\iota}^1, \dots , p_{\iota}^{\mu(\iota)}$ 
lying on $\mathbb{P}_{\ell}^2$. 
Since the indeterminacy sets of $f_{\ell}^{\pm 1}$ 
are expressed as $I(f_{\ell}^{-1})=\{p_{\ell,j}^0 \, | \, j=1,2,3 \}$ and 
$I(f_{\ell})=\{ p_{\sigma^{-1}(\ell,j)}^{\mu(\sigma^{-1}(\ell,j))} \, | \, 
j=1,2,3\}$, the action 
$F_{\ell}^* : H^2(X_{\ell};\mathbb{Z}) \to H^2(X_{\ell-1};\mathbb{Z})$ 
is given by 
\[
F_{\ell}^* : \left\{
\begin{array}{cll}
[H_{\ell}]  & \mapsto 2 [H_{\ell-1}] - \sum_{j=1}^3 
[E_{\sigma^{-1}(\ell,j),\ell-1}^{\kappa(\sigma^{-1}(\ell,j),\ell-1)}] & ~ \\[2mm]
[E_{(\ell,i),\ell}^1]  & \mapsto [H_{\ell-1}] 
- [E_{\sigma^{-1}(\ell,j),\ell-1}^{\kappa(\sigma^{-1}(\ell,j),\ell-1)}] 
- [E_{\sigma^{-1}(\ell,k),\ell-1}^{\kappa(\sigma^{-1}(\ell,k),\ell-1)}] & 
(\{i,j,k\}=\{1,2,3\}) \\[2mm]
[E_{(\ell,j),\ell}^k] & \mapsto [E_{(\ell,j),\ell-1}^{k-1}] & 
(j \in \{1,2,3\}, \, k \ge 2) \\[2mm]
[E_{(i,j),\ell}^k] & \mapsto [E_{(i,j),\ell-1}^k] & 
(i \neq \ell, \, j \in \{1,2,3\}, \, k \ge 1) \\[2mm]
\end{array}
\right. 
\]
(see Lemma \ref{lem:QA}). Now, for $\ell \ge 1$, an isomorphism 
$G_{\ell} : H^2(X_{\ell};\mathbb{Z}) \to H^2(X_{\tau};\mathbb{Z})$ is defined by 
$G_{\ell}([H_{\ell}])=[H]$ and $G_{\ell}([E_{\iota,\ell}^k])=[E_{\iota'}^{k'}]$, 
where $E_{\iota'}^{k'}$ is the exceptional divisor over 
$p_{\iota}^{m+n-\ell} \in \mathbb{P}_{n}^2$ 
if $E_{\iota,\ell}^k$ is the exceptional divisor over 
$p_{\iota}^{m} \in \mathbb{P}_{\ell}^2$. 
In this definition, 
$p_{\iota}^{m}$ should be interpreted as 
$p_{\sigma(\iota)}^{m-\mu(\iota)-1}$ provided $m > \mu(\iota)$. The isomorphism $G_{\ell}$ 
sends $[E_{\iota,\ell}^k]$ as 
\[
G_{\ell}([E_{\iota,\ell}^k]) = 
\left\{
\begin{array}{ll}
[E_{\bar{\sigma}(\iota_1)}^1]  \quad & 
(\text{ if } i_1 > \ell \text{ and  } k=\kappa(\iota,\ell)) \\[2mm]
[E_{\iota}^{k+1}] \quad & 
(\text{ otherwise, if } i > \ell ) \\[2mm]
[E_{\iota}^{k}] \quad & 
(\text{ otherwise }). \\[2mm]
\end{array}
\right. 
\]
Since $G_{n}=\mathrm{id}$, the action $F_{\tau}^*$ is expressed as 
\[
F_{\tau}^*= F_1^* \circ \cdots \circ F_n^* 
= (F_1^* \circ G_1^{-1}) \circ \widehat{F}_2^* \circ \cdots \circ 
\widehat{F}_n^*,  
\]
where $\widehat{F}_{\ell}^* = G_{\ell-1} \circ F_{\ell}^* \circ G_{\ell}^{-1}$. 
It should be noted that 
\[
\begin{array}{ll}
G_{\ell-1}( [E_{\sigma^{-1}(\ell,j),\ell-1}^{\kappa(\sigma^{-1}(\ell,j),\ell-1)}])=
G_{\ell}([E_{(\ell,j),\ell}^1]) & = 
\left\{ 
\begin{array}{ll}
[E_{\bar{\sigma}( \sigma ( \ell,j ))}^1] \quad & (\text{ if } \ell_1 > \ell \text{ and } 
\kappa((\ell,j),\ell)=1 ) \\[2mm]
[E_{(\ell,j)}^1] \quad & (\text{ otherwise }) 
\end{array}
\right. \\[6mm]
~ & = ~~~[E_{\bar{\sigma}(\ell,j)}^1], 
\end{array}
\]
as the conditions $\ell_1 > \ell$ and $\kappa((\ell,j),\ell)=1$ 
are equivalent to saying that $\kappa(\ell,j)=0$, where $\sigma(\ell,j)=(\ell_1,j_1)$. 
Therefore, for $\ell \ge 2$, one has 
\begin{equation} \label{eqn:elemquad}
\widehat{F}_{\ell}^*  : \left\{
\begin{array}{cll}
[H]  & \mapsto 2 [H] - \sum_{j=1}^3 
[E_{\bar{\sigma}(\ell,j)}^{1}] & \quad ~ \\[2mm]
[E_{\bar{\sigma}(\ell,i)}^1]  & \mapsto [H] 
- [E_{\bar{\sigma}(\ell,j)}^{1}] 
- [E_{\bar{\sigma}(\ell,k)}^{1}] & \quad 
(\{i,j,k\}=\{1,2,3\}) \\[2mm]
[E_{(i,j)}^k] & \mapsto [E_{(i,j)}^k] & \quad (\text{ otherwise }), 
\end{array}
\right. 
\end{equation}
since $G_{\ell-1}([E_{(\ell,j),\ell-1}^{k-1}])=G_{\ell}([E_{(\ell,j),\ell}^k])$ 
and $G_{\ell-1}([E_{(i,j),\ell-1}^{k}])=G_{\ell}([E_{(i,j),\ell}^k])$ 
when $i \neq \ell$. Finally, by observing that 
\[
F_{1}^* \circ G_1^{-1} : \left\{
\begin{array}{cll}
[H]  & \mapsto 2 [H] - \sum_{j=1}^3 
[E_{\sigma^{-1}(1,j)}^{\kappa(\sigma^{-1}(1,j))}] & ~ \\[2mm]
[E_{\bar{\sigma}(1,i)}^1]  & \mapsto [H] 
- [E_{\sigma^{-1}(1,j)}^{\kappa(\sigma^{-1}(1,j))}] 
- [E_{\sigma^{-1}(1,k)}^{\kappa(\sigma^{-1}(1,k))}] & 
(\{i,j,k\}=\{1,2,3\}) \\[2mm]
[E_{\bar{\sigma}(\sigma(i,j))}^1] & \mapsto [E_{(i,j)}^{\kappa(i,j)}] & 
(i_1 \neq 1 ) \\[2mm]
[E_{(i,j)}^k] & \mapsto [E_{(i,j)}^{k-1}] & 
(\text{ otherwise }), \\[2mm]
\end{array}
\right. 
\]
we define an isomorphism $G_0 : H^2(X_{\tau}; \mathbb{Z}) \circlearrowleft$ 
by 
\begin{equation} \label{eqn:elemperm}
G_0 : \left\{
\begin{array}{cll}
[H]  & \mapsto [H] & ~ \\[2mm]
[E_{\iota}^{\kappa(\iota)}]  & \mapsto [E_{\bar{\sigma}(\iota_1)}^{1}] & (\kappa(\iota) \ge 1) \\[2mm]
[E_{\iota}^k] & \mapsto [E_{\iota}^{k+1}] & 
(1 \le k \le \kappa(\iota)-1). \\[2mm]
\end{array}
\right. 
\end{equation}
Then, $\widehat{F}_1^*=G_0 \circ F_1^* \circ G_1^{-1}$ satisfies 
(\ref{eqn:elemquad}) with $\ell=1$, and $F_{\tau}^*$ satisfies 
$F_{\tau}^*=G_0^{-1} \circ \widehat{F}_1^* \circ \cdots \circ \widehat{F}_n^*$. 
From (\ref{eqn:elemquad}) and (\ref{eqn:elemperm}), one has 
$\widehat{F}_{\ell}^* = \phi_{\pi_{\tau}} \circ q_{\ell} \circ \phi_{\pi_{\tau}}^{-1}$ 
and 
$G_0^{-1} = \phi_{\pi_{\tau}} \circ r_{\tau} \circ \phi_{\pi_{\tau}}^{-1}$, 
which shows that $w_{\tau}$ is 
realized by $(\pi_{\tau},F_{\tau})$. 
\hfill $\Box$ \par\medskip 
We conclude this section by establishing the statement that a given element $w \in W_N$ can be expressed as 
$w=w_{\tau}$ for some orbit data $\tau$. 
To this end, we spend a short while working with an element $w \in W_N$ not acting by a permutation 
on a non-empty subset of the basis elements $\{ e_j \}_{j=1}^N$, namely, there is no element 
$e_j \in \{e_i \}$ with $w^\ell(e_j) \in \{e_i \}$ for any $\ell \ge 1$, 
and explain how to construct an orbit data $\tau$ with $w=w_{\tau}$ briefly. 
First, note that $w$ can be expressed as 
\begin{equation} \label{eqn:expw}
w=r \circ q_1 \circ \cdots q_n, 
\end{equation}
where $r$ is a permutation of $\{ e_j \}_{j=1}^N$ (see the proof of Proposition \ref{pro:iden}). 
Moreover, there are elements $\{ e_{m,i}^1 \}_{i=1}^3 \subset \{ e_j \}_{j=1}^N$ such that $q_m$ sends 
$\{e_j \}_{j=0}^N$ as 
\[
q_m : \left\{
\begin{array}{lll}
e_0 & \mapsto 2 e_0 - \sum_{\ell=1}^3 e_{m,\ell}^1 & ~ \\[2mm]
e_{m,i}^1 & \mapsto e_0 
- e_{m,j}^1 - e_{m,k}^1 \, & 
(\{i,j,k\}=\{1,2,3\}) \\[2mm]
e_j & \mapsto e_j & (\text{otherwise}). 
\end{array}
\right. 
\]
We notice that it may happen that $e_{i,j}^1 =e_{i',j'}^1$ when $i \neq i'$. 
For each $(i,j) \in \mathcal{K}(n)$, it turns out that there is a unique element $(i',j')$ such that 
\begin{equation} \label{eqn:orbite}
e_{i',j'}^1=q_{i'-1} \circ \cdots \circ q_1 \circ r^{-1} \circ q_n \circ \cdots \circ q_1 \circ r^{-1} \circ q_n \circ \cdots \circ q_{i+1}
(e_{i,j}^1)
\end{equation}
with minimal length, where the length is the number of automorphisms $q_m$ in (\ref{eqn:orbite}). 
Roughly speaking, $e_{i,j}^1$ and $e_{i',j'}^1$ are regarded as backward and forward indeterminacy points respectively, 
and the length is the number of points from $e_{i,j}^1$ to $e_{i',j'}^1$. 
So we define $\sigma(i,j)=(i',j')$ and $\mu(i,j)$ to be the length in (\ref{eqn:orbite}). 
Note that $\kappa(i,j)$ is the number of $r^{-1}$ in (\ref{eqn:orbite}), 
and $\sigma$ becomes a permutation of $\mathcal{K}(n)$ because of the minimality of the length. 
Then $\tau=(n,\sigma,\kappa)$ is an orbit data satisfying $w=w_{\tau}$. 
\begin{example} \label{ex:fwto}
Consider the element $w \in W_5$ given by 
\[
w : \left\{
\begin{array}{lll}
e_0 & \mapsto 3 e_0 - 2 e_1 - e_2 -e_3 -e_4 -e_5  & ~ \\[2mm]
e_1 & \mapsto 2 e_0 - e_1 -e_2 -e_3 -e_4 -e_5  & ~ \\[2mm]
e_i & \mapsto e_0 - e_1-e_{7-i} & (i \in \{2,3,4,5\}). 
\end{array}
\right. 
\]
\end{example}
It can be checked that $w$ is expressed as $w=r \circ q_1 \circ q_2$, where 
\[
r : e_2 \longleftrightarrow e_4, \quad e_3 \longleftrightarrow e_5, \quad e_i \circlearrowleft (i=0,1), 
\]
and $(e_{1,1},e_{1,2},e_{1,3})=(e_1,e_2,e_3)$, $(e_{2,1},e_{2,2},e_{2,3})=(e_1,e_4,e_5)$. 
Then we have 
\[
\begin{array}{l}
e_{2,1}=e_{1,1}, \quad e_{1,k}=r^{-1} (e_{2,k}), \quad (k \in \{1,2,3\}), \\[2mm]
e_{2,j} \neq e_{1,k}, \quad e_{1,j} \neq r^{-1} \circ q_2(e_{1,k}), \quad e_{2,k}=q_1 \circ r^{-1} \circ q_2 (e_{1,k}), \quad 
(j \in \{1,2,3\}, \, k \in \{2,3\}). 
\end{array}
\]
This means that $w$ admits an expression $w=w_{\tau}$, where $\tau=(2,\sigma,\kappa)$ is given by 
\[
\left\{
\begin{array}{l}
\phantom{(} \sigma : (1,k) \mapsto (2,k) \mapsto (1,k), \quad (k \in \{1,2,3\}), \\[2mm]
\phantom{(} \kappa(1,1)=0, \qquad \kappa(i,j)=1 ~~ (\text{otherwise}), \\[1mm]
(\mu(1,1)=\mu(2,1)=\mu(2,2)=\mu(2,3)=0, \quad 
\mu(1,2)=\mu(1,3)=2 ). 
\end{array}
\right. 
\]
\begin{proposition} \label{pro:iden}
For any $w \in W_N$, there is an orbit data 
$\tau=(n,\sigma,\kappa)$ 
such that $w= w_{\tau}: \mathbb{Z}^{1,N} \to \mathbb{Z}^{1,N}$ 
under some identification 
$\{e_j \, | \, 1 \le j \le N\}=\{e_{\iota}^k \, | \, 
\iota \in \mathcal{K}(n), \, 1 \le k \le \kappa(\iota)\}$ 
with $N=\sum_{\iota \in \mathcal{K}(n)} \kappa(\iota)$. 
\end{proposition}
{\it Proof}. 
First we prove that any element $w$ admits expression (\ref{eqn:expw}). 
Since $w$ is an element of $W_N$, it can be expressed as 
\[
w= r_0 \cdot \rho_0 \cdot r_1 \cdots \rho_0 \cdot 
r_{m-1} \cdot \rho_0 \cdot r_{m}, 
\]
where $r_\ell$ is a permutation of $\{e_j\}_{j=1}^N$. 
The expression can be written as 
\[
w  = r \cdot 
\big\{ (r_1  \cdots r_{m})^{-1} \cdot \rho_0 \cdot 
(r_1  \cdots r_{m}) \big\} \cdots 
\big\{ (r_{m-1} \cdot r_{m})^{-1} \cdot \rho_0 \cdot (r_{m-1} \cdot 
r_{m}) \big\} \cdot \big\{ r_{m}^{-1} \cdot \rho_0 \cdot r_{m} \big\}, 
\]
where $r:=r_0  \cdots r_{m}$ is also a permutation of $\{e_j\}$. 
By putting $q_{i}:= (r_{i} \cdots r_{m})^{-1} 
\cdot \rho_0 \cdot (r_{i} \cdots r_{m})$ and 
$e_{i,j}^1:=(r_{i}  \cdots r_{m})^{-1}(e_{j})$ 
for $i=1,\dots,m$ and $j=1,2,3$, 
we have expression (\ref{eqn:expw}). 
\par
Next, under the assumption that $w$ does not act by a permutation on a non-empty subset of $\{ e_j \}_{j=1}^N$, 
we show that the orbit data $\tau$ constructed above realizes $w$. 
Put 
\[
e_{i,j}^{k+1}:= q_n \circ \cdots \circ q_1 \circ r^{-1} \circ q_n \circ \cdots \circ q_1 \circ r^{-1} \circ q_n 
\circ \cdots \circ q_{i+1}(e_{i,j}^1) 
\]
for $1 \le k \le \kappa(i,j)-1$, where the number of $r^{-1}$ in the righthand side is $k$. 
Since $e_{\sigma(i,j)}^1=e_{i,j}^1$ when $\kappa(i,j)=0$, one has $e_{\bar{\sigma}(i,j)}^1=e_{i,j}^1$, 
which shows $q_m$ sends  $\{e_{\iota}^k \}$ as in (\ref{eqn:q}). 
Moreover it follows from the minimality of the length in (\ref{eqn:orbite}) that 
$e_{i,j}^{k}= r \circ q_1 \circ \cdots \circ q_n(e_{i,j}^{k+1})=r \circ q_1 \circ \cdots \circ q_{n-1}(e_{i,j}^{k+1}) 
= \cdots = r(e_{i,j}^{k+1})$, 
and also that 
$e_{i,j}^{\kappa(i,j)}= r \circ q_1 \circ \cdots \circ q_{i'-1} (e_{\sigma(i,j)}^1)= r (e_{\sigma(i,j)}^1)=
r(e_{\bar{\sigma}(\sigma(i,j))}^1)$
for $\kappa(i,j) \ge 1$, which means that $r=r_{\tau}$. 
Now we claim that $\{e_j\}=\{e_{\iota}^k \}$, that is, any element $e_{m} \in \{e_j\}$ can be expressed as $e_m = e_{\iota}^k$ 
for some $\iota \in \mathcal{K}(n)$ and $1 \le k \le \kappa(\iota)$. 
Indeed, assume the contrary that $e_m \neq e_{\iota}^k$ for any $\iota \in \mathcal{K}(n)$ and $1 \le k \le \kappa(\iota)$. 
Then $e_m$ satisfies $w^{\ell}(e_m) \notin \{e_{\iota}^1 \}$ 
and thus $w^{\ell}(e_m) \in \{e_{j} \}$  for all $\ell \ge 1$, which is a contradiction. 
Moreover, we can easily check that $e_{\iota}^k \neq e_{\iota'}^{k'}$ for any 
$(\iota,k) \neq (\iota',k')$ with $1 \le k \le \kappa(\iota)$ and $1 \le k' \le \kappa(\iota')$, which shows that 
$N=\sum_{\iota \in \mathcal{K}(n)} \kappa(\iota)$. 
These observations show that $w$ is expressed as $w=w_{\tau}$.  
\par 
Finally we consider a general element $w \in W_N$. Then $w$ can be expressed as 
$w=w_1 \cdot \hat{w}$, where $w_1$ does not act by a permutation on a non-empty subset of $\{ e_j \}_{j=1}^N$ 
and hence admits an expression $w_1=w_{\check{\tau}}$ for some orbit data 
$\check{\tau}=(m,\check{\sigma},\check{\kappa})$, and 
$\hat{w}$ is a permutation of $\{ e_j \}$, sending the basis elements, 
after reordering $\{ e_j \}$, as 
\[
\hat{w} : \left\{
\begin{array}{lll}
e_0  & \mapsto e_0 & ~ \\[2mm]
e_{\iota}^k  & \mapsto e_{\iota}^k & 
(\iota \in \mathcal{K}(m), 1 \le k \le \check{\kappa}(\iota) ) \\[2mm]
e_{m+2i,1}^k & \mapsto e_{m+2i,1}^{k-1} & 
(i=1, \dots, \ell, \quad k \in \mathbb{Z}/\hat{\kappa}(i) \mathbb{Z}) 
\end{array}
\right. 
\]
for some $\ell \ge 0$ and $\hat{\kappa} : \{1,\dots,\ell\} \to \mathbb{Z}_{\ge 1}$. 
Now, by composing automorphisms of the form $q_{m+2i-1} \circ q_{m+2i}$ with 
$e_{m+2i-1,j}^1=e_{m+2i,j}^1$ for all $i=1,2,3$, which are also permutations of $\{e_j \}_{j=1}^N$, 
we add elements $e_{m+2i,1}^1$ and construct $\hat{w}$. 
Namely, $\tau=(n,\sigma,\kappa)$ is defined by $n:=m+2\ell$ and 
\[
\sigma(i,j) := 
\left\{
\begin{array}{ll}
(i+1,j) & 
\left(
\begin{array}{cl}
\text{either } & j=1 \text{ and } i=m+1,m+3,\dots, m+ 2 \ell-1, \\
\text{ or } & j=2,3 \text{ and } i=m,m+1,\dots,m+ 2 \ell -1  
\end{array} \right) \\[3mm]
(i-1,j) & (j=1 \text{ and } i=m+2,m+4,\dots, m +2 \ell  ) \\[2mm]
\check{\sigma}(m,j) & (j=2,3  \text{ and } i=m+2 \ell ) \\[2mm] 
\check{\sigma}(i,j) & (\text{otherwise}), 
\end{array}
\right. 
\]
\[
\kappa(i,j) := 
\left\{
\begin{array}{ll}
0 & 
\left(
\begin{array}{cl}
\text{either } & j=1 \text{ and } i=m+1,m+3,\dots, m+ 2 \ell-1, \\
\text{ or } & j=2,3 \text{ and } i=m,m+1,\dots,m+ 2 \ell -1  
\end{array} \right) \\[3mm]
\hat{\kappa}((i-m)/2) & (j=1 \text{ and } i=m+2,m+4,\dots, m +2 \ell  ) \\[2mm]
\check{\kappa}(m,j) & (j=2,3  \text{ and } i=m+2 \ell ) \\[2mm] 
\check{\kappa}(i,j) & (\text{otherwise}). 
\end{array}
\right. 
\]
A little calculation shows that $w=w_{\check{\tau}} \cdot \hat{w}$ can be expressed 
as $w=w_{\tau}=r_{\tau} \circ q_1 \circ \cdots \circ q_n$. 
Thus the proposition is established. 
\hfill $\Box$ \par\medskip 
\section{Tentative Realizability} \label{sec:tenta}
As mentioned in the previous section, 
an automorphism can be constructed in terms of a realization of an orbit data. 
At this stage, of particular interest is 
the existence of such a realization. 
In this and next sections, we investigate this existence 
by restricting our attention to birational maps preserving a cuspidal cubic. 
The aim of this section is to 
define a concept of tentative realization of an orbit data $\tau$, 
which is a necessary condition for realization, 
and to show that the tentative realization of $\tau$ exists under the 
condition that some finitely many roots determined by $\tau$ 
are not periodic ones of the Weyl group element $w_{\tau}$ 
(see condition (\ref{eqn:TR})). 
\par
In this section, we mainly consider the smooth points of the cuspidal cubic, 
which is also described in a more general context as follows. 
Let $X$ be a smooth surface, 
$C$ be a curve in $X$, and $x$ be a proper point of the smooth locus $C^*$ 
of $C$. Moreover, put $(X_0,C_0^*,x_0):=(X,C^*,x)$, and for $m > 0$, 
inductively determine $(X_m,C_m^*,x_m)$ from the blowup 
$\pi_m : X_m \to X_{m-1}$ of $x_{m-1} \in C_{m-1}^*$, the strict transform 
$C_m^*$ of $C_{m-1}^*$ under $\pi_m$, and a unique point 
$x_{m} \in C_m^* \cap E_m$, where $E_m$ stands for the exceptional curve 
of $\pi_m$. 
The unique point 
$x_m$ is called the {\sl point in the $m$-th infinitesimal 
neighbourhood of $x$ on $C^*$}, or an {\sl $m$-th point on $C^*$}. 
Moreover, if a cluster $I$ consists of proper or infinitely near points on $C^*$, 
then we say that $I$ is {\sl a cluster in $C^*$}.  
\par
Now let $C$ be a cubic curve on $\mathbb{P}^2$ with a cusp singularity. 
In what follows, a coordinate on $\mathbb{P}^2$ is chosen so that 
$C=\{[x:y:z] \in \mathbb{P}^2 \, | \, y z^2=x^3\} \subset \mathbb{P}^2$ 
with a cusp $[0:1:0]$. Then the smooth locus $C^*=C \setminus \{[0:1:0]\}$ 
is parametrized as $\mathbb{C} \ni t \mapsto [t:t^3:1] \in C^*$. 
We denote by $\mathcal{B}(C)$ the set of birational self-maps $f$ of 
$\mathbb{P}^2$ such that $f(C):=\overline{f(C \setminus I(f))}=C$ and 
$I(f) \subset C^*$, and denote by 
$\mathcal{Q}(C) \subset \mathcal{B}(C)$ and 
$\mathcal{L}(C) \subset \mathcal{B}(C)$ the subsets consisting of 
the quadratic maps in $\mathcal{B}(C)$ and of 
the linear maps in $\mathcal{B}(C)$, respectively. 
Any map $f \in \mathcal{B}(C)$ 
restricted to $C^*$ is an automorphism of $C^*$ expressed as 
\[
f|_{C^*} : C^* \ni [t:t^3:1] \mapsto 
[\delta(f) \cdot t+ c(f): (\delta(f) \cdot t+c(f))^3:1] \in C^*, 
\]
for some $\delta(f) \in \mathbb{C}^{\times}$ and $c(f) \in \mathbb{C}$. 
The value $\delta(f)$ is called the {\sl determinant} of $f$. 
It is independent of the choice of 
coordinates. Moreover, when $f \in \mathcal{Q}(C)$, 
it turns out that the indeterminacy cluster $I(f^{- 1})$ is also contained in $C^*$ (see Lemma \ref{lem:quad}). 
\par
We give the following definition 
for an $n$-tuple $\overline{f}=(f_1,\dots,f_n) \in \mathcal{Q}(C)^n$ 
of quadratic birational maps $f_i$ preserving $C$. 
\begin{definition} \label{def:tenta} 
An $n$-tuple $\overline{f}=(f_1,\dots,f_n) \in \mathcal{Q}(C)^n$ 
is called a {\sl tentative realization} of an orbit data 
$\tau=(n,\sigma,\kappa)$ if 
$p_{\iota}^{\mu(\iota)} 
\approx p_{\sigma(\iota)}^+$ 
for any $\iota \in \mathcal{K}(n)$, 
where $p_{\iota}^m$ is given in (\ref{eqn:traj}) 
with $f_{\ell}$ restricted to $C$ and thus is well-defined. 
\end{definition}
We should note that a realization $\overline{f}$ of $\tau$ is of course a 
tentative realization of $\tau$, and thus the existence of a 
tentative realization is of interest to us. 
\par
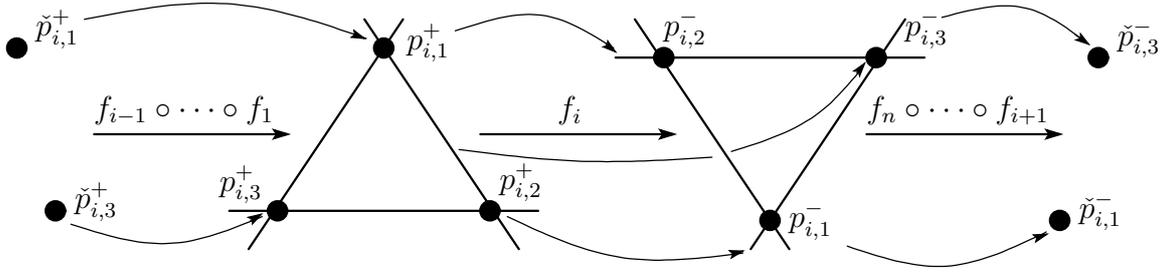
\begin{figure}[t]
\begin{center}
\unitlength 0.1in
\begin{picture}( 57.4500, 13.7400)(  7.5500,-16.9000)
%
\special{pn 13}%
\special{pa 1200 1000}%
\special{pa 2200 1000}%
\special{fp}%
\special{sh 1}%
\special{pa 2200 1000}%
\special{pa 2134 980}%
\special{pa 2148 1000}%
\special{pa 2134 1020}%
\special{pa 2200 1000}%
\special{fp}%
%
\special{pn 13}%
\special{pa 2800 400}%
\special{pa 2000 1600}%
\special{fp}%
%
\special{pn 13}%
\special{pa 2600 400}%
\special{pa 3400 1600}%
\special{fp}%
%
\special{pn 13}%
\special{pa 4600 1600}%
\special{pa 5400 400}%
\special{fp}%
%
\special{pn 13}%
\special{pa 4800 1600}%
\special{pa 4000 400}%
\special{fp}%
%
\special{pn 13}%
\special{pa 3200 1000}%
\special{pa 4200 1000}%
\special{fp}%
\special{sh 1}%
\special{pa 4200 1000}%
\special{pa 4134 980}%
\special{pa 4148 1000}%
\special{pa 4134 1020}%
\special{pa 4200 1000}%
\special{fp}%
%
\special{pn 13}%
\special{pa 5200 1000}%
\special{pa 6200 1000}%
\special{fp}%
\special{sh 1}%
\special{pa 6200 1000}%
\special{pa 6134 980}%
\special{pa 6148 1000}%
\special{pa 6134 1020}%
\special{pa 6200 1000}%
\special{fp}%
%
\special{pn 20}%
\special{sh 1.000}%
\special{ar 2700 550 46 46  0.0000000 6.2831853}%
%
\special{pn 20}%
\special{sh 1.000}%
\special{ar 2150 1400 46 46  0.0000000 6.2831853}%
%
\special{pn 20}%
\special{sh 1.000}%
\special{ar 3250 1400 46 46  0.0000000 6.2831853}%
%
\special{pn 20}%
\special{sh 1.000}%
\special{ar 4150 600 46 46  0.0000000 6.2831853}%
%
\special{pn 20}%
\special{sh 1.000}%
\special{ar 4700 1450 46 46  0.0000000 6.2831853}%
%
\special{pn 20}%
\special{sh 1.000}%
\special{ar 5250 600 46 46  0.0000000 6.2831853}%
\put(12.0000,-9.5000){\makebox(0,0)[lb]{$f_{i-1} \circ \cdots \circ f_1$}}%
\put(52.0000,-9.5000){\makebox(0,0)[lb]{$f_{n} \circ \cdots \circ f_{i+1}$}}%
\put(28.2000,-6.2000){\makebox(0,0)[lb]{$p_{i,1}^+$}}%
\put(33.0000,-13.5000){\makebox(0,0)[lb]{$p_{i,2}^+$}}%
\put(18.5000,-13.6000){\makebox(0,0)[lb]{$p_{i,3}^+$}}%
\put(48.0000,-15.5000){\makebox(0,0)[lb]{$p_{i,1}^-$}}%
\put(41.5000,-5.5000){\makebox(0,0)[lb]{$p_{i,2}^-$}}%
\put(54.0000,-5.5000){\makebox(0,0)[lb]{$p_{i,3}^-$}}%
%
\special{pn 13}%
\special{pa 1900 1400}%
\special{pa 3500 1400}%
\special{fp}%
%
\special{pn 13}%
\special{pa 3900 600}%
\special{pa 5500 600}%
\special{fp}%
%
\special{pn 8}%
\special{pa 3340 1430}%
\special{pa 3370 1444}%
\special{pa 3400 1458}%
\special{pa 3430 1472}%
\special{pa 3458 1486}%
\special{pa 3488 1500}%
\special{pa 3518 1512}%
\special{pa 3548 1526}%
\special{pa 3578 1538}%
\special{pa 3608 1550}%
\special{pa 3638 1562}%
\special{pa 3668 1574}%
\special{pa 3698 1584}%
\special{pa 3728 1594}%
\special{pa 3758 1604}%
\special{pa 3788 1612}%
\special{pa 3818 1620}%
\special{pa 3850 1628}%
\special{pa 3880 1634}%
\special{pa 3912 1640}%
\special{pa 3942 1646}%
\special{pa 3974 1650}%
\special{pa 4004 1652}%
\special{pa 4036 1654}%
\special{pa 4068 1656}%
\special{pa 4100 1658}%
\special{pa 4132 1658}%
\special{pa 4162 1656}%
\special{pa 4194 1656}%
\special{pa 4228 1654}%
\special{pa 4260 1652}%
\special{pa 4292 1650}%
\special{pa 4324 1648}%
\special{pa 4356 1644}%
\special{pa 4388 1642}%
\special{pa 4420 1638}%
\special{pa 4454 1634}%
\special{pa 4486 1630}%
\special{pa 4518 1626}%
\special{pa 4550 1622}%
\special{pa 4560 1620}%
\special{sp}%
%
\special{pn 8}%
\special{pa 3070 460}%
\special{pa 3102 450}%
\special{pa 3134 440}%
\special{pa 3168 430}%
\special{pa 3200 420}%
\special{pa 3232 410}%
\special{pa 3264 402}%
\special{pa 3294 394}%
\special{pa 3326 388}%
\special{pa 3358 382}%
\special{pa 3390 378}%
\special{pa 3420 376}%
\special{pa 3452 374}%
\special{pa 3482 376}%
\special{pa 3512 378}%
\special{pa 3542 382}%
\special{pa 3572 390}%
\special{pa 3602 398}%
\special{pa 3632 408}%
\special{pa 3660 420}%
\special{pa 3690 434}%
\special{pa 3718 448}%
\special{pa 3748 464}%
\special{pa 3776 480}%
\special{pa 3804 496}%
\special{pa 3832 514}%
\special{pa 3862 532}%
\special{pa 3890 550}%
\special{pa 3890 550}%
\special{sp}%
%
\special{pn 8}%
\special{pa 3090 1080}%
\special{pa 3122 1084}%
\special{pa 3154 1088}%
\special{pa 3186 1092}%
\special{pa 3218 1096}%
\special{pa 3250 1100}%
\special{pa 3282 1104}%
\special{pa 3314 1108}%
\special{pa 3346 1112}%
\special{pa 3376 1114}%
\special{pa 3408 1118}%
\special{pa 3440 1122}%
\special{pa 3472 1124}%
\special{pa 3504 1128}%
\special{pa 3536 1130}%
\special{pa 3568 1132}%
\special{pa 3600 1134}%
\special{pa 3632 1136}%
\special{pa 3664 1138}%
\special{pa 3696 1140}%
\special{pa 3728 1142}%
\special{pa 3760 1142}%
\special{pa 3792 1142}%
\special{pa 3824 1142}%
\special{pa 3856 1144}%
\special{pa 3888 1142}%
\special{pa 3920 1142}%
\special{pa 3952 1142}%
\special{pa 3984 1142}%
\special{pa 4016 1140}%
\special{pa 4048 1140}%
\special{pa 4080 1138}%
\special{pa 4112 1136}%
\special{pa 4144 1136}%
\special{pa 4176 1134}%
\special{pa 4208 1132}%
\special{pa 4240 1130}%
\special{pa 4272 1128}%
\special{pa 4304 1126}%
\special{pa 4336 1124}%
\special{pa 4368 1122}%
\special{pa 4400 1120}%
\special{pa 4400 1120}%
\special{sp}%
%
\special{pn 8}%
\special{pa 4500 1090}%
\special{pa 4532 1080}%
\special{pa 4564 1070}%
\special{pa 4594 1060}%
\special{pa 4626 1050}%
\special{pa 4656 1040}%
\special{pa 4686 1028}%
\special{pa 4716 1016}%
\special{pa 4746 1004}%
\special{pa 4776 992}%
\special{pa 4806 978}%
\special{pa 4834 964}%
\special{pa 4862 948}%
\special{pa 4888 932}%
\special{pa 4916 916}%
\special{pa 4940 898}%
\special{pa 4966 878}%
\special{pa 4990 858}%
\special{pa 5014 838}%
\special{pa 5036 814}%
\special{pa 5058 792}%
\special{pa 5080 768}%
\special{pa 5102 744}%
\special{pa 5122 720}%
\special{pa 5144 696}%
\special{pa 5164 670}%
\special{pa 5180 650}%
\special{sp}%
%
\special{pn 8}%
\special{pa 4520 1620}%
\special{pa 4560 1610}%
\special{fp}%
\special{sh 1}%
\special{pa 4560 1610}%
\special{pa 4490 1608}%
\special{pa 4508 1624}%
\special{pa 4500 1646}%
\special{pa 4560 1610}%
\special{fp}%
%
\special{pn 8}%
\special{pa 5150 680}%
\special{pa 5190 630}%
\special{fp}%
\special{sh 1}%
\special{pa 5190 630}%
\special{pa 5134 670}%
\special{pa 5158 672}%
\special{pa 5164 696}%
\special{pa 5190 630}%
\special{fp}%
%
\special{pn 8}%
\special{pa 3850 530}%
\special{pa 3900 570}%
\special{fp}%
\special{sh 1}%
\special{pa 3900 570}%
\special{pa 3860 514}%
\special{pa 3858 538}%
\special{pa 3836 544}%
\special{pa 3900 570}%
\special{fp}%
%
\special{pn 20}%
\special{sh 1.000}%
\special{ar 1000 1400 46 46  0.0000000 6.2831853}%
%
\special{pn 20}%
\special{sh 1.000}%
\special{ar 800 550 46 46  0.0000000 6.2831853}%
%
\special{pn 20}%
\special{sh 1.000}%
\special{ar 6400 600 46 46  0.0000000 6.2831853}%
%
\special{pn 20}%
\special{sh 1.000}%
\special{ar 6200 1450 46 46  0.0000000 6.2831853}%
\put(11.0000,-14.5000){\makebox(0,0)[lb]{$\check{p}_{i,3}^+$}}%
\put(9.0000,-5.5000){\makebox(0,0)[lb]{$\check{p}_{i,1}^+$}}%
\put(65.0000,-6.0000){\makebox(0,0)[lb]{$\check{p}_{i,3}^-$}}%
\put(63.0000,-15.0000){\makebox(0,0)[lb]{$\check{p}_{i,1}^-$}}%
%
\special{pn 8}%
\special{pa 1080 1470}%
\special{pa 1112 1480}%
\special{pa 1144 1490}%
\special{pa 1176 1500}%
\special{pa 1206 1510}%
\special{pa 1238 1518}%
\special{pa 1270 1528}%
\special{pa 1302 1536}%
\special{pa 1332 1544}%
\special{pa 1364 1550}%
\special{pa 1396 1556}%
\special{pa 1426 1562}%
\special{pa 1458 1566}%
\special{pa 1490 1568}%
\special{pa 1520 1570}%
\special{pa 1552 1572}%
\special{pa 1584 1570}%
\special{pa 1614 1568}%
\special{pa 1646 1566}%
\special{pa 1676 1562}%
\special{pa 1706 1556}%
\special{pa 1738 1548}%
\special{pa 1768 1540}%
\special{pa 1800 1532}%
\special{pa 1830 1522}%
\special{pa 1860 1512}%
\special{pa 1892 1502}%
\special{pa 1922 1490}%
\special{pa 1952 1478}%
\special{pa 1982 1466}%
\special{pa 2014 1454}%
\special{pa 2044 1442}%
\special{pa 2070 1430}%
\special{sp}%
%
\special{pn 8}%
\special{pa 1150 400}%
\special{pa 1182 394}%
\special{pa 1214 388}%
\special{pa 1246 382}%
\special{pa 1278 376}%
\special{pa 1310 370}%
\special{pa 1342 364}%
\special{pa 1374 358}%
\special{pa 1406 354}%
\special{pa 1438 348}%
\special{pa 1470 344}%
\special{pa 1502 338}%
\special{pa 1534 334}%
\special{pa 1566 330}%
\special{pa 1598 326}%
\special{pa 1630 324}%
\special{pa 1662 322}%
\special{pa 1694 320}%
\special{pa 1724 318}%
\special{pa 1756 316}%
\special{pa 1788 316}%
\special{pa 1820 316}%
\special{pa 1852 318}%
\special{pa 1882 318}%
\special{pa 1914 320}%
\special{pa 1946 324}%
\special{pa 1976 328}%
\special{pa 2008 332}%
\special{pa 2040 336}%
\special{pa 2070 342}%
\special{pa 2102 348}%
\special{pa 2132 354}%
\special{pa 2164 362}%
\special{pa 2194 370}%
\special{pa 2226 378}%
\special{pa 2256 386}%
\special{pa 2288 394}%
\special{pa 2318 404}%
\special{pa 2350 412}%
\special{pa 2380 422}%
\special{pa 2412 432}%
\special{pa 2442 442}%
\special{pa 2472 454}%
\special{pa 2504 464}%
\special{pa 2534 474}%
\special{pa 2566 484}%
\special{pa 2596 496}%
\special{pa 2610 500}%
\special{sp}%
%
\special{pn 8}%
\special{pa 5590 390}%
\special{pa 5624 382}%
\special{pa 5656 372}%
\special{pa 5688 364}%
\special{pa 5722 356}%
\special{pa 5754 350}%
\special{pa 5786 344}%
\special{pa 5818 338}%
\special{pa 5850 334}%
\special{pa 5882 330}%
\special{pa 5914 330}%
\special{pa 5944 330}%
\special{pa 5974 332}%
\special{pa 6004 338}%
\special{pa 6034 344}%
\special{pa 6064 352}%
\special{pa 6092 364}%
\special{pa 6120 376}%
\special{pa 6148 392}%
\special{pa 6176 406}%
\special{pa 6204 424}%
\special{pa 6232 442}%
\special{pa 6258 462}%
\special{pa 6286 482}%
\special{pa 6312 502}%
\special{pa 6340 522}%
\special{pa 6350 530}%
\special{sp}%
%
\special{pn 8}%
\special{pa 5100 1584}%
\special{pa 5132 1596}%
\special{pa 5164 1606}%
\special{pa 5196 1616}%
\special{pa 5226 1626}%
\special{pa 5258 1636}%
\special{pa 5290 1644}%
\special{pa 5322 1654}%
\special{pa 5352 1662}%
\special{pa 5384 1668}%
\special{pa 5416 1674}%
\special{pa 5446 1680}%
\special{pa 5478 1684}%
\special{pa 5508 1688}%
\special{pa 5540 1690}%
\special{pa 5570 1690}%
\special{pa 5602 1690}%
\special{pa 5632 1688}%
\special{pa 5662 1684}%
\special{pa 5694 1678}%
\special{pa 5724 1672}%
\special{pa 5754 1664}%
\special{pa 5786 1656}%
\special{pa 5816 1646}%
\special{pa 5846 1636}%
\special{pa 5876 1626}%
\special{pa 5906 1614}%
\special{pa 5936 1600}%
\special{pa 5966 1588}%
\special{pa 5996 1574}%
\special{pa 6026 1560}%
\special{pa 6056 1546}%
\special{pa 6088 1532}%
\special{pa 6118 1518}%
\special{pa 6130 1510}%
\special{sp}%
%
\special{pn 8}%
\special{pa 2550 480}%
\special{pa 2600 500}%
\special{fp}%
\special{sh 1}%
\special{pa 2600 500}%
\special{pa 2546 458}%
\special{pa 2550 480}%
\special{pa 2532 494}%
\special{pa 2600 500}%
\special{fp}%
%
\special{pn 8}%
\special{pa 2010 1450}%
\special{pa 2066 1418}%
\special{fp}%
\special{sh 1}%
\special{pa 2066 1418}%
\special{pa 1998 1434}%
\special{pa 2020 1444}%
\special{pa 2020 1468}%
\special{pa 2066 1418}%
\special{fp}%
%
\special{pn 8}%
\special{pa 6080 1540}%
\special{pa 6134 1522}%
\special{fp}%
\special{sh 1}%
\special{pa 6134 1522}%
\special{pa 6064 1526}%
\special{pa 6082 1540}%
\special{pa 6076 1562}%
\special{pa 6134 1522}%
\special{fp}%
%
\special{pn 8}%
\special{pa 6280 480}%
\special{pa 6350 530}%
\special{fp}%
\special{sh 1}%
\special{pa 6350 530}%
\special{pa 6308 476}%
\special{pa 6308 500}%
\special{pa 6284 508}%
\special{pa 6350 530}%
\special{fp}%
\put(36.0000,-9.5000){\makebox(0,0)[lb]{$f_i$}}%
\end{picture}%
\end{center}
\caption{The points $\check{p}_{i,j}^+ \in I(f)$ and 
$\check{p}_{i,j}^- \in I(f^{-1})$} 
\label{fig:comp}
\end{figure}

Now, we describe a quadratic birational map $f \in \mathcal{Q}(C)$ 
in terms of the behavior of $f|_{C^*}$. The following proposition 
states that the configuration of $I(f^{-1})$ on $C^*$ and 
the determinant $\delta(f)$ of $f$ determine the map 
$f \in \mathcal{Q}(C)$ uniquely (see \cite{D}). 
\begin{lemma}  \label{lem:quad}
A birational map $f$ belongs to $\mathcal{Q}(C)$ 
if and only if there exists $d \in \mathbb{C}^{\times}$ and 
$b=(b_{\ell})_{\ell=1}^3 \in \mathbb{C}^3$ with $b_1+b_2+b_3 \neq 0$ 
such that $f$ can be expressed as $f=f_{d,b}$, 
where $f_{d,b} \in \mathcal{Q}(C)$ is a unique map determined by the 
following properties. 
\begin{enumerate}
\item $\delta(f_{d,b})=d$. 
\item $p_{\ell}^{-} \approx [b_{\ell}:b_{\ell}^3:1] \in C^*$ for 
$I(f_{d,b}^{-1})=\{p_1^{-},p_2^{-},p_3^{-}\}$. 
\end{enumerate}
Moreover, the map $f_{d,b} \in \mathcal{Q}(C)$ satisfies the following. 
\begin{enumerate}
\item $c(f_{d,b}) = - \frac{1}{3} (b_1+b_2+b_3)\in \mathbb{C}^{\times}$. 
\item $p_{\ell}^{+} \approx [a_{\ell}:a_{\ell}^3:1] \in C^*$ for 
$I(f_{d,b})=\{p_1^{+},p_2^{+},p_3^{+}\}$, where 
$\displaystyle a_{\ell}:=\frac{1}{d} \big\{b_{\ell}-\frac{2}{3} 
(b_1+b_2+b_3) \big\}$. 
\end{enumerate}
\end{lemma}
\begin{remark} \label{rem:explquad}
The quadratic map 
$f_{d,b}([x:y:z])=[f_1:f_2:f_3]$ 
mentioned in Lemma \ref{lem:quad} is explicitly written as 
\begin{equation*}
\left\{
\begin{array}{cl}
f_1([x:y:z]) & =  (d/3) \cdot 
\big\{ (\nu_1^2-3 \nu_2) x^2 +\nu_1 \nu_3 z^2 -3 x y +2 \nu_1 y z 
-(\nu_1 \nu_2 - 3 \nu_3) z x \big\} \\[2mm]
f_2([x:y:z]) & =  (d/3)^3 \cdot 
\big\{ \nu_1 (\nu_1^3 -9 \nu_1 \nu_2 +27 \nu_3) x^2 -27 y^2 
+\nu_1^3 \nu_3 z^2 +9 (2 \nu_1^2 -3 \nu_2) x y \\[1mm]
~ & \phantom{=}~~~~~~~~~~~~ + (8 \nu_1^3 -27 \nu_1 \nu_2 +27 \nu_3) y z 
-\nu_1^2 (\nu_1 \nu_2 - 9 \nu_3) z x \big\} \\[2mm]
f_3([x:y:z]) & =  \nu_1 x^2 + \nu_3 z^2  - y z -\nu_2  z x, \\[2mm]
\end{array}
\right.
\end{equation*}
where $\nu_{\ell}=\nu_{\ell}(d,b)$ is given by 
\[\nu_1=a_1+a_2+a_3,\quad \nu_2=a_1 a_2 + a_2 a_3 +a_3 a_1, \quad 
\nu_3=a_1 a_2 a_3. \]
\end{remark}
In a similar manner, any linear map $f \in \mathcal{L}(C)$ is determined 
uniquely by the determinant $\delta(f)$ of $f$ (see \cite{D}). 
\begin{lemma} \label{lem:linear}
For any $d \in \mathbb{C}^{\times}$, there is a unique linear map 
$f \in \mathcal{L}(C)$ such that $\delta(f)=d$. In particular, the map 
$f \in \mathcal{L}(C)$ with $\delta(f)=1$ is the identity. Moreover, 
for any $f \in \mathcal{L}(C)$, 
the automorphism $f|_{C^*}$ restricted to $C^*$ is given by
\[
f|_{C^*} : [t:t^3:1] \mapsto [\delta(f) \cdot t : (\delta(f) \cdot t )^3:1]. 
\]
\end{lemma}
Next, let us consider the composition 
$f=f_n \circ f_{n-1} \circ \cdots \circ f_1 : \mathbb{P}^2 \to \mathbb{P}^2$ 
of quadratic birational maps 
$\overline{f}=(f_1,\dots,f_n) \in \mathcal{Q}(C)^n$. 
Put $I(f_i^{\pm 1})=\{ p_{i,1}^{\pm},p_{i,2}^{\pm},p_{i,3}^{\pm} \}$ and
\begin{equation} \label{eqn:Indp}
\check{p}_{i,j}^{+}:=f_1^{-1}|_{C} \circ \cdots \circ 
f_{i-1}^{-1}|_{C} (p_{i,j}^{+}), \qquad 
\check{p}_{i,j}^{-}:=f_n|_{C} \circ \cdots \circ f_{i+1}|_{C} 
(p_{i,j}^{-}) 
\end{equation} 
(see Figure \ref{fig:comp}). Then it is easy to see that 
$I(f^{\pm 1}) \subset \{\check{p}_{i,j}^{\pm} \, | \, 
(i,j) \in \mathcal{K}(n) \}$. Moreover, let $\delta(\overline{f})$ be the 
{\sl determinant} of $\overline{f}$ defined by 
$\delta(\overline{f})=\prod_{i=1}^n \delta(f_i)$ or, in other words, 
$\delta(\overline{f})=\delta(f)$. 
\begin{proposition} \label{pro:birat} 
Let $\overline{f}=(f_1,\dots,f_n) \in \mathcal{Q}(C)^n$ be an $n$-tuple of 
quadratic birational maps in $\mathcal{Q}(C)$ with 
$d=\delta(\overline{f}) \neq 1$, and let 
$\check{p}_{i,j}^{\pm}$ be the points 
given in (\ref{eqn:Indp}). Then 
there is a unique pair $(v,s)$ of values 
$v=(v_{\iota})_{\iota \in 
\mathcal{K}(n)} \in \mathbb{C}^{3n}$ 
and $s=(s_i)_{i=1}^n \in (\mathbb{C}^{\times})^n$ satisfying 
\begin{equation} \label{eqn:BS}
v_{i,1}+v_{i,2}+v_{i,3}   =  - \sum_{k=1}^{i-1} s_k 
+(d-2) \cdot s_i - d \sum_{k=i+1}^n s_k, \quad (1 \le i \le n), 
\end{equation}
such that 
the composition $f=f_n \circ \cdots \circ f_1$ satisfies 
\begin{enumerate}
\item $f|_{C^*} : C^* \ni [t+ \frac{1}{3} c_s:
(t+ \frac{1}{3} c_s)^3:1] \mapsto [d \cdot t+ \frac{1}{3} c_s
:(d \cdot t+ \frac{1}{3}c_s)^3:1] \in C^*$ with 
\begin{equation} \label{eqn:fix}
~~~~~~~~~~~~~~~~~~~~~~~~~~~
c_s := \sum_{k=1}^{n} s_k, 
~~~~~~~~~~~~~~~~~~~~~~~~~~~
\end{equation}
\item 
$\check{p}_{i,j}^{-} \approx [v_{i,j}+\frac{1}{3} c_s:
(v_{i,j}+ \frac{1}{3} c_s)^3:1] \in C^*$, 
\item 
$\check{p}_{i,j}^{+} \approx [u_{i,j}+\frac{1}{3} c_s:
(u_{i,j}+\frac{1}{3} c_s)^3:1] \in C^*$, 
where 
\begin{equation} \label{eqn:uv}
u_{i,j}:= \frac{1}{d} \big\{ v_{i,j} 
- (d-1) \cdot s_i \big\}. 
\end{equation}
\end{enumerate}
Conversely, for any $d \in \mathbb{C} \setminus \{ 0,1\}$, 
$v \in \mathbb{C}^{3n}$ and $s \in (\mathbb{C}^{\times})^n$ 
satisfying equation 
(\ref{eqn:BS}), there exists an $n$-tuple 
$\overline{f}=(f_1,\dots,f_n) \in \mathcal{Q}(C)^n$ such that 
the composition $f=f_n \circ \cdots \circ f_1$ satisfies  
$\delta(f)=d$ and conditions (1)--(3). 
Moreover, $\overline{f}$ is 
determined uniquely by $(d,v,s)$ in the sense that if 
$\overline{f}=(f_1,\dots,f_n)$ and $\overline{f}'=(f_1',\dots,f_n')$ 
are determined by $(d,v,s)$, then there are linear maps 
$g_1, \dots, g_{n-1} \in \mathcal{L}(C)$ such that 
the following diagram commutes: 
\begin{equation*} 
\begin{CD}
\mathbb{P}_0^2 @> f_1' >> \mathbb{P}_1^2 @> f_2' >>  
\cdots @> f_{n-1}' >> \mathbb{P}_{n-1}^2 @> f_n' >> \mathbb{P}_n^2 \\
@|  @V g_1 VV @. @V g_{n-1} VV @| \\
\mathbb{P}_0^2 @> f_1 >> \mathbb{P}_1^2 @> f_2 >>  
\cdots @> f_{n-1} >> \mathbb{P}_{n-1}^2 @> f_n >> \mathbb{P}_n^2. 
\end{CD}
\end{equation*}
\end{proposition}
{\it Proof}.
From Lemma \ref{lem:quad}, each map $f_i \in \mathcal{Q}(C)$ is given by 
$f_i=f_{d_i,(b_{i,j})}$ for some $d_i \in \mathbb{C}^{\times}$ and  
$(b_{i,j})_{j=1}^{3} \in \mathbb{C}^3$ with 
$b_{i}:=b_{i,1}+b_{i,2}+b_{i,3} \neq 0$. 
Then the maps $f_i$ and $f$ restricted to $C^*$ can be 
expressed as $f_i |_{C^*}( [t:t^3:1]) = [y_i(t):y_i(t)^3:1]$ 
and $f |_{C^*}([t:t^3:1]) = [y(t):y(t)^3:1]$, respectively, 
where $y_i, \, y : \mathbb{C} \to \mathbb{C}$ are given by 
\[
y_i(t)=d_i \cdot t -\frac{1}{3} b_{i}, 
\]
and $y:=y_n \circ y_{n-1} \circ \cdots \circ y_1$. 
Now we put 
$\check{d}_i:=d_{i+1} \cdot d_{i+2} \cdots d_n$, 
$a_{i,j}:=( b_{i,j}-\frac{2}{3} b_i)/d_i$ and 
\[
\check{a}_{i,j}:=y_1^{-1} \circ \cdots \circ y_{i-1}^{-1} (a_{i,j}), 
\qquad 
\check{b}_{i,j}:= y_n \circ \cdots \circ y_{i+1} (b_{i,j}), \qquad 
\check{b}_{i}:= \check{b}_{i,1}+\check{b}_{i,2}+\check{b}_{i,3}. 
\]
Then it follows that 
$\check{p}_{i,j}^{+} \approx 
[\check{a}_{i,j}:\check{a}_{i,j}^3:1]$,  
$\check{p}_{i,j}^{-} \approx 
[\check{b}_{i,j}:\check{b}_{i,j}^3:1]$ 
and $d=\check{d}_0$. A little calculation shows that 
\begin{equation*}
\begin{array}{rl}
y(t) & \displaystyle = 
d \cdot t- \frac{1}{3} \sum_{k=1}^n 2^{k-1} \cdot \check{b}_k, \\[2mm]
\check{b}_{i,j} & = \displaystyle
\check{d}_{i} \cdot b_{i,j} -\frac{1}{3} \sum_{k=i+1}^n 
\check{d}_{k} \cdot b_k = \check{d}_{i} \cdot b_{i,j} 
-\frac{d-1}{3} \sum_{k=i+1}^n s_k, \\[2mm]
\check{a}_{i,j} & \displaystyle = 
\frac{1}{d} \Bigl( \check{d}_i \cdot b_{i,j} - \check{d}_i \cdot b_i 
+ \frac{1}{3} \sum_{k=1}^i \check{d}_k \cdot b_k \Bigr) 
=\frac{1}{d} \Big\{ \check{d}_i \cdot b_{i,j} - (d-1) \cdot s_i 
+ \frac{d-1}{3} \sum_{k=1}^i s_k \Big\}, \\[2mm] 
\end{array}
\end{equation*}
where $s_i := \check{d}_i \cdot b_i/(d-1) \neq0$. If we put 
\[
v_{i,j} := 
\check{d}_i \cdot b_{i,j} - \frac{1}{3} 
\Big( \sum_{k=1}^i s_k + d \cdot \sum_{k=i+1}^n s_k \Big), 
\]
then we have 
\[
v_{i,1}+v_{i,2}+v_{i,3} = \displaystyle 
\check{d}_i \cdot b_{i} - 
\Big( \sum_{k=1}^i s_k + d \cdot \sum_{k=i+1}^n 
s_k \Big) = \displaystyle  - \sum_{k=1}^{i-1} s_k 
+(d-2) s_i - d \sum_{k=i+1}^n s_k, 
\]
which shows that equation (\ref{eqn:BS}) holds. Moreover, since 
$\check{b}_i = \check{d}_{i} \cdot b_{i} - (d-1) \cdot \sum_{k=i+1}^n 
s_k = (d-1) \cdot \{s_{i}-\sum_{k=i+1}^n s_{k} \}$ 
and thus 
$\sum_{k=1}^n 2^{k-1} \cdot \check{b}_k = (d-1) \cdot c_s$, the map 
$y(t)= d \cdot t - (d-1) \cdot c_s/3$ 
has the unique fixed point $c_s/3$ under the assumption that $d \neq 1$. 
Finally, we have 
\[
\begin{array}{rl}
\check{b}_{i,j} = & \displaystyle 
\check{d}_{i} \cdot b_{i,j} 
-\frac{d-1}{3} \hspace{-1.5mm} \sum_{k=i+1}^n \hspace{-0.5mm} s_k = 
v_{i,j} + \frac{1}{3} 
\Big( \sum_{k=1}^i \hspace{-0.5mm} s_k + d \cdot \hspace{-1mm} 
\sum_{k=i+1}^n \hspace{-0.5mm} s_k \Big) 
- \frac{d-1}{3} \hspace{-1.5mm} \sum_{k=i+1}^n \hspace{-0.5mm} 
s_k = v_{i,j} +  \frac{1}{3} c_s, \\[2mm] 
d \cdot \check{a}_{i,j} = & 
\displaystyle \check{d}_{i} \cdot b_{i,j} 
- (d-1) \cdot s_i + \frac{d-1}{3} \sum_{k=1}^i s_k \\[2mm]
= & \displaystyle v_{i,j}  + \frac{1}{3} 
\Big( \sum_{k=1}^i s_k + d \cdot \sum_{k=i+1}^n s_k \Big) 
- (d-1) \cdot s_i + \frac{d-1}{3} \sum_{k=1}^i s_k \\[2mm]
= & \displaystyle v_{i,j} - (d-1) \cdot s_i  +  \frac{d}{3} c_s. 
\end{array}
\]
Thus conditions (1)--(3) hold. 
\par 
Conversely, for any $d \neq 1$, 
$(s_{i})$ and $(v_{i,j})$ satisfying (\ref{eqn:BS}), 
the maps $(f_i)=(f_{d_i,(b_{i,j})})$ with 
\[
\begin{array}{l}
d_1 \cdots d_n =d, \\[2mm] 
\displaystyle b_{i,j} = \frac{1}{d_{i+1} \cdots d_{n}} \Big\{ v_{i,j} + 
\frac{1}{3} \Big(\sum_{k=1}^i s_k + d \cdot \sum_{k=i+1}^n s_k \Big) \Big\} 
\end{array}
\]
give the birational map $f=f_n \circ \cdots \circ f_1$ satisfying 
$\delta(f)=d$ and conditions (1)--(3). 
Moreover, assume that there are two $n$-tuples 
$\overline{f}=(f_1, \dots,f_n)$ and $\overline{f}'=(f_1', \dots,f_n')$ in 
$\mathcal{Q}(C)^n$ such that 
$f=f_n \circ \cdots \circ f_1$ and $f'=f_n' \circ \cdots \circ f_1'$ 
satisfy $\delta(f)=\delta(f')=d$ and conditions (1)--(3) for $(d,v,s)$. Put 
$g_i:=f_{i+1}^{-1} \circ \cdots \circ f_{n}^{-1} \circ 
f_n' \circ \cdots \circ f_{i+1}' : \mathbb{P}_i^2 \to \mathbb{P}_i^2$. 
Then one has $f_i'=g_{i}^{-1} \circ f_i \circ g_{i-1}$, where 
$g_n=\text{id}$. It follows from 
condition (2) that $I(f_n^{-1})=I((f_n')^{-1})$, which means that 
$g_{n-1}=f_n^{-1} \circ f_n'$ is an automorphism of $\mathbb{P}_{n-1}^2$ 
preserving $C$, and thus $g_{n-1} \in \mathcal{L}(C)$. 
In a similar manner, under the assumption that $g_{i} \in \mathcal{L}(C)$ for 
some $1 \le i \le n$, condition (2) shows that $I(f_i^{-1})=g_iI((f_i')^{-1})$ 
and hence $g_{i-1}=f_i^{-1} \circ g_{i} \circ f_i' \in \mathcal{L}(C)$. 
Moreover, it follows from condition (1) that the determinant of $g_0$ is given by 
$\delta(g_0)=\delta(f') \cdot \delta(f)^{-1}=1$, which means that 
$g_0=\text{id}$ (see Lemma \ref{lem:linear}). 
This completes the proof. 
\hfill $\Box$ \par\medskip 
\begin{remark} \label{rem:checkp}
From the definition of $\check{p}_{\iota}^{\pm}$ given in (\ref{eqn:Indp}), 
we have the following relations: 
\begin{enumerate}
\item $f|_{C}^{k}(\check{p}_{\iota}^-) \approx \check{p}_{\iota'}^-$, 
or equivalently $d^{k} \cdot v_{\iota}=v_{\iota'}$,  
if and only if $p_{\iota}^m \approx p_{\iota'}^-$, 
where $m=\theta_{i,i'}(k)$ with $\theta_{i,i'}(k)$ given in 
(\ref{eqn:vartheta}). 
\item $f|_{C}^{k}(\check{p}_{\iota}^-) \approx \check{p}_{\iota'}^+$, 
or equivalently $d^{k} \cdot v_{\iota}=u_{\iota'}$, 
if and only if $p_{\iota}^m \approx p_{\iota'}^+$, 
where $m=\theta_{i,i'-1}(k+1)$. 
In particular, it follows from (\ref{eqn:mu}) that 
$p_{\iota}^{\mu(\iota)} \approx p_{\sigma(\iota)}^+$ if and only if 
$f|_{C}^{\kappa(\iota)-1}(\check{p}_{\iota}^-) \approx \check{p}_{\sigma(\iota)}^+$. 
\end{enumerate}
Indeed, for example, assertion (1) can be established from the relation 
\[
\begin{array}{rl}
p_{\iota'}^- = & f_{i'+1}^{-1}|_{C} \circ \cdots \circ 
f_{n}^{-1}|_{C} (\check{p}_{\iota'}^-) \approx 
f_{i'+1}^{-1}|_{C} \circ \cdots \circ f_{n}^{-1}|_{C} 
( f|_{C}^k(\check{p}_{\iota}^-)) \\[2mm] 
= & f_{i'+1}^{-1}|_{C} \circ \cdots \circ f_{n}^{-1}|_{C} \circ
( f|_{C})^k \circ f_n|_C \circ \cdots \circ 
f_{i+1}|_{C} (p_{\iota}^-) = 
p_{\iota}^{\theta_{i,i'}(k)}. 
\end{array}
\]
\end{remark}
Assume that there is a tentative realization $\overline{f}$ of 
$\tau$. Then, the relation 
$p_{\iota}^{\mu(\iota)} \approx p_{\sigma(\iota)}^+$ yields 
$f|_C^{\kappa(\iota)-1}(\check{p}_{\iota}^-)\approx \check{p}_{\sigma(\iota)}^+$ 
and thus $d^{\kappa(\iota)-1} \cdot v_{\iota}=u_{\sigma(\iota)}=u_{\iota_1}$ 
for any $\iota \in \mathcal{K}(n)$. 
Therefore, in view of (\ref{eqn:uv}), the pair 
$(v,s) \in \mathbb{C}^{3n} \times (\mathbb{C}^{\times})^{n}$ satisfies 
\begin{equation} \label{eqn:expab}
v_{\iota_1}  =  d^{\kappa(\iota)} \cdot 
v_{\iota} +(d-1) \cdot s_{i_1} 
\qquad (\iota =(i,j) \in \mathcal{K}(n)),  
\end{equation}
which is equivalent, when $d$ is not a root of unity, 
to the expression 
\begin{equation} \label{eqn:expre} 
v_{\iota} = v_{\iota}(d) = - 
\frac{ d^{ \varepsilon_{|\iota|}} 
\cdot (d-1)}{d^{ \varepsilon_{|\iota|}}-1} 
\bigl( d^{-\varepsilon_{1}} \cdot s_{i_1} + 
d^{-\varepsilon_{2}} \cdot s_{i_2} + \cdots + 
d^{-\varepsilon_{|\iota|}} \cdot s_{i_{|\iota|}} \bigr), 
\end{equation}
where $|\iota|:=\# \{ \iota_k \, | \, k \ge 0 \}$ and 
$\varepsilon_{\ell}:=\varepsilon_{\ell}(\iota)=\sum_{k=0}^{\ell-1} 
\kappa(\iota_k)$. 
Conversely, if there is a pair 
$(d,v,s) \in (\mathbb{C} \setminus \{0,1\}) \times \mathbb{C}^{3n} 
\times (\mathbb{C}^{\times})^{n}$ satisfying 
(\ref{eqn:BS}) and (\ref{eqn:expab}), then, by virtue of 
Proposition \ref{pro:birat}, 
there is a tentative realization $\overline{f} \in \mathcal{Q}(C)$ 
of $\tau$ with $\delta(\overline{f})=d$. Therefore, it is important 
to consider whether the system of equations 
\begin{equation} \label{eqn:system} 
\left\{
\begin{array}{rcl}
v_{i,1}+v_{i,2}+v_{i,3}  & = & - \sum_{k=1}^{i-1} s_k 
+(d-2) \cdot s_i - d \sum_{k=i+1}^n s_k, \quad (1 \le i \le n)  \\[2mm]
v_{\iota_1} & = & d^{\kappa(\iota)} \cdot 
v_{\iota} +(d-1) 
\cdot s_{i_1} \qquad (\iota \in \mathcal{K}(n)) 
\end{array}
\right.
\end{equation}
for $(d,v,s) \in (\mathbb{C} \setminus \{0,1\}) \times \mathbb{C}^{3n} 
\times (\mathbb{C}^{\times})^{n}$ can admit solutions, which are closely 
related to the eigenvalue problem of the Weyl group element $w_{\tau}$. 
Namely, as is mentioned in the following proposition, 
solutions of (\ref{eqn:system}) will appear 
in the coefficients of eigenvectors of $w_{\tau}$. 
\begin{proposition} \label{pro:eigen}
Let $\tau$ be an orbit data, and $d$ be a complex number different from 
$0$ and $1$. Then, a vector $y \neq 0$ in $\mathbb{Z}^{\tau} 
\otimes_{\mathbb{Z}} \mathbb{C}$ expressed as 
\begin{equation*} 
y = v_0 \cdot e_0 +
\sum v_{\iota}^{k} \cdot e_{\iota}^k \in 
\mathbb{Z}^{\tau} \otimes_{\mathbb{Z}} \mathbb{C} 
\end{equation*}
is an eigenvector of $w_{\tau}$ corresponding to the eigenvalue $d$ 
if and only if there is a pair 
$(v,s) \neq (0,0) \in \mathbb{C}^{3n} \times \mathbb{C}^{n}$ 
satisfying equations (\ref{eqn:system}) 
such that the following conditions hold: 
\begin{enumerate}
\item $v_{\iota}^k=d^{k-1} \cdot v_{\iota}$ for any 
$\iota \in \mathcal{K}(n)$ and 
$1 \le k \le \kappa(\iota)$. 
\item $ v_0 = c_s$, where $c_s$ is given in (\ref{eqn:fix}).  
\end{enumerate}
Moreover, for the eigenvector $y$ corresponding to $d$, a pair 
$(v,s) \neq (0,0) \in \mathbb{C}^{3n} \times \mathbb{C}^{n}$ 
satisfying equations (\ref{eqn:system}) and conditions 
(1)--(2) is uniquely determined. 
\end{proposition}
{\it Proof}. Assume that $y$ is an eigenvector corresponding to $d$. 
It is easy to see that the coefficient of $e_{\iota}^{k}$ in $w_{\tau}(y)$ 
is $v_{\iota}^{k+1}$ for any $1 \le k \le \kappa(\iota)-1$. 
Hence, one has $v_{\iota}^{k+1}= d \cdot v_{\iota}^k$, 
or $v_{\iota}^k=d^{k-1} \cdot v_{\iota}^1$. 
Moreover, we determine $(v,s) \in \mathbb{C}^{3n} \times \mathbb{C}^n$ 
as follows. 
Put $v_{n,i}=v_{n,i}^1$ for $i \in \{1,2,3\}$, and 
$s_n=(v_0+v_n)/(d-1)$, where $v_n:=v_{n,1}+v_{n,2}+v_{n,3}$. 
For $1 \le \ell \le n-1$, assume that $v_{j,i}$ 
and $s_{j}$ with $j \ge \ell +1$ are already determined. Then, put 
\[
v_{\ell,i}=\left\{
\begin{array}{ll}
v_{\ell,i}^1 \quad & (\kappa(\ell,i) \ge 1) \\[2mm]
v_{(\ell,i)_1}-(d-1) \cdot s_{\ell_1} \quad & 
(\kappa(\ell,i) =0), 
\end{array}
\right.
\]
where $\sigma(\ell,i)=(\ell,i)_1=(\ell_1,i_1)$ 
(note that $\ell_1\ge \ell+1$ if $\kappa(\ell,i) =0$), 
and 
\[
s_{\ell}=\left\{
\begin{array}{ll}
(v_0+v_{\ell})/(d-1)+\sum_{k=\ell+1}^n s_k \quad & (\ell \ge 2) \\[2mm]
v_0-\sum_{k=2}^n s_k \quad & 
(\ell =1), 
\end{array}
\right.
\]
where $v_{\ell}:=v_{\ell,1}+v_{\ell,2}+v_{\ell,3}$. At this stage, 
it is easily checked that 
$(v,s)$ satisfies conditions (1)--(2), equation 
(\ref{eqn:BS}) for $i \ge 2$ and equation (\ref{eqn:expab}) for 
$\kappa(\iota)=0$. 
Now we claim that 
the following relation holds 
for $1 \le \ell \le n+1$: 
\begin{equation} \label{eqn:wexp}
q_\ell \circ \cdots \circ q_n (y) = 
v^\ell \cdot e_0 + \hspace{-3mm} 
\sum_{\substack{\ i < \ell \\ \kappa(\iota) \ge 1}} \hspace{-2mm} 
v_{\iota} \cdot e_{\iota}^1 
+ \hspace{-3mm} 
\sum_{\substack{\ i < \ell \le 
i_1 \\ \kappa(\iota)=0}} \hspace{-2mm} 
v_{\iota} \cdot e_{\bar{\sigma}(\iota)}^1 
+ \hspace{-3mm} 
\sum_{\substack{\ \ell \le i \\ \kappa(\sigma^{-1}(\iota)) \ge 1}} 
\hspace{-3mm} \Big\{ v_{\iota} - (d-1) \cdot s_{i} \Big\} \cdot 
e_{\bar{\sigma}(\iota)}^1 
+ \sum_{k \ge 2} d^{k-1} \cdot v_{\iota} 
\cdot e_{\iota}^k, 
\end{equation}
where $v^\ell:=\sum_{k=1}^{\ell-1} s_k +d \sum_{k=\ell}^{n} s_k$. 
Indeed, if $\ell=n+1$, the relation is trivial. Assume that the 
relation holds when $\ell+1 \ge 2$. Then the automorphism $q_\ell$ changes only 
the coefficients $v^{\ell+1}$ and $v_{\ell,i}$ 
in $q_{\ell+1} \circ \cdots \circ q_n (y)$ as follows: 
\begin{equation*}
\displaystyle q_\ell \Big( v^{\ell+1} \cdot e_0 + 
\sum_{i=1}^3 v_{\ell,i} \cdot e_{\bar{\sigma}(\ell,i)}^1 \Big) = 
\big( 2 v^{\ell+1} + v_\ell \big) \cdot e_0 
- \sum_{i=1}^3 \big(v^{\ell+1} + v_{\ell,j} + 
v_{\ell,k} \big) \cdot e_{\bar{\sigma}(\ell,i)}^1,  
\end{equation*}
where $\{i,j,k\}=\{1,2,3\}$. Therefore, 
when $\ell \ge 2$, equation (\ref{eqn:wexp}) holds from 
the facts that $2 v^{\ell+1} + v_\ell= v^{\ell}$, 
$v^{\ell+1} + v_{\ell,j} + v_{\ell,k}=
(d -1) \cdot s_\ell- v_{\ell,i}$, and that 
$\big\{ v_{\ell,i} - (d-1) \cdot s_{\ell} \big\} \cdot e_{\bar{\sigma}(\ell,i)}^1= 
v_{\sigma^{-1}(\ell,i)} \cdot e_{\bar{\sigma}(\sigma^{-1}(\ell,i))}^1$ if 
$\kappa(\sigma^{-1}(\ell,i))=0$. 
Moreover, since the coefficient of $e_0$ in $w_{\tau} (y)$ is 
$d \cdot v_0$ and $e_0$ is fixed by $r_{\tau}$, the coefficient of 
$e_0$ in $q_1 \circ \cdots \circ q_n (y)$ is expressed 
as $2 v^2 + v_1=d \cdot v_0=v^1$, which 
yields $v_{1}=(d-2) \cdot s_1 - d \sum_{k=2}^n s_k$. Thus, the coefficient 
of $e_{1,i}$ in $q_1 \circ \cdots \circ q_n (y)$ 
is given by $v_{1,i}-(d-1) \cdot s_1$ 
and equation (\ref{eqn:wexp}) holds when $\ell=1$. 
The claim follows from these observations. In particular, 
$(v,s)$ satisfies equation (\ref{eqn:BS}) for $i =1$. 
\par
By the above claim, we have 
\[
q_1 \circ \cdots \circ q_n (y) = 
v^1 \cdot e_0 + \sum_{\kappa(\iota) \ge 1} 
\Big\{ v_{\iota_1} - (d-1) \cdot s_{i_1} \Big\} \cdot 
e_{\bar{\sigma}(\iota_1)}^1 
+ \sum_{k \ge 2} d^{k-1} \cdot v_{\iota}
\cdot e_{\iota}^k. 
\]
Thus, the coefficient of $e_{0}$ in $w_{\tau}(y)$ is 
$v^1= d \cdot c_s$. 
Similarly, the coefficient of 
$e_{\iota}^{\kappa(\iota)}$ in 
$w_{\tau}(y)$ is given by 
$v_{\iota_1} - (d-1) \cdot s_{i_1}$. 
This means that $v_{\iota_1} - (d-1) \cdot s_{i_1} = 
d \cdot v_{\iota}^{\kappa(\iota)}=
d^{\kappa(\iota)} \cdot v_{\iota}$ and that 
$(v,s)$ satisfies equation (\ref{eqn:expab}) for $\kappa(\iota) \ge 1$. 
It should be noted that $v \neq 0$, 
since if $v = 0$ then one has $s=0$ from (\ref{eqn:expab}) 
and so $y=0$, which is a contradiction. 
\par
Moreover, the above argument shows that, for the eigenvector $y$, 
$(v,s)$ is determined uniquely by equation (\ref{eqn:BS}) for $i \ge 2$, 
equation (\ref{eqn:expab}) for $\kappa(\iota)=0$, condition (1) for $k=1$ 
and condition (2), which gives the uniqueness of $(v,s)$. 
\par
Conversely, we can easily check that, for a pair 
$(v,s) \neq (0,0) \in \mathbb{C}^{3n} \times \mathbb{C}^{n}$ 
satisfying equations (\ref{eqn:system}), the vector $y$ given by  
conditions (1)--(2) is an eigenvector of $w_{\tau}$ corresponding to $d$. 
The proof is complete.
\hfill $\Box$ \par\medskip 
Let $w : \mathbb{Z}^{1,N} \to \mathbb{Z}^{1,N}$ be 
a general lattice automorphism in $W_N$. 
It is known that the characteristic polynomial $\chi_{w}(t)$ of $w$ 
can be expressed as 
\[
\chi_{w}(t) = 
\left\{
\begin{array}{ll}
R_w(t) \quad & (\lambda(w)=1) \\[2mm]
R_w(t) S_w(t) \quad & (\lambda(w)>1), 
\end{array}
\right.
\]
where $R_w(t)$ is a product of cyclotomic polynomials, and $S_w(t)$ 
is a Salem polynomial, namely, the minimal polynomial of a Salem number. 
Here, a {\sl Salem number} is an algebraic integer $\delta > 1$ 
such that its conjugates include $\delta^{-1} < 1$ and the conjugates 
other than $\delta^{\pm 1}$ lie on the unit circle. 
Therefore, if $w$ satisfies $\lambda(w) >1$ and $d$ satisfies 
$S_{w}(d)=0$, then there is a unique eigenvector, up to constant multiple, 
corresponding to $d$. Moreover, an eigenvalue $d$ with $|d| >1$ is unique 
and is a Salem number $d = \lambda(w) > 1$. 
\par
To simplify the notation, we put $S_{\tau}(t):=S_{w_{\tau}}(t)$ and 
$\lambda(\tau):=\lambda(w_{\tau})$. Then, 
a corollary of Proposition \ref{pro:eigen} can be established. 
\begin{corollary} \label{cor:sol}
Assume that $d$ is not a root of unity and there is a solution 
$(v,s) \neq (0,0) \in \mathbb{C}^{3n} \times \mathbb{C}^{n}$ of equations 
(\ref{eqn:system}). Then, $d$ is a root of $S_{\tau}(t)=0$, 
and $v$ and $s$ are nonzero. 
Conversely, if $d$ is a root of $S_{\tau}(t)=0$, 
then there is a unique solution 
$(v,s) \in (\mathbb{C}^{3n} \setminus \{0\}) \times 
(\mathbb{C}^{n} \setminus \{0\})$ 
of equations (\ref{eqn:system}), up to a constant multiple. 
\end{corollary}
{\it Proof}. 
First assume that $d$ is not a root of unity and there is a solution 
$(v,s) \neq (0,0) \in \mathbb{C}^{3n} \times \mathbb{C}^{n}$ of (\ref{eqn:system}). 
Then $d$ is a root of $S_{\tau}(t)=0$ from 
Proposition \ref{pro:eigen}. Moreover $v$ and $s$ are nonzero. 
Indeed, if $v=0$ then $s=0$ from (\ref{eqn:expab}), and if $s=0$ then $v=0$ 
from (\ref{eqn:expre}). Conversely, if 
$d$ is a root of $S_{\tau}(t)=0$, then a solution 
$(v,s) \neq (0,0) \in \mathbb{C}^{3n} \times \mathbb{C}^{n}$ of (\ref{eqn:system}), 
which satisfies $v \neq 0$ and $s \neq 0$ from the above argument, 
is unique, as an eigenvector corresponding to $d$ is unique. 
\hfill $\Box$ \par\medskip 
For a root $d$ of $S_{\tau}(t)=0$, let 
$(v,s) \in (\mathbb{C}^{3n} \setminus \{0\}) \times 
(\mathbb{C}^{n} \setminus \{0\})$ be a solution of (\ref{eqn:system}) as is 
mentioned in Corollary \ref{cor:sol}. If $s$ satisfies $s_{\ell} \neq 0$ for any 
$\ell =1,\dots,n$, then there is a tentative realization 
$\overline{f} \in \mathcal{Q}(C)^n$ of $\tau$ such that 
$\delta(\overline{f})=d$ from the above argument. Moreover, 
the composition $f=f_n \circ \cdots \circ f_1$ 
is unique up to conjugacy by a linear map in $\mathcal{L}(C)$, 
as a solution of (\ref{eqn:system}) is unique up to a constant multiple. 
Summing up these discussions, we have the following proposition. 
\begin{proposition} \label{pro:tenta0}
Let $\tau$ be an orbit data with $\lambda(\tau)>1$, $d$ be 
a root of $S_{\tau}(t)=0$ and $s \neq 0$ be the unique solution 
of equations (\ref{eqn:system}) (see Corollary \ref{cor:sol}). 
Then $s$ satisfies $s_\ell \neq 0$ for any $1 \le \ell \le n$ 
if and only if there is a tentative realization $\overline{f}$ of $\tau$ 
such that $\delta(\overline{f})=d$. Moreover, the tentative realization 
$\overline{f}$ of $\tau$ is uniquely determined in the sense that 
if there are two tentative realizations 
$\overline{f}=(f_1,\dots,f_n)$ and $\overline{f}'=(f_1',\dots,f_n')$ of $\tau$ 
such that $\delta(\overline{f})=\delta(\overline{f}')=d$, 
then there are linear maps 
$g_1, \dots, g_n \in \mathcal{L}(C)$ such that 
the following diagram commutes: 
\begin{equation*} 
\begin{CD}
\mathbb{P}_0^2 @> f_1' >> \mathbb{P}_1^2 @> f_2' >>  
\cdots @> f_{n-1}' >> \mathbb{P}_{n-1}^2 @> f_n' >> \mathbb{P}_n^2 \\
@V g_0 VV  @V g_1 VV @. @V g_{n-1} VV @V g_n VV \\
\mathbb{P}_0^2 @> f_1 >> \mathbb{P}_1^2 @> f_2 >>  
\cdots @> f_{n-1} >> \mathbb{P}_{n-1}^2 @> f_n >> \mathbb{P}_n^2. 
\end{CD}
\end{equation*}
where $g_0:=g_n$. 
\end{proposition}
\begin{remark} \label{rem:inf2}
As is seen in Proposition \ref{pro:tenta0}, the tentative realization 
$\overline{f}$ of $\tau$ with $\delta(\overline{f})=d$ 
is unique. 
However, when $p_{\iota}^- \approx p_{\iota'}^-$ 
for some 
$\iota \neq \iota' \in \{(i,1),(i,2),(i,3)\}$, 
there remains an ambiguity about how to label indeterminacy 
points, namely, about a choice between  
$p_{\iota}^{\pm} < p_{\iota'}^{\pm}$ and  
$p_{\iota}^{\pm} > p_{\iota'}^{\pm}$ (see also Remark \ref{rem:labeling}). 
\end{remark}
Proposition \ref{pro:tenta0} raises a question as to whether, 
for a given orbit data $\tau$, the solution $s$ of (\ref{eqn:system}) 
satisfies $s_\ell \neq 0$ for any $\ell=1,\dots,n$. 
This question can be solved in terms of the absence of periodic roots. 
To see this, we need some preliminaries. 
\par
Let $w \in W_N$ be a general Weyl group element with $\lambda(w) >1$. 
Then there is a direct sum decomposition of the real vector space: 
\[
~~~~~~~~~~~~~~~~~~~~~~~~~~~
\mathbb{R}^{1,N} :=\mathbb{Z}^{1,N} \otimes_{\mathbb{Z}} \mathbb{R}
= V_w \oplus V_w^c,
~~~~~~~~~~~~~~~~~~~~~~~~~~~
\]
such that the decomposition is preserved by $w$, and $S_w(t)$ and $R_w(t)$ 
are the characteristic polynomials of $w|_{V_w}$ and $w|_{V_w^c}$ 
, respectively. 
We notice that $V_w^c$ is the orthogonal complement of $V_w$ with respect to 
the Lorentz inner product. Moreover, let $\ell_{w}$ be the 
the minimal positive integer satisfying $d^{\ell_{w}}=1$ for any 
root $d$ of the equation $R_{w}(t)=0$. 
Then we have the following lemma. 
\begin{lemma} \label{lem:per}
Assume that $\delta= \lambda(w) > 1$, and let $d$ be an 
eigenvalue of $w$ that is not a root of unity. 
Then, for a vector $z \in \mathbb{Z}^{1,N}$, the following are equivalent. 
\begin{enumerate}
\item $(z, y_d)=0$, where $y_d$ 
is the eigenvector of $w$ corresponding to the eigenvalue $d$. 
\item $z \in V_w^c \cap \mathbb{Z}^{1,N}$.
\item $z$ is a periodic vector of $w$ with period $\ell_{w}$. 
\item $z$ is a periodic vector of $w$ with some period $k$. 
\end{enumerate}
\end{lemma}
{\it Proof}. $(1) \Longrightarrow (2)$. First, we notice that 
$y_d$ can be chosen so that 
$y_d \in \mathbb{Z}^{1,N} \otimes_{\mathbb{Z}} \mathbb{Z}[d]$. 
The coefficient of $e_i$ in $y_d$, and thus 
that in $y_{d'}$ for any conjugate $d'$, are expressed as 
$(y_d)_i=\upsilon_i(d)$ and 
$(y_{d'})_i=\upsilon_i(d')$ 
for some $\upsilon_i(x) \in \mathbb{Z}[x]$. 
Since $z=(z_i) \in \mathbb{Z}^{1,N}$ and so 
$(z,y_x)=z_0 \cdot \upsilon_0(x) - \sum_{i \neq 0} z_i \cdot \upsilon_i(x) \in \mathbb{Z}[x]$, 
we have $(z, y_{d'})=0$ from the relation $(z, y_d)=0$. 
Thus it follows that $z \in V_w^c \cap \mathbb{Z}^{1,N}$. 
\par 
$(2) \Longrightarrow (3)$. For any eigenvalues $d, d'$, we have 
\[
(y_d,y_{d'})=(w(y_d),w(y_{d'}))
=d \cdot d' \cdot (y_d,y_{d'}), 
\]
which means that 
$(y_d,y_{d'})=0$ if $d \cdot d' \neq 1$. 
In particular, one has 
$(y_{\delta},y_{\delta})=(y_{1/\delta},y_{1/\delta})=0$. 
Moreover, since 
$y_{\delta}, y_{1/\delta} \in \mathbb{R}^{1,N}$ 
are linearly independent over $\mathbb{R}$, 
$(y_{\delta}, y_{1/\delta})$ is nonzero, and thus either 
$(y_{\delta}+y_{1/\delta},y_{\delta}+y_{1/\delta})$ or 
$(y_{\delta}-y_{1/\delta},y_{\delta}-y_{1/\delta})$ is 
positive. As $\mathbb{R}^{1,N}$ has signature $(1,N)$ and 
$V_w$ has signature $(1,s)$ for some $s \ge 1$, 
$V_w^c$ is negative definite. This shows that $w|_{V_w^c}$ has 
finite order. Since any eigenvalue $d$ of $w|_{V_w^c}$ satisfies 
$d^{\ell_w}=1$, we have $w^{\ell_w}(z)=z$. 
\par 
$(4) \Longrightarrow (2)$. Assume that $w^{k}(z)=z$ for some $k \ge 1$. We 
express $z$ as $z=z' + z''$ for some $z' \in V_w$ and $z'' \in V_w^c$, and 
then express $z'$ as $z'=\sum_{S_w(d)=0} z_d \cdot y_d$ for some 
$z_d \in \mathbb{C}$. Under the assumption that $w^{k}(z)=z$, one has 
$\sum_{S_w(d)=0} z_d \cdot y_d= z'=w^{k} (z')=
\sum_{S_w(d)=0} d^{k} \cdot z_d \cdot y_d$. This means that 
$z_d=d^{k} \cdot z_d$ for any $d$ with $S_w(d)=0$. Since 
$d$ is not a root of unity, $z_d$ is zero for any $d$. Therefore, 
we have $z'=0$ and $z=z'' \in V_w^c$, and the assertion is established. 
\par 
Finally, as assertions 
$(3) \Rightarrow (4)$ and $(2) \Rightarrow (1)$ are 
obvious, the lemma is proved. 
\hfill $\Box$ \par\medskip 
Now, for an orbit data $\tau$, let $P(\tau)$ be the set of periodic roots 
with period $\ell_{w_{\tau}}$, that is, 
\[
P(\tau):=\big\{ \alpha \in \Phi_{N} \, \big| \, 
w_{\tau}^{\ell_{w_{\tau}}}(\alpha) 
= \alpha \big\}. 
\]
Moreover, we define a finite subset of the root system by 
\[
\Gamma_1(\tau)  := \big\{ \alpha_{\ell}^c \, \big| \, \ell=1, \dots,n 
\big\} \subset \Phi_N, 
\]
where $\alpha_{\ell}^c$ is given by 
\[
\alpha_{\ell}^c  :=  q_n \circ \cdots \circ q_{\ell+1} 
(e_0 -e_{\bar{\sigma}(\ell,1)}^1-e_{\bar{\sigma}(\ell,2)}^1-e_{\bar{\sigma}(\ell,3)}^1). 
\]
\begin{lemma} \label{lem:vanish} 
Let $d$ be a root of $S_{\tau}(t)=0$ and 
$s=(s_{\ell})$ be the solution of (\ref{eqn:system}). 
Then for each $1 \le \ell \le n$, 
$\alpha_{\ell}^c$ belongs to $P(\tau)$ 
if and only if $s_\ell=0$. 
\end{lemma}
{\it Proof}. Assume that $\alpha_{\ell}^c \in P(\tau)$, 
which is equivalent to saying that $(\alpha_{\ell}^c,y_d)=0$ 
from Lemma \ref{lem:per}. 
By (\ref{eqn:wexp}), we have 
\[
\begin{array}{ll}
(\alpha_{\ell}^c,y_d) &
= (e_0-e_{\bar{\sigma}(\ell,1)}^1-e_{\bar{\sigma}(\ell,2)}^1-e_{\bar{\sigma}(\ell,3)}^1, 
q_{\ell+1} \circ \cdots \circ q_n (y_d)) \\[2mm] 
~&= \displaystyle
\Big( \sum_{k=1}^{\ell} s_k + d \sum_{k=\ell+1}^{n} s_k \Big) 
+ \sum_{i=1}^3 v_{\ell,i} \\[2mm] 
~ &= \displaystyle 
\Big( \sum_{k=1}^{\ell} s_k + d \sum_{k=\ell+1}^{n} s_k \Big) 
+ \Big( -\sum_{k=1}^{\ell-1} s_k +(d-2) s_{\ell} - 
d \sum_{k=\ell+1}^{n} s_k  \Big) \\[2mm] 
~&=\displaystyle (d-1) s_{\ell}. 
\end{array}
\]
Thus the equation $(\alpha_{\ell}^c,y_d)=0$ is equivalent to 
saying that $s_{\ell}=0$, since $d \neq 1$. 
\hfill $\Box$ \par\medskip 
Propositions \ref{pro:tenta}, \ref{pro:nonzero} and \ref{pro:det} mentioned 
below run parallel with Theorems \ref{thm:main1}--\ref{thm:main3} 
in terms of condition (\ref{eqn:TR}) (see Proposition \ref{pro:tenta}). 
Namely, Proposition \ref{pro:tenta} states that $\tau$ admits 
a tentative realization if and only if $\tau$ satisfies condition (\ref{eqn:TR}), 
Proposition \ref{pro:nonzero} shows that the sibling $\check{\tau}$ of any 
orbit data satisfies the condition, and finally 
Proposition \ref{pro:det} gives a sufficient condition for (\ref{eqn:TR}). 
\begin{proposition} \label{pro:tenta}
Let $\tau$ be an orbit data with $\lambda(\tau) > 1$ and 
$d$ be a root of $S_{\tau}(t)=0$. 
Then, $\tau$ satisfies the condition 
\begin{equation} \label{eqn:TR}
\Gamma_1(\tau) \cap P(\tau) = \emptyset, 
\end{equation} 
if and only if there is a tentative realization 
$\overline{f} \in \mathcal{Q}(C)^n$ of $\tau$ such that $\delta(\overline{f}) = d$. 
Moreover, the tentative realization $\overline{f} \in \mathcal{Q}(C)^n$ of 
$\tau$ with $\delta(\overline{f}) = d$ is uniquely determined. 
\end{proposition}
{\it Proof}. 
This proposition follows easily from Proposition \ref{pro:tenta0} and 
Lemma \ref{lem:vanish}. 
\hfill $\Box$ \par\medskip 
\begin{proposition} \label{pro:nonzero}
For any orbit data $\tau=(n,\sigma,\kappa)$ with 
$\lambda(\tau) > 1$, 
there is a data $\check{\tau}=(\check{n}, \check{\sigma}, 
\check{\kappa})$ with $\check{n} \le n$ 
such that $\check{\tau}$ satisfies condition (\ref{eqn:TR}) and 
$\lambda(\tau)=\lambda(\check{\tau})$. 
\end{proposition}
{\it Proof}. Let $d$ be a root of $S_{\tau}(t)=0$ and 
$(v,s) \in \mathbb{C} \times (\mathbb{C}^{3n} \setminus \{0\}) 
\times (\mathbb{C}^{n} \setminus \{0\})$ 
be the solution of 
(\ref{eqn:system}) for $\tau$ as in Corollary \ref{cor:sol}. 
If $s_{\ell} \neq 0$ for any $\ell=1,\dots,n$, then putting 
$\check{\tau}=\tau$ leads to the proposition. Otherwise, 
assume that $s_{\ell} = 0$ for some $\ell$. Then we put 
$\check{n}:=n-1$, and for any 
$\iota \in \mathcal{K}(\check{n}) \cong \{(i,j) 
\in \mathcal{K}(n) \, | \, i \neq \ell \}$, 
choose $\nu(\iota)$ so that 
$i_1=\cdots =i_{\nu(\iota)-1} = \ell$ but 
$i_{\nu(\iota)} \neq \ell$. 
A new orbit data $\check{\tau}=(\check{n}, \check{\sigma}, \check{\kappa})$ 
is defined by 
$\check{\sigma}(\iota):=\sigma^{\nu(\iota)}(\iota)$ and 
$\check{\kappa}(\iota):=\sum_{k=0}^{\nu(\iota)-1} \kappa(\sigma^k(\iota))$ 
for any $\iota \in \mathcal{K}(\check{n})$. 
Then, since $v_{\iota_1}= d^{\kappa(\iota)} 
\cdot v_{\iota} + (d-1) \cdot s_{i_1}$ 
and $s_{\ell}=0$, we have 
$v_{\iota_{\check{\sigma}}}= d^{\check{\kappa}(\iota)} 
\cdot v_{\iota} + (d-1) \cdot s_{i_{\check{\sigma}}}$ 
for any $\iota \in \mathcal{K}(\check{n})$, 
where $\check{\sigma}(\iota)=\iota_{\check{\sigma}}=(i_{\check{\sigma}},j_{\check{\sigma}})$. 
Moreover, as $v$ satisfies (\ref{eqn:BS}), the new vector 
$\check{v}=(v_{\iota})_{\iota \in \mathcal{K}(\check{n})}$ 
also satisfies (\ref{eqn:BS}) with $n=\check{n}$ and 
$\check{s}=(s_1,\dots,s_{\ell-1}, s_{\ell+1}, 
\dots, s_n) \neq 0$. 
Hence, $(d,\check{v},\check{s})$ is a solution of (\ref{eqn:system}) 
for $\check{\tau}$, which means that 
$S_{\tau}(t)=S_{\check{\tau}}(t)$ and thus 
$\lambda(\tau)=\lambda(\check{\tau})$. 
\par
Therefore, either $\check{s}_{\ell} \neq 0$ for any $\ell$, 
or we can repeat the above argument to eliminate 
$\check{s}_{\ell}=0$ from $\check{s}$. Since each step 
reduces $n$ by $1$, $\check{\tau}$ satisfies 
$\check{s}_{\ell} \neq 0$ 
for any $\ell$ after finitely many steps. 
\hfill $\Box$ \par\medskip 
\begin{proposition} \label{pro:det} 
For any orbit data $\tau$ satisfying conditions $(1)$ and $(2)$ in Theorem 
\ref{thm:main3}, there is an estimate $2^n-1 < \lambda(\tau) < 2^n$. 
Moreover, $\tau$ satisfies condition (\ref{eqn:TR}). 
\end{proposition}
The proof of this proposition is given in Section \ref{sec:proof}. 
\section{Realizability} \label{sec:real}
The aim of this section is to construct a realization of an orbit data 
$\tau$ and to establish our main theorems. 
In the previous section, we construct a tentative realization of 
$\overline{f} \in \mathcal{Q}(C)^n$ of $\tau$ 
under condition (\ref{eqn:TR}). However, $\overline{f}$ does not 
necessarily become a realization of $\tau$. 
We give two such examples 
after stating a preliminary lemma. To this end, 
for $\iota=(i,j), \iota'=(i',j') \in 
\mathcal{K}(n)$, let us define a root 
$\alpha_{\iota,\iota'}^k$ by 
\begin{equation*} \label{eqn:root2}
\alpha_{\iota,\iota'}^k := 
q_n \circ \cdots \circ q_{i'+1} ( e_{\bar{\sigma}(\iota)}^{k+1} - 
e_{\bar{\sigma}(\iota')}^1) \in \Phi_N. 
\end{equation*}
\begin{lemma} \label{lem:indA}
Assume that $\tau$ satisfies condition (\ref{eqn:TR}). Then, 
for a tentative realization $\overline{f}$ of $\tau$ mentioned 
in Proposition \ref{pro:tenta}, the following are equivalent. 
\begin{enumerate}
\item $p_{\iota}^{m} \approx p_{\iota'}^-$ for $m \ge 0$. 
\item $p_{\iota}^{m+\ell} 
\approx p_{\iota'}^{\ell}$ for some $\ell \ge 0$.  
\item $p_{\iota}^{m+\ell} 
\approx p_{\iota'}^{\ell}$ for any $\ell \ge 0$. 
\item $d^k \cdot v_{\iota} =v_{\iota'}$, where 
$m=\theta_{i,i'}(k) \ge 0$. 
\item $\alpha_{\iota,\iota'}^k \in P(\tau)$.  
\end{enumerate}
\end{lemma}
{\it Proof}. One can easily check that 
assertions (1)--(4) are equivalent (see also Remark \ref{rem:checkp}). 
Moreover, by virtue of (\ref{eqn:wexp}), 
the condition $\alpha_{\iota,\iota'}^k \in P(\tau)$, 
or $(\alpha_{\iota,\iota'}^k,y_{d})=0$ from Lemma \ref{lem:per}, 
is equivalent to saying that  
\[
0=(\alpha_{\iota,\iota'}^k,y_{d})
=(e_{\bar{\sigma}(\iota)}^{k+1}-e_{\bar{\sigma}(\iota')}^1,q_{i'+1} 
\circ \cdots \circ q_n (y_d))=d^{k} \cdot 
v_{\iota}-v_{\iota'},
\]
where the last equality follows from the fact that the coefficients 
of $e_{\bar{\sigma}(\iota)}^{k+1}$ and $e_{\bar{\sigma}(\iota')}^1$ in 
$q_{i'+1} \circ \cdots \circ q_n (y_d)$ are 
$d^{k} \cdot v_{\iota}$ and $v_{\iota'}$ 
respectively, since $\theta_{i,i'}(k) \ge 0$ 
(see also (\ref{eqn:wexp})). 
The lemma is established. 
\hfill $\Box$ \par\medskip 
\begin{remark} \label{rem:RIN}
If $\overline{f}$ is a realization of $\tau$ and 
$p_{\iota}^{m} \approx p_{\iota'}^-$ for a positive integer $m >0$, 
then it follows that $p_{\iota}^{m} > p_{\iota'}^-$ and 
$p_{\iota}^{m+\ell} > p_{\iota'}^{\ell}$ for any $\ell \ge 0$ with 
$m+\ell \le \mu(\iota)$ and $\ell \le \mu(\iota')$ 
(see Remark \ref{rem:propmap}). 
\end{remark}
\begin{example} \label{ex:ex1}
Consider the orbit data $\tau=(2,\sigma,\kappa)$ given by 
\begin{equation*} 
\left\{
\begin{array}{l}
\sigma : (1,1) \mapsto (1,2) \mapsto (2,2) \mapsto (2,1) \mapsto (1,1), \quad 
(1,3) \mapsto (2,3) \mapsto (1,3) \\[2mm]
\phantom{(}\kappa(1,1)=\kappa(2,2)=4, \quad \kappa(1,2)=\kappa(1,3)=0, \quad 
\kappa(2,1)=1, \quad \kappa(2,3)=3 \\
(\mu(1,1)=\mu(2,2)=7, \quad \mu(1,2)=\mu(1,3)=0, \quad 
\mu(2,1)= 0, \quad 
\mu(2,3)=4 ). 
\end{array}
\right. 
\end{equation*}
Then, equations (\ref{eqn:system}) for $\tau$ admit a solution 
$(d,v,s)$ with  
$d =\lambda(\tau) \approx 1.582$, 
which is the unique root of $t^6-t^4-2t^3-t^2+1=0$ in $|t|>1$, 
$v=(v_{1,1},v_{1,2},v_{1,3},v_{2,1},v_{2,2},v_{2,3}) \approx 
(1,7.269,8.048,0,1,1.779)$ and $s =(s_1,s_2) \approx (1.717,-10.765)$.  
Proposition \ref{pro:tenta0} assures that there is a tentative realization 
$\overline{f}=(f_1,f_2) \in \mathcal{Q}(C)^2$ of $\tau$, all the 
indeterminacy points of which are proper as 
$v_{k,i} \neq v_{k,j}$ for any $i \neq j$. 
However, $\overline{f}$ is not a realization of $\tau$. 
Indeed, assume the contrary that $\overline{f}$ is a realization of $\tau$. 
The fact that $v_{1,1}=v_{2,2}$ and Lemma \ref{lem:indA} yield 
$p_{1,1}^1 \approx p_{2,2}^-$, which means that 
$p_{1,1}^1 = f_2(p_{1,1}^-) > p_{2,2}^-$ and thus 
$p_{1,1}^{\ell} > p_{2,2}^{\ell-1}$ for each 
$\ell=1,\dots, \mu(1,1)=7$ 
(see Remark \ref{rem:RIN}). 
On the other hand, since 
$p_{1,1}^{\mu(1,1)} \approx p_{\sigma(1,1)}^{+}=p_{1,2}^{+}$
and $p_{1,2}^{+}$ is proper, one has 
$p_{1,1}^{\mu(1,1)} \neq p_{1,2}^{+}$ but  
$p_{2,2}^{\mu(1,1)-1}=p_{1,2}^{+}$. 
Hence, $\overline{f}$ is not a realization of $\tau$. 
This argument implies that there should not be a periodic root 
$\alpha_{\iota,\iota'}^k \in P(\tau)$ 
with $m=\theta_{i,i'}(k) > 0$ and 
$\mu(\iota) < m+\mu(\iota')$ (see Figure \ref{fig:real}). 
We remark that the solutions of (\ref{eqn:system}) for $\tau$ 
are the same as the ones for 
another orbit data $\check{\tau}=(2,\check{\sigma},\kappa)$ 
given by 
\[
\check{\sigma} : (1,\ell) \mapsto (2,\ell) \mapsto (1,\ell), \quad 
(\ell \in \{1,2,3\}), 
\]
and that $\overline{f}$ is a realization of $\check{\tau}$. In particular, 
one has $\lambda(\check{\tau})=\lambda(\tau) > 1$. 
\end{example}
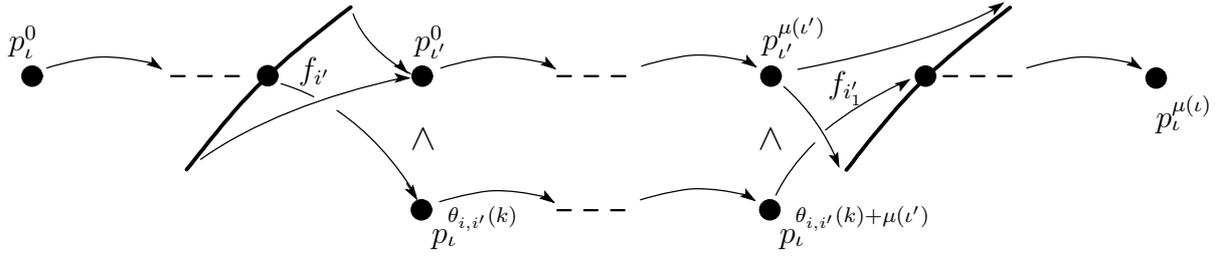
\begin{figure}[t]
\begin{center}
\unitlength 0.1in
\begin{picture}( 59.8500, 11.2300)(  0.6000,-14.4300)
\put(0.6000,-6.1000){\makebox(0,0)[lb]{$p_{\iota}^0$}}%
%
\special{pn 20}%
\special{sh 1.000}%
\special{ar 1400 698 46 46  0.0000000 6.2831853}%
%
\special{pn 20}%
\special{sh 1.000}%
\special{ar 2200 700 46 46  0.0000000 6.2831853}%
%
\special{pn 20}%
\special{sh 1.000}%
\special{ar 4006 700 46 46  0.0000000 6.2831853}%
%
\special{pn 8}%
\special{pa 2296 658}%
\special{pa 2326 648}%
\special{pa 2358 640}%
\special{pa 2390 632}%
\special{pa 2420 624}%
\special{pa 2452 618}%
\special{pa 2484 612}%
\special{pa 2514 606}%
\special{pa 2546 602}%
\special{pa 2578 600}%
\special{pa 2610 600}%
\special{pa 2642 602}%
\special{pa 2672 604}%
\special{pa 2704 608}%
\special{pa 2736 612}%
\special{pa 2768 618}%
\special{pa 2800 624}%
\special{pa 2832 630}%
\special{pa 2864 638}%
\special{pa 2868 638}%
\special{sp}%
%
\special{pn 8}%
\special{pa 2854 638}%
\special{pa 2886 660}%
\special{fp}%
\special{sh 1}%
\special{pa 2886 660}%
\special{pa 2842 606}%
\special{pa 2842 630}%
\special{pa 2820 640}%
\special{pa 2886 660}%
\special{fp}%
%
\special{pn 8}%
\special{pa 3336 648}%
\special{pa 3366 638}%
\special{pa 3398 630}%
\special{pa 3430 622}%
\special{pa 3460 614}%
\special{pa 3492 608}%
\special{pa 3524 602}%
\special{pa 3554 596}%
\special{pa 3586 592}%
\special{pa 3618 590}%
\special{pa 3650 590}%
\special{pa 3682 592}%
\special{pa 3712 594}%
\special{pa 3744 598}%
\special{pa 3776 602}%
\special{pa 3808 608}%
\special{pa 3840 614}%
\special{pa 3872 620}%
\special{pa 3904 628}%
\special{pa 3908 628}%
\special{sp}%
%
\special{pn 8}%
\special{pa 3894 628}%
\special{pa 3926 650}%
\special{fp}%
\special{sh 1}%
\special{pa 3926 650}%
\special{pa 3882 596}%
\special{pa 3882 620}%
\special{pa 3860 630}%
\special{pa 3926 650}%
\special{fp}%
%
\special{pn 8}%
\special{pa 256 658}%
\special{pa 286 648}%
\special{pa 318 640}%
\special{pa 350 632}%
\special{pa 380 624}%
\special{pa 412 618}%
\special{pa 444 612}%
\special{pa 474 606}%
\special{pa 506 602}%
\special{pa 538 600}%
\special{pa 570 600}%
\special{pa 602 602}%
\special{pa 632 604}%
\special{pa 664 608}%
\special{pa 696 612}%
\special{pa 728 618}%
\special{pa 760 624}%
\special{pa 792 630}%
\special{pa 824 638}%
\special{pa 828 638}%
\special{sp}%
%
\special{pn 8}%
\special{pa 814 638}%
\special{pa 846 660}%
\special{fp}%
\special{sh 1}%
\special{pa 846 660}%
\special{pa 802 606}%
\special{pa 802 630}%
\special{pa 780 640}%
\special{pa 846 660}%
\special{fp}%
%
\special{pn 20}%
\special{pa 980 1188}%
\special{pa 1000 1164}%
\special{pa 1020 1138}%
\special{pa 1040 1112}%
\special{pa 1060 1088}%
\special{pa 1080 1062}%
\special{pa 1100 1036}%
\special{pa 1120 1012}%
\special{pa 1140 988}%
\special{pa 1160 962}%
\special{pa 1180 938}%
\special{pa 1202 914}%
\special{pa 1222 888}%
\special{pa 1244 864}%
\special{pa 1264 840}%
\special{pa 1286 818}%
\special{pa 1308 794}%
\special{pa 1330 770}%
\special{pa 1352 748}%
\special{pa 1374 726}%
\special{pa 1396 702}%
\special{pa 1420 680}%
\special{pa 1444 658}%
\special{pa 1466 638}%
\special{pa 1490 616}%
\special{pa 1514 594}%
\special{pa 1538 574}%
\special{pa 1564 554}%
\special{pa 1588 532}%
\special{pa 1612 512}%
\special{pa 1638 492}%
\special{pa 1662 472}%
\special{pa 1688 452}%
\special{pa 1712 432}%
\special{pa 1738 412}%
\special{pa 1762 392}%
\special{pa 1788 372}%
\special{pa 1814 352}%
\special{pa 1830 340}%
\special{sp}%
%
\special{pn 20}%
\special{sh 1.000}%
\special{ar 180 700 46 46  0.0000000 6.2831853}%
%
\special{pn 8}%
\special{ar 3150 1800 2520 1200  3.7294249 4.2893551}%
%
\special{pn 8}%
\special{pa 2120 710}%
\special{pa 2130 710}%
\special{fp}%
\special{sh 1}%
\special{pa 2130 710}%
\special{pa 2064 690}%
\special{pa 2078 710}%
\special{pa 2064 730}%
\special{pa 2130 710}%
\special{fp}%
%
\special{pn 8}%
\special{ar 2830 250 1010 600  2.3241020 2.9313652}%
%
\special{pn 8}%
\special{pa 2100 660}%
\special{pa 2120 680}%
\special{fp}%
\special{sh 1}%
\special{pa 2120 680}%
\special{pa 2088 620}%
\special{pa 2082 642}%
\special{pa 2060 648}%
\special{pa 2120 680}%
\special{fp}%
%
\special{pn 8}%
\special{ar 1200 2400 1250 1700  4.9291370 5.0701908}%
%
\special{pn 8}%
\special{ar 1200 2400 1250 1700  5.1755551 5.5978748}%
%
\special{pn 8}%
\special{pa 2120 1260}%
\special{pa 2170 1340}%
\special{fp}%
\special{sh 1}%
\special{pa 2170 1340}%
\special{pa 2152 1274}%
\special{pa 2142 1296}%
\special{pa 2118 1294}%
\special{pa 2170 1340}%
\special{fp}%
%
\special{pn 20}%
\special{pa 4396 1190}%
\special{pa 4416 1166}%
\special{pa 4436 1140}%
\special{pa 4454 1114}%
\special{pa 4474 1090}%
\special{pa 4494 1064}%
\special{pa 4514 1038}%
\special{pa 4534 1014}%
\special{pa 4554 990}%
\special{pa 4576 964}%
\special{pa 4596 940}%
\special{pa 4616 916}%
\special{pa 4638 890}%
\special{pa 4658 866}%
\special{pa 4680 842}%
\special{pa 4700 820}%
\special{pa 4722 796}%
\special{pa 4744 772}%
\special{pa 4766 750}%
\special{pa 4788 728}%
\special{pa 4812 704}%
\special{pa 4834 682}%
\special{pa 4858 660}%
\special{pa 4882 640}%
\special{pa 4906 618}%
\special{pa 4930 596}%
\special{pa 4954 576}%
\special{pa 4978 556}%
\special{pa 5002 534}%
\special{pa 5028 514}%
\special{pa 5052 494}%
\special{pa 5078 474}%
\special{pa 5102 454}%
\special{pa 5128 434}%
\special{pa 5152 414}%
\special{pa 5178 394}%
\special{pa 5204 374}%
\special{pa 5228 354}%
\special{pa 5244 342}%
\special{sp}%
%
\special{pn 20}%
\special{sh 1.000}%
\special{ar 2196 1398 46 46  0.0000000 6.2831853}%
%
\special{pn 20}%
\special{sh 1.000}%
\special{ar 4000 1398 46 46  0.0000000 6.2831853}%
%
\special{pn 8}%
\special{pa 2290 1356}%
\special{pa 2322 1346}%
\special{pa 2354 1338}%
\special{pa 2384 1330}%
\special{pa 2416 1322}%
\special{pa 2448 1316}%
\special{pa 2478 1310}%
\special{pa 2510 1304}%
\special{pa 2542 1300}%
\special{pa 2572 1298}%
\special{pa 2604 1298}%
\special{pa 2636 1300}%
\special{pa 2668 1302}%
\special{pa 2700 1306}%
\special{pa 2732 1310}%
\special{pa 2764 1316}%
\special{pa 2794 1322}%
\special{pa 2826 1328}%
\special{pa 2858 1336}%
\special{pa 2862 1336}%
\special{sp}%
%
\special{pn 8}%
\special{pa 2848 1336}%
\special{pa 2880 1358}%
\special{fp}%
\special{sh 1}%
\special{pa 2880 1358}%
\special{pa 2836 1304}%
\special{pa 2836 1328}%
\special{pa 2814 1338}%
\special{pa 2880 1358}%
\special{fp}%
%
\special{pn 8}%
\special{pa 3330 1346}%
\special{pa 3362 1336}%
\special{pa 3394 1328}%
\special{pa 3424 1320}%
\special{pa 3456 1312}%
\special{pa 3488 1306}%
\special{pa 3518 1300}%
\special{pa 3550 1294}%
\special{pa 3582 1290}%
\special{pa 3612 1288}%
\special{pa 3644 1288}%
\special{pa 3676 1290}%
\special{pa 3708 1292}%
\special{pa 3740 1296}%
\special{pa 3772 1300}%
\special{pa 3804 1306}%
\special{pa 3834 1312}%
\special{pa 3866 1318}%
\special{pa 3898 1326}%
\special{pa 3902 1326}%
\special{sp}%
%
\special{pn 8}%
\special{pa 3888 1326}%
\special{pa 3920 1348}%
\special{fp}%
\special{sh 1}%
\special{pa 3920 1348}%
\special{pa 3876 1294}%
\special{pa 3876 1318}%
\special{pa 3854 1328}%
\special{pa 3920 1348}%
\special{fp}%
%
\special{pn 20}%
\special{sh 1.000}%
\special{ar 4806 700 46 46  0.0000000 6.2831853}%
%
\special{pn 8}%
\special{ar 3436 1580 1000 1070  5.4017947 5.8772234}%
%
\special{pn 8}%
\special{pa 4336 1120}%
\special{pa 4356 1170}%
\special{fp}%
\special{sh 1}%
\special{pa 4356 1170}%
\special{pa 4350 1102}%
\special{pa 4336 1120}%
\special{pa 4312 1116}%
\special{pa 4356 1170}%
\special{fp}%
%
\special{pn 8}%
\special{ar 1206 80 4190 820  0.3257927 0.7933462}%
%
\special{pn 8}%
\special{ar 6606 1600 2640 1260  3.3649130 3.5946037}%
%
\special{pn 8}%
\special{ar 6576 1600 2640 1240  3.6641690 3.9102326}%
%
\special{pn 20}%
\special{sh 1.000}%
\special{ar 6000 710 46 46  0.0000000 6.2831853}%
%
\special{pn 8}%
\special{pa 5330 658}%
\special{pa 5362 648}%
\special{pa 5394 640}%
\special{pa 5424 632}%
\special{pa 5456 624}%
\special{pa 5488 618}%
\special{pa 5518 612}%
\special{pa 5550 606}%
\special{pa 5582 602}%
\special{pa 5612 600}%
\special{pa 5644 600}%
\special{pa 5676 602}%
\special{pa 5708 604}%
\special{pa 5740 608}%
\special{pa 5772 612}%
\special{pa 5804 618}%
\special{pa 5834 624}%
\special{pa 5866 630}%
\special{pa 5898 638}%
\special{pa 5902 638}%
\special{sp}%
%
\special{pn 8}%
\special{pa 5888 638}%
\special{pa 5920 660}%
\special{fp}%
\special{sh 1}%
\special{pa 5920 660}%
\special{pa 5876 606}%
\special{pa 5876 630}%
\special{pa 5854 640}%
\special{pa 5920 660}%
\special{fp}%
%
\special{pn 8}%
\special{pa 5146 360}%
\special{pa 5206 320}%
\special{fp}%
\special{sh 1}%
\special{pa 5206 320}%
\special{pa 5138 340}%
\special{pa 5162 350}%
\special{pa 5162 374}%
\special{pa 5206 320}%
\special{fp}%
%
\special{pn 8}%
\special{pa 4690 740}%
\special{pa 4720 720}%
\special{fp}%
\special{sh 1}%
\special{pa 4720 720}%
\special{pa 4654 740}%
\special{pa 4676 750}%
\special{pa 4676 774}%
\special{pa 4720 720}%
\special{fp}%
%
\special{pn 13}%
\special{pa 900 700}%
\special{pa 1300 700}%
\special{da 0.070}%
%
\special{pn 13}%
\special{pa 2900 700}%
\special{pa 3300 700}%
\special{da 0.070}%
%
\special{pn 13}%
\special{pa 2900 1400}%
\special{pa 3300 1400}%
\special{da 0.070}%
%
\special{pn 13}%
\special{pa 4900 700}%
\special{pa 5300 700}%
\special{da 0.070}%
\put(15.6000,-7.5000){\makebox(0,0)[lb]{$f_{i'}$}}%
\put(22.5000,-16.0000){\makebox(0,0)[lb]{$p_{\iota}^{\theta_{i,i'}(k)}$}}%
\put(21.6000,-6.1000){\makebox(0,0)[lb]{$p_{\iota'}^0$}}%
\put(40.5000,-16.0000){\makebox(0,0)[lb]{$p_{\iota}^{\theta_{i,i'}(k)+\mu(\iota')}$}}%
\put(39.6000,-6.1000){\makebox(0,0)[lb]{$p_{\iota'}^{\mu(\iota')}$}}%
\put(43.0000,-8.5000){\makebox(0,0)[lb]{$f_{i_1'}$}}%
\put(60.0000,-10.0000){\makebox(0,0)[lb]{$p_{\iota}^{\mu(\iota)}$}}%
\put(21.5000,-11.0000){\makebox(0,0)[lb]{\rotatebox[origin=c]{90}{\large{$>$}}}}%
\put(39.5000,-11.0000){\makebox(0,0)[lb]{\rotatebox[origin=c]{90}{\large{$>$}}}}%
\end{picture}%
\end{center}
\caption{Orbit segments starting at two indeterminacy points} 
\label{fig:real}
\end{figure}
\begin{example} \label{ex:ex2}
Consider the orbit data $\tau=(1,\sigma,\kappa)$ given by 
\begin{equation*} 
\left\{
\begin{array}{l}
\sigma : (1,1) \mapsto (1,2) \mapsto (1,1), \quad 
(1,3) \mapsto (1,3) \\[2mm]
\kappa(1,1)=\kappa(1,2)=4, \quad \kappa(1,3)=3.
\end{array}
\right. 
\end{equation*}
Then,  equations (\ref{eqn:system}) for $\tau$ admit a solution 
$(d,v,s)$ with 
$d =\lambda(\tau) \approx 1.582$, 
which is the unique root of $t^6-t^4-2t^3-t^2+1=0$ in $|t|>1$, 
$v=(v_{1,1},v_{1,2},v_{1,3}) \approx 
(-0.190,-0.190,-0.338)$ and $s =(s_1) = (1)$.  
Proposition \ref{pro:tenta0} assures that there is a tentative realization 
$\overline{f}=(f_1) \in \mathcal{Q}(C)$ of $\tau$, whose 
indeterminacy points satisfy $p_{1,1}^{\pm} \approx p_{1,2}^{\pm}$. 
However, $\overline{f}$ is not a realization of $\tau$. 
Indeed, assume the contrary that $\overline{f}$ is a realization of $\tau$. 
If it is assumed that $p_{1,2}^{\pm} > p_{1,1}^{\pm}$ 
(see Remark \ref{rem:inf2}), then 
it follows that $p_{1,2}^{\ell} > p_{1,1}^{\ell}$ for any $\ell =0,\dots, 3$. 
On the other hand, since 
$p_{1,1}^{3} \approx p_{1,2}^{3} \approx p_{1,1}^{+} \approx p_{1,2}^{+}$
and $p_{1,2}^{+} > p_{1,1}^{+}$, one has 
$p_{1,1}^{3}=p_{1,1}^{+}$ and $p_{1,2}^{3}=p_{1,2}^{+}$. 
The case $p_{1,1}^{\pm} > p_{1,2}^{\pm}$ is the same. 
Hence, $\overline{f}$ is not a realization of $\tau$. 
This argument implies that there should not be a periodic root 
$\alpha_{(1,i),(1,j)}^0 \in P(\tau)$ with 
$\mu(1,i) = \mu(1,j)$ and $(1,j)=\sigma(1,i)$ for some $i \neq j$. 
We remark that the solutions of (\ref{eqn:system}) for $\tau$ 
are the same as the ones for 
another orbit data $\check{\tau}=(2,\check{\sigma},\kappa)$ 
given by 
\[
\check{\sigma} : (1,\ell) \mapsto (1,\ell), \quad 
(\ell \in \{1,2,3\}), 
\]
and that $\overline{f}$ is a realization of $\check{\tau}$. In particular, 
one has $\lambda(\check{\tau})=\lambda(\tau) > 1$. 
\end{example}
Based on the above two examples, we give the following definition. 
\begin{definition} \label{def:roots2}
We define subsets $\Gamma_{2}^{1}(\tau)$, $\Gamma_{2}^{2}(\tau)$ and 
$\Gamma_{2}(\tau)$ of the set 
\begin{equation*} \label{eqn:roots2}
\overline{\Gamma}_{2}(\tau) = 
\big\{ \alpha_{\iota,\iota'}^k  \, \big|
 \, \iota=(i,j),\iota'=(i',j') \in 
\mathcal{K}(n), 0 \le 
\theta_{i,i'}(k) \le \mu(\iota) \big\} \subset \Phi_N. 
\end{equation*}
\begin{enumerate}
\item Let $\Gamma_{2}^{1}(\tau)$ be the set of all roots 
$\alpha_{\iota,\iota'}^k \in \overline{\Gamma}_{2}(\tau)$ 
with $\theta_{i,i'}(k)>0$ such that 
$\mu(\sigma^{\ell}(\iota))=\mu(\sigma^{\ell}(\iota'))+ \delta_{\ell,0} 
\cdot \theta_{i,i'}(k)$ for $\ell=0,\dots, h-1$ and 
$\mu(\sigma^{h}(\iota))<\mu(\sigma^{h}(\iota'))+ \delta_{h,0} 
\cdot \theta_{i,i'}(k)$ with some $h \ge 0$, 
where $\delta_{i,j}$ stands for the Kronecker delta. 
\item Let $\Gamma_{2}^{2}(\tau)$ be the set of all roots 
$\alpha_{\iota,\iota'}^k \in \overline{\Gamma}_{2}(\tau)$ 
with $k=0$, $i=i'$ and $j \neq j'$ 
such that 
$\mu(\sigma^{\ell}(\iota))=\mu(\sigma^{\ell}(\iota'))$ 
for any $\ell \ge 0$ 
and $\iota' = \sigma^{h}(\iota)$ for some $h \ge 1$. 
\item The set $\Gamma_{2}(\tau)$ 
is defined by the union 
$\Gamma_{2}(\tau):= \Gamma_{2}^{1}(\tau) \cup \Gamma_{2}^{2}(\tau)$. 
\end{enumerate}
\end{definition}
\begin{example} \label{ex:gamma12}
Consider the orbit datum $\tau$ and $\check{\tau}$ given in Example \ref{ex:ex2}. 
Then we have
\[
\begin{array}{l}
\Gamma_2^1(\tau) = \{ e_i^k-e_j^1 \, | \, i,j \in \{1,2,3\}, k=2,3,4 \}, \quad 
\Gamma_2^2(\tau) = \{ e_1^1-e_2^1, e_2^1-e_1^1 \}, \\[2mm]
\Gamma_2^1(\check{\tau}) = \Gamma_2^1(\tau), \quad 
\Gamma_2^2(\check{\tau}) = \emptyset .
\end{array}
\]
Moreover, it can be checked that $w_{\tau}$ admits periodic roots $e_1^1-e_2^1$ and $e_2^1-e_1^1$ in $\Gamma_{2}(\tau)$, 
and $w_{\check{\tau}}$ admits no periodic roots in $\Gamma_{2}(\check{\tau})$. 
\end{example}
From the above preliminaries, we have the following three propositions, 
which also run parallel with Theorems \ref{thm:main1}--\ref{thm:main3} 
in a similar way to 
Propositions \ref{pro:tenta}, \ref{pro:nonzero} and \ref{pro:det}. 
\begin{proposition} \label{pro:ind} 
Let $\tau$ be an orbit data satisfying $\lambda(\tau) > 1$ and 
condition (\ref{eqn:TR}), and $\overline{f} \in \mathcal{Q}(C)^n$ be 
the tentative realization 
of $\tau$ such that $\delta(\overline{f}) =d$ is a root of 
$S_{\tau}(t)=0$ as is mentioned in Proposition \ref{pro:tenta}. 
Then, $\overline{f}$ is a realization of $\tau$ if and only if 
$\tau$ satisfies the condition 
\begin{equation} \label{eqn:R}
\Gamma_{2}(\tau) \cap P(\tau) = \emptyset. 
\end{equation}
\end{proposition}
{\it Proof}. 
First, assume that $\overline{f} \in  \mathcal{Q}(C)^n$ is a 
realization of $\tau$, and also assume the contrary that there is a root 
$\alpha_{\iota,\iota'}^k \in 
\Gamma_2(\tau) \cap P(\tau) \neq \emptyset$. 
If $\alpha_{\iota,\iota'}^k \in \Gamma_2^1(\tau)$, then 
Lemma \ref{lem:indA} shows that 
$p_{\iota}^m \approx p_{\iota'}^-$ with $m=\theta_{i,i'}(k)>0$ 
and thus $p_{\iota}^m > p_{\iota'}^-$ (see Remark \ref{rem:RIN}). 
Let $h \ge 0$ be the integer given in item (1) of Definition 
\ref{def:roots2}. 
When $h=0$, the relation 
$p_{\sigma(\iota)}^{+}=p_{\iota}^{\mu(\iota)} > 
p_{\iota'}^{\mu(\iota)-m}$ means that 
$p_{\iota'}^{\mu(\iota)-m}=p_{\iota''}^+$ for some $\iota''$ and thus 
$\mu(\iota')=\mu(\iota)-m$, 
since the indeterminacy set is a cluster. 
But this is impossible because $\mu(\iota)-  m < \mu(\iota')$. 
In a similar manner, when $h \ge 1$, it follows that 
$p_{\sigma^{\ell}(\iota)}^{\pm} > 
p_{\sigma^{\ell}(\iota')}^{\pm}$ for each $1 \le \ell \le h$. 
Moreover, the relation 
$p_{\sigma^{h+1}(\iota)}^{+}=p_{\sigma^{h}(\iota)}^{\mu(\sigma^{h}(\iota))} > 
p_{\sigma^{h}(\iota')}^{\mu(\sigma^{h}(\iota))}$ means that 
$\mu(\sigma^{h}(\iota)) = \mu(\sigma^{h}(\iota'))$, which is a contradiction. 
Hence we have $\Gamma_2^1(\tau) \cap P(\tau) = \emptyset$.
On the other hand, if 
$\alpha_{\iota,\iota'}^k \in \Gamma_2^2(\tau)$, 
then the relation $<$ can not be well defined. 
Indeed, assuming that $p_{\iota}^- < p_{\iota'}^-$, one has 
$p_{\sigma^{\ell}(\iota)}^{\pm}  < 
p_{\sigma^{\ell}(\iota')}^{\pm}$ for any $\ell \ge 0$. 
Since $\sigma^h(\iota)=\iota'$ and 
$\mu(\sigma^{\ell}(\iota))=\mu(\sigma^{\ell}(\iota'))$ for 
any $\ell \ge 0$, it is easily seen that $\iota'$ satisfies either 
$\sigma^h(\iota')=\iota$, 
or $\sigma^h(\iota')=\iota''$ and $\sigma^h(\iota'')=\iota$
with $\{\iota,\iota',\iota''\}=\{(i,1),(i,2),(i,3)\}$. 
The former case is impossible because it follows in this case 
that $p_{\iota'}^- < p_{\iota}^-$, which contradicts the assumption that 
$p_{\iota}^- < p_{\iota'}^-$. 
The latter case is also impossible, because it follows in this case that 
$p_{\iota'}^- < p_{\iota''}^-$ and $p_{\iota''}^- < p_{\iota}^-$, 
which is also a contradiction. 
Similarly, the assumption $p_{\iota}^- > p_{\iota'}^-$ yields 
a contradiction.
Therefore, we have $\Gamma_2^2(\tau) \cap P(\tau) = \emptyset$. 
These observations show that condition (\ref{eqn:R}) holds.
\par
Conversely, assume that condition (\ref{eqn:R}) holds. For two 
indeterminacy points $p_{k,i}^{\pm}$ and $p_{k,j}^{\pm}$ satisfying 
$p_{k,i}^{\pm} \approx p_{k,j}^{\pm}$, we fix 
$p_{k,i}^{\pm} < p_{k,j}^{\pm}$ if either 
$\mu(k,i) < \mu(k,j)$, or $\mu(k,i) = \mu(k,j)$ 
and $p_{\sigma(k,i)}^{\pm} < p_{\sigma(k,j)}^{\pm}$ 
(see Remark \ref{rem:inf2}). 
In other words, when there is $0 \le h \le \infty$ such that 
$\mu(\sigma^{\ell}(k,i)) = \mu(\sigma^{\ell}(k,j))$ for 
$0 \le \ell < h$ 
and $\mu(\sigma^{h}(k,i)) < \mu(\sigma^{h}(k,j))$, we fix 
$p_{\sigma^{\ell}(k,i)}^{\pm} < p_{\sigma^{\ell}(k,j)}^{\pm}$ 
for each $0 \le \ell < h$. This is well-defined because of the absence of 
roots in $\Gamma_2^2(\tau) \cap P(\tau)$. 
Now, in order to check condition (\ref{eqn:orbit}), 
for each $\iota \in \mathcal{K}(n)$, 
by putting  
\[
\begin{array}{l}
Q_{\iota}^+(m):= 
\left\{
\begin{array}{ll}
\{ p_{\iota'}^+ \, | \, p_{\iota'}^+ \approx p_{\iota}^m \} \quad 
& (0 \le m \le \mu(\iota)-1) \\[2mm]
\{ p_{\iota'}^+ \, | \, p_{\iota'}^+ 
\hspace{0.3em}\raisebox{0.4ex}{$<$}\hspace{-0.75em}\raisebox{-1.0ex}{$\neq$}\hspace{0.3em}  
p_{\sigma(\iota)}^+ \} \quad 
& (m = \mu(\iota)), 
\end{array}
\right. \\[6mm]
Q_{\iota}^-(m):=
\left\{
\begin{array}{ll}
\{ p_{\iota'}^- \, | \, p_{\iota'}^- 
\hspace{0.3em}\raisebox{0.4ex}{$<$}\hspace{-0.75em}\raisebox{-1.0ex}{$\neq$}\hspace{0.3em}
p_{\iota}^- \} \quad 
& ~~(m=0) \\[2mm]
\{ p_{\iota'}^- \, | \, p_{\iota'}^- \approx p_{\iota}^{m} \} \quad 
& ~~(1 \le m \le \mu(\iota)), 
\end{array}
\right. 
\end{array}
\]
we will show the condition 
\begin{equation} \label{eqn:card}
\sum_{k=0}^m (\# Q_{\iota}^-(k) - \# Q_{\iota}^+(k)) \ge 0 \quad
(0 \le m < \mu(\iota)), \qquad 
\sum_{k=0}^{\mu(\iota)} (\# Q_{\iota}^-(k) - \# Q_{\iota}^+(k))=0. 
\end{equation}
In our situation, it should be noted that 
\[
\# Q_{\iota}^{\pm}(m)=
\max \{ \ell \ge 0 \, | \, \text{there is an } (\ell-1) \text{-th point } 
p_{\iota'}^{\pm} \text{ with } p_{\iota}^m \approx p_{\iota'}^{\pm} \}. 
\]
Hence, for each $0 \le m \le \mu(\iota)-1$, 
if $p_{\iota}^{m}$ is an $\ell$-th point with 
$\ell \ge \# Q_{\iota}^+(m)$, which means that 
$p_{\iota}^{m} \neq p_{\iota'}^+$ for any $\iota' \in \mathcal{K}(n)$, 
then $p_{\iota}^{m+1}$ is a 
$(\# Q_{\iota}^{-}(m+1) - \#Q_{\iota}^+(m)+\ell)$-th point (see Remark \ref{rem:propmap}). 
As $p_{\iota}^0$ is a $\# Q_{\iota}^-(0)$-th point, 
$p_{\iota}^m$ is a 
$(\sum_{k=0}^m \# Q_{\iota}^-(k) -\sum_{k=0}^{m-1} \# Q_{\iota}^+(k))$-th point with 
$(\sum_{k=0}^m \# Q_{\iota}^-(k) -\sum_{k=0}^{m-1} \# Q_{\iota}^+(k)) 
\ge \# Q_{\iota}^+(m)$ under condition (\ref{eqn:card}), which shows that 
$p_{\iota}^{m} \neq p_{\iota'}^+$ for any $0 \le m \le \mu(\iota)-1$ 
and $\iota' \in \mathcal{K}(n)$. 
Finally, since $p_{\iota}^{\mu(\iota)} \approx p_{\sigma(\iota)}^+$ are 
$(\sum_{k=0}^{\mu(\iota)} \# Q_{\iota}^-(k) -
\sum_{k=0}^{\mu(\iota)-1} \# Q_{\iota}^+(k))=\# Q_{\iota}^+(\mu(\iota))$-th points, we have 
$p_{\iota}^{\mu(\iota)} = p_{\sigma(\iota)}^+$. 
\par
Condition (\ref{eqn:card}) is an immediate consequence of the following 
two assertions: 
\begin{enumerate}
\item For any $p_{\iota'}^- \in Q_{\iota}^-(m_1)$ , there is a unique 
$m_2$ with $m_1 \le m_2 \le \mu(\iota)$ and 
$p_{\sigma(\iota')}^+ \in Q_{\iota}^+(m_2)$. 
\item For any $p_{\iota'}^+ \in Q_{\iota}^+(m_2)$ , there is a unique 
$m_1$ with $0 \le m_1 \le m_2$ and 
$p_{\sigma^{-1}(\iota')}^- \in Q_{\iota}^-(m_1)$. 
\end{enumerate}
Indeed, these two assertions lead to the bijection 
\[
\mathcal{F}_{\iota} : 
\{(m_1,p_{\iota'}^-) \, | \, 0 \le m_1 \le \mu(\iota), 
p_{\iota'}^- \in Q_{\iota}^-(m_1)\} 
\to 
\{(m_2,p_{\iota'}^+) \, | \, 0 \le m_2 \le \mu(\iota), 
p_{\iota'}^+ \in Q_{\iota}^+(m_2)\} 
\]
defined by $\mathcal{F}_{\iota}(m_1,p_{\iota'}^-)=(m_2,p_{\sigma(\iota')}^+)$. 
Since $m_1 \le m_2$, the bijection $\mathcal{F}_{\iota}$
 yields condition (\ref{eqn:card}). 
\par
To finish the proof of the proposition, 
we only prove assertion (1) as assertion (2) can be 
treated in a similar manner. For uniqueness, 
assuming that there are $m_2 < m$ such that 
$p_{\sigma(\iota')}^+ \in Q_{\iota}^+(m_2) \cap Q_{\iota}^+(m)$, namely 
$p_{\sigma(\iota')}^+ \approx p_{\iota}^{m_2} \approx p_{\iota}^{m}$, 
one has $\alpha_{\iota,\iota}^k \in \Gamma_{2}^1(\tau) \cap P(\tau)$ with 
$m-m_2=n \cdot k > 0$, which is a contradiction. 
\par 
Next we show the existence of $m_2$. Assume that 
$p_{\iota'}^- \in Q_{\iota}^-(m_1)$ for some $0 \le m_1 \le \mu(\iota)$, 
which means that $\alpha_{\iota,\iota'}^{k_1} \in P(\tau)$ 
with $m_1=\theta_{i,i'}(k_1) \ge 0$. If 
$\mu(\iota) > \mu(\iota')+m_1$, then one has 
$p_{\sigma(\iota')}^+ \in Q_{\iota}^+(\mu(\iota')+m_1)$. 
On the other hand, if $\mu(\iota) < \mu(\iota')+m_1$, then 
the root $\alpha_{\iota,\iota'}^{k_1}$ belongs to 
$\Gamma_2^1(\tau) \cap P(\tau)$, which is a contradiction. 
Finally, suppose that $\mu(\iota) = \mu(\iota')+m_1$, or in other words 
$p_{\sigma(\iota)}^{\pm} \approx p_{\sigma(\iota')}^{\pm}$. 
If $m_1=0$, then the assumption $p_{\iota}^{-} > p_{\iota'}^{-}$ leads to 
$p_{\sigma(\iota)}^{\pm} > p_{\sigma(\iota')}^{\pm}$ and thus 
$p_{\sigma(\iota')}^{+} \in Q_{\iota}^+(\mu(\iota))$. 
On the other hand, if $m_1>0$, then we also have 
$p_{\sigma(\iota')}^{+} \in Q_{\iota}^+(\mu(\iota))$. 
Indeed, assume the contrary that 
$p_{\sigma(\iota')}^{+} \notin Q_{\iota}^+(\mu(\iota))$ 
or $p_{\sigma(\iota')}^{+} > p_{\sigma(\iota)}^{+}$. Fix $1 \le h \le \infty$ 
satisfying $\mu(\sigma^{\ell}(\iota))=\mu(\sigma^{\ell}(\iota'))$ for 
$1 \le \ell < h$ and $\mu(\sigma^{h}(\iota))<\mu(\sigma^{h}(\iota'))$. 
If $h=\infty$, then the existence of $L>0$ with 
$\sigma^{L}=\mathrm{id}$ shows that 
$\mu(\iota)=\mu(\sigma^{L}(\iota))=\mu(\sigma^{L}(\iota'))=\mu(\iota')$, 
which is a contradiction. If $h < \infty$, then the root 
$\alpha_{\iota,\iota'}^{k_1}$ belongs to 
$\Gamma_2^1(\tau) \cap P(\tau)$, which is also a contradiction. 
Assertion (1) is established by combining all these observations. 
Therefore the proof of the proposition is complete. 
\hfill $\Box$ \par\medskip 
\begin{proposition} \label{pro:vari}
Let $\tau$ be an orbit data satisfying $\lambda(\tau) > 1$ and 
condition (\ref{eqn:TR}), and $\overline{f}$ be a tentative realization 
mentioned in Proposition \ref{pro:tenta}. Then, there is an orbit 
data $\check{\tau}$ such that 
$\lambda(\tau)=\lambda(\check{\tau})$ and 
$\overline{f}$ is a realization of $\check{\tau}$. 
In particular, $\check{\tau}$ satisfies condition (\ref{eqn:R}). 
\end{proposition}
{\it Proof}. 
Let $d$ be a root of $S_{\tau}(t)=0$, 
$(v,s) \in (\mathbb{C}^{3n} \setminus \{0\}) 
\times (\mathbb{C}^n \setminus \{0\})$ 
be a solution of (\ref{eqn:system}) as in 
Corollary \ref{cor:sol}, and $u_{\iota}$ be given in (\ref{eqn:uv}). 
If $\Gamma_2(\tau) \cap P(\tau)=\emptyset$, then putting 
$\check{\tau}=\tau$ leads to the proposition. Otherwise, we make a decomposition 
\[
P_2(\tau):=\Gamma_2(\tau) \cap P(\tau)=
P_2^1(\tau) \cup P_2^2(\tau) \cup P_2^3(\tau), 
\]
and divide the proof into three steps, where $P_2^\ell(\tau)$ is given by 
\[
P_2^{\ell}(\tau) := \left\{
\begin{array}{ll}
\{ \alpha_{\iota',\iota''}^{k} \in P_2(\tau) \, | \, 
\iota'=\iota'' \} \quad & (\ell=1) \\[2mm] 
\{ \alpha_{\iota',\iota''}^{k} \in P_2(\tau) \, | \, 
\iota' \neq \iota'', \mu(\iota'') - \mu(\iota') +
\theta_{i',i''}(k) >0 \} \quad & (\ell=2) \\[2mm]
\{ \alpha_{\iota',\iota''}^{k} \in P_2(\tau) \, | \, 
\iota' \neq \iota'', \mu(\iota'') - \mu(\iota') +
\theta_{i',i''}(k) =0 \} \quad & (\ell=3). 
\end{array}
\right.
\]
It should be noted that 
$\mu(\iota'')-\mu(\iota')+\theta_{i',i''}(k)$ 
turns out to be 
nonnegative provided $\alpha_{\iota',\iota''}^{k} \in  \Gamma_2(\tau)$, 
and that $\mu(\iota'')-\mu(\iota')+\theta_{i',i''}(k)=0$ 
if and only if $p_{\sigma(\iota')}^{\pm} \approx p_{\sigma(\iota'')}^{\pm}$. 
Now we fix a root $\alpha_{\iota',\iota''}^{k} \in P_2(\tau)$, 
which means that $d^k \cdot v_{\iota'}=v_{\iota''}$ and 
$p_{\iota'}^m \approx p_{\iota''}^-$ 
with $0 \le m=\theta_{i',i''}(k) \le \mu(\iota')$ 
by Lemma \ref{lem:indA}, 
and put $m'=\mu(\iota'')-\mu(\iota')+m \ge 0$. 
\\ \underline{Step1} : 
First suppose that $\alpha_{\iota',\iota'}^{k} \in P_2^1(\tau)$, 
which belongs to $\Gamma_2^1(\tau)$ and thus satisfies 
$\theta_{i',i'}(k) > 0$ or $k \ge 1$. 
As $d^k \cdot v_{\iota'}=v_{\iota'}$ and $d$ is not a root of unity, we have 
$v_{\iota'}=0$ and thus $d^\ell \cdot v_{\iota'}=v_{\iota'}$ 
for any $\ell \ge 0$. Hence, it follows that 
$p_{\iota'}^0 \approx p_{\iota'}^n \approx \cdots 
\approx p_{\iota'}^{k'' \cdot n}$ 
and $p_{\sigma(\iota')}^+ \approx p_{\iota'}^{\mu(\iota')} \approx 
p_{\iota'}^{\mu(\iota')-k'' \cdot n}$, where 
$k'' \ge 1$ is chosen so that $0 \le \mu(\iota')-k'' \cdot n <n$. 
A new orbit data $\check{\tau}=(n,\sigma,\check{\kappa})$ is defined by 
$\check{\kappa}(\iota'):=\kappa(\iota')-k''$, and 
$\check{\kappa}(\iota):=\kappa(\iota)$ if $\iota \neq \iota'$. 
Then $\alpha_{\iota',\iota'}^{k}$ is not defined for $\check{\tau}$. 
The relation $d^{\check{\kappa}(\iota)-1} \cdot v_{\iota}
=u_{\sigma(\iota)}=u_{\iota_1}$ shows that 
$v_{\iota_1} = d^{\check{\kappa}(\iota)} 
\cdot v_{\iota} +(d-1) \cdot s_{i_1}$ 
for any $\iota \in \mathcal{K}(n)$, and thus 
$(d,v,s)$ satisfies (\ref{eqn:system}) for $\check{\tau}$. 
Therefore, $\overline{f}$ is also a tentative realization of 
$\check{\tau}$. 
Moreover, it follows from the relation $S_{\tau}(t)=S_{\check{\tau}}(t)$ 
that $\lambda(\tau)=\lambda(\check{\tau})$. 
\par
Since $\sum_{\iota \in \mathcal{K}(n)} \kappa(\iota)$ is finite, 
we can assume that $\# P_2^1(\tau)=0$ by repeating this argument. 
In particular, an integer $\overline{k} \ge 0$ with 
$\alpha_{\iota,\overline{\iota}}^{\overline{k}} \in P_2(\tau)$ 
for given $\iota,\overline{\iota} \in \mathcal{K}(n)$ is at most unique, since if 
$\alpha_{\iota,\overline{\iota}}^{k'}, 
\alpha_{\iota,\overline{\iota}}^{k''} \in P_2(\tau)$ 
for some $k'<k''$, then 
$\alpha_{\iota,\iota}^{\overline{k}} \in P_2^1(\tau)$ with 
$\overline{k}=k''-k'>0$. 
\\ \underline{Step2} : 
Next we assume that  $\alpha_{\iota',\iota''}^{k} \in P_2^2(\tau)$, 
which also belongs to $\Gamma_2^1(\tau)$ and thus satisfies 
$m=\theta_{i',i''}(k) > 0$. 
Then a new orbit data $\check{\tau}=(n,\check{\sigma},\check{\kappa})$
is defined by 
\begin{equation} \label{eqn:neworbit}
\check{\sigma}(\iota):= 
\left\{
\begin{array}{ll}
\sigma(\iota'') & (\iota=\iota') \\[2mm]
\sigma(\iota') & (\iota=\iota'') \\[2mm]
\sigma(\iota) & (\text{otherwise}), \\[2mm]
\end{array}
\right. \quad 
\check{\mu}(\iota):= 
\left\{
\begin{array}{ll}
\mu(\iota')+m' & (\iota=\iota') \\[2mm]
\mu(\iota'')-m' & (\iota=\iota'') \\[2mm]
\mu(\iota) & (\text{otherwise}). \\[2mm]
\end{array}
\right. 
\end{equation}
Since $p_{\iota'}^{m} \approx p_{\iota''}^-$, 
one has 
$p_{\iota'}^{\check{\mu}(\iota')} = 
p_{\iota'}^{\mu(\iota')+m'} = 
p_{\iota'}^{\mu(\iota'')+m} 
\approx p_{\iota''}^{\mu(\iota'')}
\approx p_{\sigma(\iota'')}^+
=p_{\check{\sigma}(\iota')}^+$ and 
$p_{\iota''}^{\check{\mu}(\iota'')} = 
p_{\iota''}^{\mu(\iota'')-m'}= 
p_{\iota''}^{\mu(\iota')-m} \approx 
p_{\iota'}^{\mu(\iota')} 
\approx p_{\sigma(\iota')}^+
=p_{\check{\sigma}(\iota'')}^+$, 
which yield 
$d^{\check{\kappa}(\iota')-1} \cdot v_{\iota'}=
u_{\check{\sigma}(\iota')}$ 
and 
$d^{\check{\kappa}(\iota'')-1} \cdot v_{\iota''}=
u_{\check{\sigma}(\iota'')}$. 
For $\iota \neq \iota',\iota''$, the equation 
$d^{\kappa(\iota)-1} \cdot v_{\iota} 
= u_{\sigma(\iota)}$ 
leads to 
$d^{\check{\kappa}(\iota)-1} \cdot v_{\iota}
=u_{\check{\sigma}(\iota)}$. 
This shows that 
$v_{\iota_{\check{\sigma}}} 
= d^{\check{\kappa}(\iota)} 
\cdot v_{\iota} +(d-1) 
\cdot s_{i_{\check{\sigma}}}$ 
for any $\iota \in \mathcal{K}(n)$, where 
$\check{\sigma}(\iota)=\iota_{\check{\sigma}}=(i_{\check{\sigma}},j_{\check{\sigma}})$, and 
$(d,v,s)$ satisfies (\ref{eqn:system}) for $\check{\tau}$, 
which means 
that $\overline{f}$ is also a tentative realization of $\check{\tau}$ and 
that $\lambda(\tau)=\lambda(\check{\tau})$. 
\par
Now, we will show that $\#P_2^2(\check{\tau}) < \#P_2^2(\tau)$. 
Indeed, for 
$\iota,\overline{\iota} \notin \{\iota',\iota''\}$, 
it follows that $\alpha_{\iota,\overline{\iota}}^{k'} \in P_2^2(\check{\tau})$ 
if and only if $\alpha_{\iota,\overline{\iota}}^{k'} \in P_2^2(\tau)$. 
Hence, if one assumes that 
$\alpha_{\iota,\overline{\iota}}^{k'} \in P_2^2(\check{\tau})$ 
and $\alpha_{\iota,\overline{\iota}}^{k'} \notin P_2^2(\tau)$, 
then $\iota, \overline{\iota}$ must satisfy 
$\{\iota,\overline{\iota} \} \cap \{\iota',\iota'' \} \neq \emptyset$. 
When $\overline{\iota}=\iota'$ and $\iota \neq \iota''$, the relation 
$\check{\mu}(\iota) < \check{\mu}(\iota')+\theta_{i,i'}(k')$ 
yields $\mu(\iota) < \mu(\iota'')+\theta_{i,i''}(k'')$, 
and the relation 
$\mu(\iota) > \mu(\iota')+\theta_{i,i'}(k')$ yields 
$\check{\mu}(\iota) > \check{\mu}(\iota')+\theta_{i,i''}(k'')$, 
where $k''=k+k'$. 
This means that $\alpha_{\iota,\iota''}^{k''} \notin P_2^2(\check{\tau})$ 
and $\alpha_{\iota,\iota''}^{k''} \in P_2^2(\tau)$. 
On the other hand, when $\overline{\iota}=\iota'$ and $\iota = \iota''$, 
the relation $\alpha_{\iota'',\iota'}^{k'} \in P_2^2(\check{\tau})$ 
gives $d^{k'} \cdot v_{\iota''}=v_{\iota'}$ 
with $0 < \theta_{i'',i'}(k') \le \check{\mu}(\iota'')$. 
Combining it with $d^{k} \cdot v_{\iota'}=v_{\iota''}$, one has 
the equation $d^{k''} \cdot v_{\iota'}=v_{\iota'}$ for $k''=k+k'$. 
Since $\mu(\iota')=\mu(\iota'')-m'+m = \check{\mu}(\iota'')+m \ge 
\theta_{i'',i'}(k')+m=\theta_{i',i'}(k'') > 0$, 
we have $\alpha_{\iota',\iota'}^{k''} \in P_2^1(\tau)$, 
which contradicts the assumption that $\# P_2^1(\tau)=0$. 
In a similar manner, if 
$\alpha_{\iota,\overline{\iota}}^{k'} \in P_2^2(\check{\tau})$ 
and $\alpha_{\iota,\overline{\iota}}^{k'} \notin P_2^2(\tau)$, then it can be seen that 
$\alpha_{\overline{\iota}',\overline{\iota}''}^{k''} \notin P_2^2(\check{\tau})$ 
and $\alpha_{\overline{\iota}',\overline{\iota}''}^{k''} \in P_2^2(\tau)$, 
where $(\overline{\iota}',\overline{\iota}'',k'')=(\iota',\overline{\iota},k+k')$ 
when $\iota=\iota''$, 
$(\overline{\iota}',\overline{\iota}'',k'')=(\iota,\iota',|k-k'|)$ 
when $\overline{\iota}=\iota''$, and 
$(\overline{\iota}',\overline{\iota}'',k'')=(\iota'',\overline{\iota},|k-k'|)$ 
when $\iota=\iota'$. 
Finally, it follows from $\check{\mu}(\iota')> 
\check{\mu}(\iota'')+ \theta_{i',i''}(k)$ 
that $\alpha_{\iota',\iota''}^{k} \notin P_2^2(\check{\tau})$. 
Since $\alpha_{\iota',\iota''}^{k} \in P_2^2(\tau)$, 
these observations show that $\#P_2^2(\check{\tau}) < \#P_2^2(\tau)$. 
\par
If $\alpha_{\iota,\iota}^{k'} \notin P_2^1(\check{\tau})$, 
which means in fact that $\iota=\iota'$, then 
we further take a new orbit data so that $\# P_2^1(\check{\tau})=0$. 
Hence, this step yields the orbit data $\check{\tau}$ satisfying either 
$\#P_2^2(\check{\tau}) < \#P_2^2(\tau)$ and 
$\sum_{\iota \in \mathcal{K}(n)} \check{\kappa}(\iota)=
\sum_{\iota \in \mathcal{K}(n)} \kappa(\iota)$, or 
$\sum_{\iota \in \mathcal{K}(n)} \check{\kappa}(\iota) < 
\sum_{\iota \in \mathcal{K}(n)} \kappa(\iota)$. 
So we can also assume that $\# P_2^1(\tau)=\# P_2^2(\tau)=0$ 
by repeating this argument. 
\\ \underline{Step3} : 
Finally, suppose that $\alpha_{\iota',\iota''}^{k} \in P_2^3(\tau)$, 
which can be chosen so that $h=h(\iota',\iota'') \ge 0$ 
is minimal among all elements in $P_2^3(\tau)$, 
where $1 \le h \le \infty$ is 
determined by the relations 
$\mu(\sigma^{\ell}(\iota')) = \mu(\sigma^{\ell}(\iota''))$ 
for $1 \le \ell \le h-1$ and 
$\mu(\sigma^{h}(\iota')) < \mu(\sigma^h(\iota''))$. 
Then we define a new orbit data $\check{\tau}=(n,\check{\sigma},\kappa)$ 
as in (\ref{eqn:neworbit}) with $m'=0$. 
It can be checked that 
$\overline{f}$ is also a tentative realization of $\check{\tau}$ and 
$\lambda(\tau)=\lambda(\check{\tau})$, and that 
$\# P_2^1(\check{\tau})=\# P_2^2(\check{\tau})=0$ as $\mu$ is invariant 
in this procedure. 
\par
Now we claim that 
$\mu(\check{\sigma}^{\ell}(\iota'))=\mu(\sigma^{\ell}(\iota''))$ and 
$\mu(\check{\sigma}^{\ell}(\iota''))=\mu(\sigma^{\ell}(\iota'))$ 
for any $1 \le \ell \le h$. 
In particular, $\alpha_{\iota',\iota''}^{k}$ does not belong to $P_2^3(\check{\tau})$ 
since $\mu(\check{\sigma}^{\ell}(\iota')) = \mu(\check{\sigma}^{\ell}(\iota''))$ 
for $1 \le \ell < h$ and 
$\mu(\check{\sigma}^h(\iota')) > \mu(\check{\sigma}^h(\iota''))$. 
Indeed, when $s' , s'' \ge 1$ are the minimal integers such that 
$\sigma^{s'}(\iota') \in \{\iota',\iota''\}$ and 
$\sigma^{s''}(\iota'') \in \{\iota',\iota''\}$, one has 
$\check{\sigma}^{\ell}(\iota'')=\sigma^{\ell}(\iota')$ or 
$\mu(\check{\sigma}^{\ell}(\iota''))=\mu(\sigma^{\ell}(\iota'))$ 
for $1 \le \ell \le s'$ and 
$\check{\sigma}^{\ell}(\iota')=\sigma^{\ell}(\iota'')$ or 
$\mu(\check{\sigma}^{\ell}(\iota'))=\mu(\sigma^{\ell}(\iota''))$ 
for $1 \le \ell \le s''$. 
Moreover, assuming that 
$s' < \ell' \le s''$ and that the claim is verified for $\ell \le \ell'-1$, 
we fix $k' \ge 1$ with $1 \le \ell'' := \ell' - k' \cdot s' \le s'$. 
If $\check{\sigma}^{s'}(\iota'')=\sigma^{s'}(\iota')=\iota'$, 
then it follows that $\sigma^{\ell'}(\iota')=\sigma^{\ell''}(\iota')$ and 
$\mu(\check{\sigma}^{\ell'}(\iota''))=\mu(\check{\sigma}^{\ell'-s'}(\iota')) 
= \mu(\check{\sigma}^{\ell'-s'}(\iota''))= \cdots =\mu(\check{\sigma}^{\ell''}(\iota'))$, 
which show that $\mu(\check{\sigma}^{\ell'}(\iota''))=\mu(\sigma^{\ell'}(\iota'))$. 
The case $\check{\sigma}^{s'}(\iota'')=\sigma^{s'}(\iota')=\iota''$ 
can be treated in the same manner. 
Furthermore, when $s'' < \ell' \le s'$ or $s', s'' < \ell'$, 
the claim also can be verified for $\ell=\ell'$ in a similar way under 
the assumption that the claim is already verified for $\ell \le \ell'-1$. 
\par
We also claim that $\# P_2^3(\check{\tau};\ell)=0$ for 
$1 \le \ell < h$ and $\# P_2^3(\check{\tau};h) < \# P_2^3(\tau;h)$, where 
$P_2^3(\tau;\ell) := \{ \alpha_{\iota,\overline{\iota}}^{k'} \in 
P_2^3(\tau) \, | \, h(\iota,\overline{\iota})=\ell\}$. 
Indeed, since $\mu(\sigma^\ell(\iota))=\mu(\check{\sigma}^\ell(\iota))$ 
for any $\iota \in \mathcal{K}(n)$ and $\ell < h$ by the above claim, 
it follows from the assumption $\# P_2^3(\tau;\ell)=0$ that 
$\# P_2^3(\check{\tau};\ell)=0$ for $1 \le \ell < h$. 
Moreover, if $\alpha_{\iota,\overline{\iota}}^{k'} \in P_2^3(\check{\tau};h)$ 
and $\alpha_{\iota,\overline{\iota}}^{k'} \notin P_2^3(\tau;h)$, 
then $\iota, \overline{\iota}$ should satisfy 
$\{\iota,\overline{\iota} \} \cap \{\iota',\iota'' \} \neq \emptyset$. 
A little calculation shows that they further satisfy 
either $\iota=\iota''$ or $\overline{\iota}=\iota'$, and that 
$\alpha_{\iota',\overline{\iota}}^{k''} \in P_2^3(\tau;h)$ does not belong to 
$P_2^3(\check{\tau};h)$ 
when $\iota=\iota''$, and 
$\alpha_{\iota,\iota''}^{k''} \in P_2^3(\tau;h)$ does not belong to 
$P_2^3(\check{\tau};h)$ when $\overline{\iota}=\iota'$, 
where $k''=k'+k$. Finally, since $\alpha_{\iota',\iota''}^{k} \in P_2^3(\tau;h)$ 
and $\alpha_{\iota',\iota''}^{k} \notin P_2^3(\check{\tau};h)$, 
we can show that $\# P_2^3(\check{\tau};h) < \# P_2^3(\tau;h)$. 
\par 
Therefore, by repeating this argument, we can conclude that 
$\# P_2^1(\check{\tau})=\# P_2^2(\check{\tau})=\# P_2^3(\check{\tau})=0$ 
and establish the proposition. 
\hfill $\Box$ \par\medskip 
\begin{proposition} \label{pro:cri}
Let $\tau$ be an orbit data satisfying conditions $(1)$ and $(2)$ in Theorem 
\ref{thm:main3}. Then we have 
$\overline{\Gamma}_2(\tau) \cap P(\tau) = \{ 
\alpha_{\iota,\iota'}^0 \, | \, 
i_m=i_m', \, 
\kappa(\sigma^m(\iota))=\kappa(\sigma^m(\iota')), \, 
m \ge 0 \}$. 
In addition, if $\tau$ satisfies condition $(3)$ in Theorem \ref{thm:main3}, 
then it also satisfies condition (\ref{eqn:R}). 
\end{proposition}
The proof of this proposition is given in Section \ref{sec:proof}. We are now 
in a position to establish the main theorems. To this end, we make the 
definition of the set $\Gamma(\tau)$. 
\begin{definition} \label{def:root}
The finite subset $\Gamma(\tau)$ of the root system $\Phi_{N}$ is 
defined by 
\[
\Gamma(\tau) :=\Gamma_1(\tau) \cup \Gamma_2(\tau) \subset \Phi_N. 
\]
\end{definition}
\begin{theorem} \label{thm:main}
Let $\tau$ be an orbit data with $\lambda(\tau) > 1$ 
and $d$ be a root of $S_{\tau}(t)=0$. Then, $\tau$ satisfies the condition
\begin{equation} \label{eqn:main}
\Gamma(\tau) \cap P(\tau) = \emptyset, 
\end{equation} 
if and only if there is a realization 
$\overline{f}=(f_1,\dots,f_n) \in \mathcal{Q}(C)^n$ 
of $\tau$ such that $\delta(\overline{f}) = d$. 
The realization $\overline{f} \in \mathcal{Q}(C)^n$ of $\tau$ with 
$\delta(\overline{f}) = d$ is uniquely determined. 
Moreover, the blowup $\pi_{\tau} : X_{\tau} \to \mathbb{P}^2$ of 
$N=\sum_{\iota \in \mathcal{K}(n)} \kappa(\iota)$ points 
$\{ p_{\iota}^m \, | \, \iota =(i,j)\in \mathcal{K}(n), \, 
m=\theta_{i,0}(k), 1 \le k \le \kappa(\iota) \}$ on $C^*$ 
lifts $f=f_n \circ \cdots \circ f_1$ to the automorphism 
$F_{\tau} : X_{\tau} \to X_{\tau}$. 
Finally, $(\pi_{\tau},F_{\tau})$ realizes $w_{\tau}$ and $F_{\tau}$ 
has positive entropy 
$h_{\mathrm{top}}(F_{\tau})= \log \lambda(\tau) > 0$. 
\end{theorem}
This theorem is an immediate consequence of Propositions 
\ref{pro:auto}, \ref{pro:tenta} and \ref{pro:ind}. 
\begin{remark} \label{rem:ARC}
Lemmas \ref{lem:vanish} and \ref{lem:indA} provide another realizability 
condition instead of (\ref{eqn:main}). 
Namely, for an orbit data $\tau$ with $\lambda(\tau) >1$, 
let $d$ be a root of $S_{\tau}(t)=0$ and 
$(v,s) \in (\mathbb{C}^{3n} \setminus \{0\}) 
\times (\mathbb{C}^n \setminus \{0\})$ 
be a solution of (\ref{eqn:system}) as in 
Corollary \ref{cor:sol}. Then $\tau$ satisfies condition (\ref{eqn:main}) 
if and only if 
\begin{enumerate}
\item $s_{\ell} \neq 0$ for any $1 \le \ell \le n$, and 
\item $d^{k} \cdot v_{\iota} \neq v_{\iota'}$ for any $(k,\iota,\iota')$ 
with $\alpha_{\iota,\iota'}^{k} \in \Gamma_2(\tau)$. 
\end{enumerate}
\end{remark}
{\it Proofs of Theorems \ref{thm:main1}--\ref{thm:main3}}. 
Theorem \ref{thm:main1} is an immediate consequence of Theorem \ref{thm:main}, 
since $\alpha$ is a periodic root of $w_{\tau}$ with period some $k \ge 1$ 
if and only if $\alpha$ is a periodic root of $w_{\tau}$ with period 
$\ell_{w_{\tau}}$ (see Lemma \ref{lem:per}). 
Moreover, Theorem \ref{thm:main2} follows from Propositions \ref{pro:nonzero} 
and \ref{pro:vari}, and Theorem \ref{thm:main3} follows from Propositions 
\ref{pro:det} and \ref{pro:cri}. 
\hfill $\Box$ \\[4mm] 
{\it Proof of Theorem \ref{thm:main0}}. 
For any value $\lambda \neq 1 \in \Lambda$, Theorem \ref{thm:main2} 
and Proposition \ref{pro:iden} show that there is an orbit data $\tau$ such 
that $\lambda=\lambda(\tau)$ and 
$\tau$ satisfies the realizability condition (\ref{eqn:condi}). In 
particular, the automorphism $F_{\tau}$ mentioned in Theorem \ref{thm:main1} 
has entropy $h_{\mathrm{top}}(F_{\tau})= \log \lambda >0$. Note that 
when $\lambda=1 \in \Lambda$, the automorphism 
$\mathrm{id}_{\mathbb{P}^2} : \mathbb{P}^2 \to \mathbb{P}^2$ satisfies 
$\lambda(\mathrm{id}_{\mathbb{P}^2}^*)=\lambda=1$ and 
$h_{\mathrm{top}}(\mathrm{id}_{\mathbb{P}^2})=0$. 
On the other hand, it follows from Proposition \ref{pro:expent} that 
the entropy of any automorphism 
$F: X \to X$ is given by $h_{\mathrm{top}}(F)= \log \lambda$ for some 
$\lambda \in \Lambda$. 
Therefore, Theorem \ref{thm:main0} is proved. 
\hfill $\Box$ \par\medskip 
\begin{example} \label{ex:auto}
We consider the orbit data $\tau=(n,\sigma,\kappa)$ given by $n=2$, 
$\sigma=\text{id}$, $\kappa(1,\ell)=3$ and $\kappa(2,\ell)=4$ 
for any $\ell=1,2,3$. 
Then $\tau$ satisfies the assumptions in Theorem \ref{thm:main3}, and thus 
$w_{\tau}$ is realized by a pair $(\pi_{\tau},F_{\tau})$, 
where $\pi_{\tau} : X_{\tau} \to \mathbb{P}^2$ is a blowup of $21$ points. 
We can check that equations (\ref{eqn:system}) admit a solution $(d,v,s)$
with $d=\lambda(\tau)\approx 3.87454251$, which is a root of 
$t^6-4t^5+t^4-2t^3+t^2-4t+1=0$, 
$v=(v_{1,1},v_{1,2},v_{1,3},v_{2,1},v_{2,2},v_{2,3}) \approx 
(1.749,1.749,1.749,0.233,0.233,0.233)$ and 
$s=(s_1,s_2) \approx (-100,-52.274)$. 
Therefore, the entropy of $F_{\tau}$ is given by 
$h_{\mathrm{top}}(F_{\tau})=\log \lambda(\tau) \approx 1.35442759$. 
Moreover, for any $\iota \neq \iota'$ with $i=i' \in \{ 1,2 \}$, 
the equality $v_{\iota}=v_{\iota'}$ shows that 
the element $w_{\tau}$ 
admits a periodic root $\alpha_{\iota,\iota'}^0$, 
which is not contained in 
$\Gamma(\tau)$. Therefore, the automorphism $F_{\tau}$ does not appear 
in the paper of McMullen \cite{M}. On the other hand, for any data 
$\hat{\tau}=(1,\hat{\sigma},\hat{\kappa})$, let 
$F_{\hat{\tau}} : X_{\hat{\tau}} \to X_{\hat{\tau}}$ be an automorphism 
that Diller in \cite{D} constructs from a single quadratic map preserving 
a cuspidal cubic. 
We claim that $F_{\tau}$ is not topologically 
conjugate to $F_{\hat{\tau}}^m$ for any $m \ge 1$. 
Indeed, assume the contrary that 
$F_{\tau}$ is topologically conjugate to $F_{\hat{\tau}}^m$ for some data 
$\hat{\tau}$ and $m \ge 1$. Then the topological conjugacy yields 
$\lambda(\tau)=\lambda(\hat{\tau})^m$. 
Moreover, since $X_{\tau}$ is obtained by blowing up 
$21$ points, so is $X_{\hat{\tau}}$, which means that 
$\sum_{\ell=1}^3 \hat{\kappa}(1,\ell)=21$ and thus 
there are $190$ possibilities for $\hat{\kappa}$. 
As $\hat{\sigma}$ has $6$ possibilities, $\hat{\tau}$ has $1140$ 
possibilities. 
However, with the help of a computer, it may be easily seen that there are no 
data $\hat{\tau}$ and $m \ge 1$ satisfying these conditions. 
Our claim is proved. 
\end{example}
\section{Proof of Realizability with Estimates} 
\label{sec:proof}
As is seen in Section \ref{sec:real}, Propositions \ref{pro:det} and 
\ref{pro:cri} prove Theorem \ref{thm:main3}, or the realizability of 
orbit data. In this section, we establish these propositions 
by applying some estimates mentioned below. 
Let $\overline{c}_{\iota,k}(d)$ and $c_{i,j}(d)$ be 
polynomials of $d$ defined by 
\begin{eqnarray}
v_{\iota}(d)  & = & 
\sum_{k=1}^n \overline{c}_{\iota,k}(d) \cdot s_k, 
\label{eqn:c1} \\
v_i(d) &:= &
v_{i,1}(d)+ v_{i,2}(d)+v_{i,3}(d) 
=- \sum_{j=1}^n c_{i,j}(d) \cdot s_j,  
\label{eqn:c2}
\end{eqnarray}
where $v_{\iota}(d)$ is given in (\ref{eqn:expre}), 
and let $\mathcal{A}_n(d,x)$ 
be an $n \times n$ matrix having the $(i,j)$-th entry: 
\[
\mathcal{A}_n(d,x)_{i,j} = 
\left\{ 
\begin{array}{ll}
d-2 + x_{i,i} & (i=j) \\[2mm]
-1 + x_{i,j} & (i>j) \\[2mm]
-d + x_{i,j} & (i<j) 
\end{array}
\right.
\]
with 
$x=(x_1, \dots, x_n)=(x_{ij}) \in M_n(\mathbb{R})$. 
Then equations (\ref{eqn:system}) yield 
\begin{equation} \label{eqn:mat}
\mathcal{A}_{\tau}(d)\, s=0, \qquad  \qquad 
s = \left(
\begin{array}{c}
s_1 \\[-1mm]
\vdots \\[-1mm]
s_n 
\end{array}
\right), 
\end{equation}
where $\mathcal{A}_{\tau}(d) := \mathcal{A}_n(d,c(d))$ with 
$c(d):=(c_{i,j}(d))$. 
Finally, let $\chi_{\tau}(d)$ be the determinant $|\mathcal{A}_{\tau}(d)|$ 
of the matrix $\mathcal{A}_{\tau}(d)$. 
\begin{lemma} \label{lem:sol}
Assume that $d$ is not a root of unity. Then, 
$d$ is a root of $\chi_{\tau}(t)=0$ if and only if 
$d$ is a root of $S_{\tau}(t)=0$. 
\end{lemma}
{\it Proof}. 
If $d$ is a root of $\chi_{\tau}(t)=0$, 
then there is a solution $s \neq 0$ of (\ref{eqn:mat}). Thus, 
$(d,v,s)$ satisfies (\ref{eqn:system}), where 
$v=(v_{\iota})$ is given in (\ref{eqn:expre}). This means that $d$ is a root of 
$S_{\tau}(t)=0$ (see Corollary \ref{cor:sol}). 
Conversely, if $d$ is a root of $S_{\tau}(t)=0$, then 
there is a unique solution 
$(v,s) \in (\mathbb{C}^{3n} \setminus \{0\}) \times 
(\mathbb{C}^{n} \setminus \{0\})$ 
of (\ref{eqn:system}). Moreover, $s$ is a solution of (\ref{eqn:mat}) 
and thus $d$ is a root of $\chi_{\tau}(t)=0$. 
\hfill $\Box$ \par\medskip 
Now we fix an orbit data $\tau$ satisfying 
conditions $(1)$ and $(2)$ in Theorem \ref{thm:main3}. 
\begin{lemma}
If $d > 1$, then for any $k \in \{1,\dots , n\}$ and 
$\iota \in \mathcal{K}(n)$, 
we have 
\[
-\frac{1}{d^2+d+1} \le \overline{c}_{\iota,k}(d) \le 0. 
\]
\end{lemma}
{\it Proof}. In view of equation (\ref{eqn:expre}), 
$\overline{c}_{\iota,k}(d)$ may be expressed as either 
$\overline{c}_{\iota,k}(d)=0$, or 
$\overline{c}_{\iota,k}(d)=- (d-1) \cdot 
d^{\eta_1}/(d^{\eta}-1)$ with $\eta_1+3 \le \eta$, or 
$\overline{c}_{\iota,k}(d)=-(d-1) \cdot 
(d^{\eta_1}+d^{\eta_2}) /(d^{\eta}-1)$ with 
$\eta_1+3 \le \eta_2$ and $\eta_2+3 \le \eta$, or 
$\overline{c}_{\iota,k}(d)=-(d-1) \cdot 
(d^{\eta_1}+d^{\eta_2}+d^{\eta_3})/(d^{\eta}-1)$ with 
$\eta_1+3 \le \eta_2$, $\eta_2+3 \le \eta_3$ and $\eta_3+3 \le \eta$, 
since $\# \{ m \, | \, 1 \le m \le |\iota|, \, 
i_m=k \} \le \#\{(k,1), (k,2), (k,3)\}=3$. 
We only consider the case 
$\overline{c}_{\iota,k}(d)=-(d-1) \cdot 
(d^{\eta_1}+d^{\eta_2}) /(d^{\eta}-1)$ 
as the remaining cases can be treated in the same manner. 
Since $d > 1$, the inequality $\overline{c}_{\iota,k}(d) < 0$ 
is trivial. Moreover, one has 
\[
\frac{\overline{c}_{\iota,k}(d)}{d-1}=- \frac{d^{\eta_1}+
d^{\eta_2}}{d^{\eta}-1}
\ge -\frac{d^{\eta_2-3}+d^{\eta_2}}{d^{\eta}-1}
\ge -\frac{d^{\eta-6}+d^{\eta-3}}{d^{\eta}-1}
= -(1+\frac{1}{d^{\eta}-1}) (d^{-6}+d^{-3}) 
\ge -\frac{1}{d^3-1}. 
\]
Thus the lemma is established. 
\hfill $\Box$ \par\medskip 
Since $c_{i,j}(d)= - \sum_{\ell=1}^3 \overline{c}_{(i,\ell),j}(d)$ from 
(\ref{eqn:c2}), 
the above lemma leads to the inequality 
\[
0 \le c_{i,j}(d) \le h(d), \qquad h(d):=\frac{3}{1+d+d^2}. 
\]
Note that for any $d \ge 2$ and any $0 \le x_{i,j} \le h(d)$, 
each diagonal entry $\mathcal{A}_n(d,x)_{i,i}$ of $\mathcal{A}_n(d,x)$ 
is positive and each non-diagonal entry 
$\mathcal{A}_n(d,x)_{i,j}$ with $i \neq j$ is negative. 
Let $\overline{\mathcal{A}}_n(d,x)_{i,j}$ 
be the $(i,j)$-cofactor of the matrix 
$\mathcal{A}_n(d,x)$. Then, the relation 
$|\mathcal{A}_n(d,x)|=\sum_{i=1}^n \overline{\mathcal{A}}_n(d,x)_{i,j} \cdot 
\mathcal{A}_n(d,x)_{i,j}$ 
holds for any $j=1,\dots,n$, where $|\mathcal{A}_n(d,x)|$ is the 
determinant of the matrix $\mathcal{A}_n(d,x)$. 
\begin{lemma} \label{lem:esti1}
For any $n \ge 2$, the following inequalities hold: 
\[
\left\{
\begin{array}{ll}
\overline{\mathcal{A}}_n(d,x)_{i,j} > 0 & 
(d > 2^n-1, \, 0 \le x_{i,j} \le h(d)) \\[2mm]
|\mathcal{A}_n(d,x)| > 0 & 
(d > 2^n, \, 0 \le x_{i,j} \le h(d)) \\[2mm]
|\mathcal{A}_n(2^n-1,x)| < 0 & 
(0 \le x_{i,j} \le h(d)).  
\end{array}
\right.
\]
\end{lemma}
{\it Proof}. We prove the inequalities by induction on $n$. 
For $n=2$, the first inequality holds since 
\begin{equation*} 
\overline{\mathcal{A}}_2(d,x)_{i,j}=
\left\{
\begin{array}{ll}
- \mathcal{A}_2(d,x)_{j,i} > 0 & (i \neq j) \\[2mm]
\mathcal{A}_2(d,x)_{i+1,i+1} > 0 & (i = j \in \mathbb{Z}/2\mathbb{Z} ) .
\end{array}
\right.
\end{equation*}
As $h(d) < \frac{3}{13}$ when $d > 3$, 
the remaining inequalities follow from the estimates 
\begin{equation*} 
\left\{
\begin{array}{l}
|\mathcal{A}_2(d,x)|=(d-2+x_{1,1}) (d-2+x_{2,2}) - 
(1-x_{2,1}) (d - x_{1,2}) > 2^2 - 1 \cdot 4 = 0 \\[2mm]
|\mathcal{A}_2(3,x)|=(1+x_{1,1}) (1+x_{2,2}) - 
(1-x_{2,1}) (3 - x_{1,2}) < (1+\frac{3}{13})^2 - 
(1- \frac{3}{13}) (3 - \frac{3}{13}) < 0. 
\end{array}
\right.
\end{equation*}
Therefore, the lemma is proved when $n=2$. Assume that the inequalities 
hold when $n=l-1$. A little calculation shows that 
$\overline{\mathcal{A}}_{i,j} := \overline{\mathcal{A}}_l(d,x)_{i,j}$ 
can be expressed as 
\[
\overline{\mathcal{A}}_{i,j}=
\left\{
\begin{array}{ll}
\displaystyle 
- \Bigl\{ \sum_{k=1}^{i-1} \overline{\mathcal{A}}_{l-1}(d,x^i)_{k,j-1} \cdot 
\mathcal{A}_{l}(d,x)_{k,i} 
+ \sum_{k=i+1}^{l} \overline{\mathcal{A}}_{l-1}(d,x^i)_{k-1,j-1} \cdot 
\mathcal{A}_{l}(d,x)_{k,i} \Bigr\} & (i < j) \\[2mm]
\displaystyle
- \Bigl\{ \sum_{k=1}^{i-1} \overline{\mathcal{A}}_{l-1}(d,x^i)_{k,j} \cdot 
\mathcal{A}_{l}(d,x)_{k,i} 
+ \sum_{k=i+1}^{l} \overline{\mathcal{A}}_{l-1}(d,x^i)_{k-1,j} \cdot 
\mathcal{A}_{l}(d,x)_{k,i} \Bigr\} & (i > j) \\[8mm]
~~~ |\mathcal{A}_{l-1}(d,x^i)| & (i=j), 
\end{array}
\right.
\]
where $x^i$ is the $(l-1,l-1)$-matrix obtained from $x$ by removing 
the $i$-th row and column vectors. 
Therefore, 
the first assertion follows from the induction hypothesis. Moreover, 
since 
$|\mathcal{A}_l(d,x)|=\sum_{i=1}^l \overline{\mathcal{A}}_l(d,x)_{i,j} 
\cdot \mathcal{A}_l(d,x)_{i,j}$, 
the bounds 
\[
|\mathcal{A}_l(d,(x_1,\dots,x_{j-1},0,x_{j+1},\dots,x_l)| 
\le |\mathcal{A}_l(d,x)| \le 
|\mathcal{A}_l(d,(x_1,\dots,x_{j-1},\overline{h}(d),x_{j+1},\dots,x_l)| 
\]
hold for any $j$, 
where $\overline{h}(d)$ is the column vector having each component 
equal to $h(d)$. 
Thus, we have 
\[
~~~~~~~~~~~~~~~~~~~~~~~~~~~
\begin{array}{l}
|\mathcal{A}_l(d,x)| \ge |\mathcal{A}_l(d,(0,\dots,0))| = 
(d-1)^{l-1} (d-2^l) > 0 
\qquad (d > 2^l), \\[2mm] 
\displaystyle |\mathcal{A}_l(2^l-1,x)| \le 
|\mathcal{A}_l(2^l-1,(\overline{h}(2^l-1),\dots,
\overline{h}(2^l-1)))| = - \frac{(2^l-2)^{l+1}}{2^{2l}-2^l+1}  < 0,  
\end{array}
~~~~~~~~~~~~~~~~~~~~~~~~~~~
\]
which show that the assertions are verified when $n=l$. Therefore, 
the induction is complete, and the lemma is established. 
\hfill $\Box$ \par\medskip 
Let $\delta > 1$ be the root of $S_{\tau}(t)=0$ in $|t| >1$ and $s \neq 0$ be 
the solution of equation (\ref{eqn:mat}) with $d=\delta$. 
Then $(\delta,v,s)$ satisfies equations (\ref{eqn:system}), 
where $v=(v_{\iota})$ is given in (\ref{eqn:expre}) with $d=\delta$. 
We notice that $\delta$, which is also a root of 
$\chi_{\tau}(d)=|\mathcal{A}_{\tau}(d)|=0$, 
satisfies $2^n-1 < \delta < 2^n$, 
since $\chi_{\tau}(2^n-1) < 0$ and $\chi_{\tau}(2^n) > 0$ 
from Lemma \ref{lem:esti1}. 
\begin{lemma} \label{lem:esti2}
For any $1 \le i \le n-1$, the ratio $s_{i+1}/s_i$ 
satisfies 
\[
z_1(n) < \frac{s_{i+1}}{s_{i}} < z_2(n), 
\]
where 
\[
~~~~~~~~~~~~~~~~~~~~~~~~~~~
z_1(n) := \frac{2^{n-1} (2^n+2)}{2^{2 n} + 2^{n+1}+6}, \quad 
z_2(n) := \frac{2^{2 n-1} + 2^{n}+3}{2^{2 n} + 2^{n+1}+3}. 
~~~~~~~~~~~~~~~~~~~~~~~~~~~
\]
\end{lemma}
{\it Proof}. 
For each $k_1, k_2 \ge 0$ with $k_1 + k_2 \le n-2$, let 
$\mathcal{A}^{k_1, k_2}_n(\delta)$ be the $n \times n$ matrix defined 
inductively as follows. First, put 
$\mathcal{A}_n^{0, 0}(\delta):=\mathcal{A}_{\tau}(\delta)$. Next, 
let $\mathcal{A}_n^{k_1, 0}(\delta)$ be the matrix obtained 
from $\mathcal{A}_n^{k_1-1, 0}(\delta)$ by replacing the $i$-th row 
of $\mathcal{A}_n^{k_1-1, 0}(\delta)$ with the sum 
of the $i$-th row and the $k_1$-th row multiplied by 
$-\mathcal{A}_n^{k_1-1, 0}(\delta)_{i,k_1} / 
\mathcal{A}_n^{k_1-1, 0}(\delta)_{k_1,k_1}$, 
where $i$ runs from $k_1+1$ to $n$. Finally, let 
$\mathcal{A}_{n}^{k_1, k_2}(\delta)$ be the 
matrix obtained from $\mathcal{A}_n^{k_1, k_2-1}(\delta)$ by replacing the 
$i$-th row of $\mathcal{A}_n^{k_1, k_2-1}(\delta)$ with the sum 
of the $i$-th row and the $(n-k_2+1)$-th row multiplied by 
$-\mathcal{A}_n^{k_1,k_2-1}(\delta)_{i,n-k_2+1} / 
\mathcal{A}_n^{k_1,k_2-1}(\delta)_{n-k_2+1,n-k_2+1}$, where 
$i$ runs from $k_1+1$ to $n-k_2$. 
Therefore, each entry of $\mathcal{A}^{k_1, k_2}_n(\delta)$ 
may be expressed as 
\[
\mathcal{A}^{k_1, k_2}_n(\delta)_{i,j}= 
\left\{
\begin{array}{ll}
\delta_{i,j} + \xi_{i,j}^{i-1} & (i \le k_1 \text{~and~} i \le j) \\[2mm]
\delta_{i,j} + \xi_{i,j}^{k_1+k_2} & (k_1+1 \le i, j \le n-k_2) \\[2mm]
\delta_{i,j} + \xi_{i,j}^{k_1+n-i} & (n-k_2+1 \le i \text{~and~} 
k_1+1 \le j \le i) \\[2mm]
0 & (\text{otherwise}), 
\end{array}
\right.
\]
where 
\[
~~~~~~~~~~~~~~~~~~~~~~~~~~~
\delta_{i,j} = 
\left\{
\begin{array}{ll}
\delta-2 & (i=j) \\[2mm]
-1 & (i>j) \\[2mm]
-\delta & (i<j), 
\end{array}
\right.
~~~~~~~~~~~~~~~~~~~~~~~~~~~
\]
and $\xi_{i,j}^k$ is given inductively by 
\[
\xi_{i,j}^{0}=c_{i,j}(\delta), \quad 
\xi_{i,j}^{k+1}= 
\left\{
\begin{array}{ll}
\displaystyle \xi_{i,j}^{k}-\frac{(1-\xi_{i,k}^k)
(\delta - \xi_{k,j}^k)}{\delta-2 + \xi_{k,k}^{k}}  & (k < k_1) \\[2mm]
\displaystyle \xi_{i,j}^{k}-\frac{(\delta - \xi_{i,n-k+k_1}^k)
(1 - \xi_{n-k+k_1,j}^k)}
{\delta-2 + \xi_{n-k+k_1,n-k+k_1}^{k}} & (k \ge k_1). 
\end{array}
\right.
\]
Moreover, it is seen that $\xi_{i,j}^{k}$ satisfies the estimates
\[
-\frac{(2^k-1) \delta}{\delta -2^k} \le \xi_{i,j}^{k} \le 
-\overline{\xi}_k, \qquad 
\overline{\xi}_k:=\frac{(2^k-1) \delta - (2^k \delta -1) h(\delta)}
{(\delta-2^k)+(2^k-1) h(\delta)}. 
\]
Note that $s$ satisfies 
$\mathcal{A}^{k_1, k_2}_n(\delta) \, s=0$ for any $k_1, k_2 \ge 0$. 
In particular, one has 
$\mathcal{A}^{i-1, n-i-1}_n(\delta) \, s=0$, the 
$i$-th and $(i+1)$-th components of which are given by 
\[
\left\{
\begin{array}{l}
(\delta -2 + \xi_{i,i}^{n-2}) \, s_i + 
(-\delta + \xi_{i,i+1}^{n-2}) \, s_{i+1}=0 \\[2mm] 
(-1 + \xi_{i+1,i}^{n-2}) \, s_i + 
(\delta -2 + \xi_{i+1,i+1}^{n-2}) \, s_{i+1}=0.
\end{array}
\right.
\]
Therefore, we have 
\[
\frac{s_{i+1}}{s_i} = 
\frac{\delta -2 + \xi_{i,i}^{n-2}}{\delta - \xi_{i,i+1}^{n-2}} < 
\frac{\delta-2-\overline{\xi}_{n-2}}{\delta+\overline{\xi}_{n-2}} = 
\frac{2 \delta^2 - (2^n-4) \delta - (2^{n+1}-6)}{2 (\delta^2 +2 \delta +3)}, 
\]
the righthand side of which is monotone increasing with respect to $\delta$, 
and thus is less than $z_2(n)$ since $\delta < 2^n$. In a similar manner, 
we have 
\[
\frac{s_{i+1}}{s_i} = 
\frac{1 - \xi_{i+1,i}^{n-2}}{\delta -2 + \xi_{i+1,i+1}^{n-2}} > 
\frac{1+\overline{\xi}_{n-2}}{\delta-2-\overline{\xi}_{n-2}} = 
\frac{2^{n-2} (1-h(\delta))}{\delta-2^{n-1} + (2^{n-1}-1) 
h(\delta)} > z_1(n). 
\]
Thus, the lemma is established. 
\hfill $\Box$ \par\medskip 
We remark that the functions $z_1(n)$ and $z_2(n)$ satisfy
\[
0 <  z_1(n) < \frac{1}{2} < z_2(n) < 1.  
\]
\\
{\it Proof of Proposition \ref{pro:det}}. 
Recall that $\delta$ satisfies $2^n-1 < \delta < 2^n$ from Lemma \ref{lem:esti1}. 
Moreover, it follows from Lemma \ref{lem:esti2} that $s_\ell \neq 0$ 
for any $\ell$. Thus Lemma \ref{lem:vanish} yields 
$\Gamma_1(\tau) \cap P(\tau) = \emptyset$. 
\hfill $\Box$ \par\medskip 
Next we prove Proposition \ref{pro:cri}. 
\begin{lemma} \label{lem:g}
For any $n \ge 2$, we have the following two inequalities: 
\begin{enumerate}
\item $g_1(n) < 0$, where 
$\displaystyle g_1(n):= \frac{1}{\delta^3-1}+1-\delta \cdot z_1(n)^{n-1}$,  
\item $g_2(n) > 0$, where 
$\displaystyle g_2(n):= z_1(n)^{n-2}-z_2(n)^{n-1}-\frac{1}{\delta^3-1}$.
\end{enumerate}
\end{lemma}
{\it Proof}. 
First, we claim that the following inequality holds: 
\begin{equation} \label{eqn:ineqxi}
z_1(n)^{n-1} > \frac{1}{2^{n-1}} -(n-1) \Bigl(\frac{1}{2^{3 n -4}} + 
\frac{1}{2^{4 n -3}} \Bigr). 
\end{equation}
Indeed, since 
\[
\Bigl(1- \frac{1}{2^{n-1}} - \frac{1}{2^{2n-1}} \Bigr)
\Bigl(1+ \frac{1}{2^{n-1}} + \frac{6}{2^{2n}} \Bigr)
= 1 - \frac{1}{2^{4 n-2}} (2^{n+2}+3) \le 1, 
\]
one has 
\[
z_1(n) \ge \frac{1}{2} \Bigl(1 + \frac{1}{2^{n-1}} \Bigr)
\Bigl(1 - \frac{1}{2^{n-1}} - \frac{1}{2^{2n-1}} \Bigr)
= \frac{1}{2} \Bigl\{ 1-\bigl( \frac{3}{2^{2n-1}} + \frac{1}{2^{3n-2}} \bigr) 
\Bigr\} > \frac{1}{2}
\Bigl\{ 1-\bigl( \frac{1}{2^{2n-3}} + \frac{1}{2^{3n-2}} \bigr) \Bigr\}. 
\]
Therefore, the claim holds from the Bernoulli inequality, namely, 
$(1+x)^n \ge 1+nx$ for any $x \ge -1$. 
By using inequality (\ref{eqn:ineqxi}), we prove the two inequalities 
in the lemma. 
\par 
In order to prove assertion (1), we consider the function of $n$: 
\[
\check{g}_1(n) := \frac{1}{(2^n-1)^3-1}+1- (2^n-1) \cdot z_1(n)^{n-1}. 
\]
Then the inequality $g_1(n) < \check{g}_1(n)$ holds since $\delta > 2^n-1$. 
Moreover, as $\check{g}_1(2) < 0$, one has $g_1(2) <0$. On the other hand, 
when $n \ge 3$, inequality (\ref{eqn:ineqxi}) yields 
\[
\begin{array}{rl}
\check{g}_1(n) < & \displaystyle 
\frac{(2^n-1)^3}{(2^n-1)^3-1} - \frac{2^n-1}{2^{n-1}} \Bigl\{ 1- (n-1) 
\bigl( \frac{1}{2^{2n-3}} + \frac{1}{2^{3n-2}} \bigr) \Bigr\} \\[2mm]
~ < & \displaystyle \frac{2^n-1}{2^{n-1}} \Bigl( -1+ 
\frac{2^{n-1} (2^n-1)^2}{(2^n-1)^3-1} + \frac{n-1}{2^{2n-3}} 
+ \frac{n-1}{2^{3n-2}} \Bigr). 
\end{array}
\]
Since the terms $\frac{2^{n-1} (2^n-1)^2}{(2^n-1)^3-1}$, 
$\frac{n-1}{2^{2n-3}}$ and $\frac{n-1}{2^{3n-2}}$ are monotone decreasing 
with respect to $n$, the function 
$-1 + \frac{2^{n-1} (2^n-1)^2}{(2^n-1)^3-1} + 
\frac{n-1}{2^{2n-3}} + \frac{n-1}{2^{3n-2}}$ is maximized when $n=3$, which is 
negative. Therefore, we have $\check{g}_1(n) < 0$, and thus $g_1(n) < 0$. 
\par 
Finally, in order to prove assertion (2), we consider the function of $n$: 
\[
\check{g}_2(n) := z_1(n)^{n-2} - z_2(n)^{n-1} -\frac{1}{(2^n-1)^3-1}. 
\]
Then the inequality $g_2(n) > \check{g}_2(n)$ holds since $\delta > 2^n-1$. 
Moreover, as $\check{g}_2(2),\check{g}_2(3)> 0$, one has $g_2(2), g_2(3) >0$. 
On the other hand, when $n \ge 4$, $\check{g}_2(n)$ can be estimated as 
\[
\begin{array}{ll}
\check{g}_2(n) & \displaystyle 
= z_1(n)^{n-2} \bigl( 1-z_1(n) \bigr) - 
\bigl( z_2(n)^{n-1}-z_1(n)^{n-1} \bigr) 
 -\frac{1}{(2^n-1)^3-1} \\[2mm]
& \displaystyle
\ge z_1(n)^{n-2} \bigl( 1-z_1(n) \bigr) - (n-1) 
\bigl( z_2(n)-z_1(n) \bigr) 
z_2(n)^{n-2} -\frac{1}{(2^n-1)^3-1}, 
\end{array}
\]
where the last inequality follows from the general inequality 
$x^n-y^n \le n(x-y) x^{n-1}$ for any $x \ge y \ge 0$. 
Since  
$z_2(n)-z_1(n)=\frac{9}{2} \frac{2^{2n}+2^{n+1}+4}
{(2^{2n}+2^{n+1}+3)(2^{2n}+2^{n+1}+6)}<
\frac{9}{2}\frac{1}{2^{2n}+2^{n+1}+3} < \frac{9}{8} \frac{1}{2^{2 (n-1)}}$, 
and $z_2(n)=\frac{1}{2}+\frac{3}{2}\frac{1}{2^{2n}+2^{n+1}+3}$ is 
monotone decreasing with respect to $n$, and thus is less than 
$\frac{13}{24}$, we have 
\[
~~~~~~~~~~~~~~~~~~~~~~~~~~~
(n-1) \bigl( z_2(n)-z_1(n) \bigr) z_2(n)^{n-2}  
< (n-1) \frac{9}{8} \Bigl( \frac{13}{24} \Bigr)^{n-2} 
\frac{1}{2^{2 (n-1)}} < \frac{1}{2^{2 (n-1)}}, 
~~~~~~~~~~~~~~~~~~~~~~~~~~~
\]
where we use the fact that the function $(n-1) \frac{9}{8} 
\Bigl( \frac{13}{24} \Bigr)^{n-2}$ is 
monotone decreasing and is less than $1$. 
Moreover, as $1-z_1(n) > z_1(n)$, one has 
\[
\begin{array}{rl}
\check{g}_2(n) > & \displaystyle 
z_1(n)^{n-1} - \frac{1}{2^{2 (n-1)}} -\frac{1}{(2^n-1)^3-1} 
\\[2mm]
> & \displaystyle 
\frac{1}{2^{n-1}} \Bigl\{ 1 -(n-1) \bigl(\frac{1}{2^{2 n -3}} + \frac{1}
{2^{3 n -2}} \bigr) \Bigr\} - \frac{1}{2^{2 (n-1)}} -\frac{1}{(2^n-1)^3-1} 
\\[2mm]
= & \displaystyle 
\frac{1}{2^{n-1}} \Bigl( 1 - \frac{n-1}{2^{2 n -3}} - \frac{n-1}
{2^{3 n -2}} - \frac{1}{2^{n-1}} - \frac{2^{n-1}}{(2^n-1)^3-1} \Bigr). 
\end{array}
\]
Since the terms $\frac{n-1}{2^{2n-3}}$, $\frac{n-1}{2^{3n-2}}$, 
$\frac{1}{2^{n-1}}$ and $\frac{2^{n-1}}{(2^n-1)^3-1}$ are monotone decreasing 
with respect to $n$, the function 
$1 - \frac{n-1}{2^{2 n -3}} - \frac{n-1}
{2^{3 n -2}} - \frac{1}{2^{n-1}} - \frac{2^{n-1}}{(2^n-1)^3-1}$ 
is minimized when $n=4$, which is positive. 
Therefore, we have $\check{g}_2(n) > 0$ and thus $g_2(n) > 0$, 
and so the proof is complete. 
\hfill $\Box$ \par\medskip 
\begin{proposition} \label{pro:coeffi}
Assume that $v_{\iota'}(\delta)=\delta^k \cdot 
v_{\iota}(\delta)$. 
Then we have $k=0$, $i_m'=i_m$ and 
$\kappa(\sigma^m(\iota'))=\kappa(\sigma^m(\iota))$ for any $m \ge 0$. 
\end{proposition} 
{\it Proof}. 
Viewing $v_{\iota'}(\delta)/(\delta-1)$ and 
$\delta^k \cdot v_{\iota}(\delta)/(\delta-1)$ 
as functions of $\delta$ (see (\ref{eqn:expre})), we expand them 
into Taylor series around infinity: 
\begin{eqnarray*}
\frac{v_{\iota'}(\delta)}{\delta-1}  & = & 
- s_{i_1'} \cdot \delta^{-\varepsilon_1(\iota')} - 
\cdots - s_{i_{|\iota'|}'} \cdot 
\delta^{-\varepsilon_{|\iota'|}(\iota')} - 
s_{i_{|\iota'|+1}'} \cdot 
\delta^{-\varepsilon_{|\iota'|+1}(\iota')} - 
\cdots, 
\\[2mm] \hspace{-7mm}
\frac{\delta^k \cdot v_{\iota}(\delta)}{\delta-1}  & = & 
- s_{i_1} \cdot \delta^{-\varepsilon_1(\iota)+k} - \cdots - 
s_{i_{|\iota|}} \cdot 
\delta^{-\varepsilon_{|\iota|}(\iota)+k} - 
s_{i_{|\iota|+1}} \cdot 
\delta^{-\varepsilon_{|\iota|+1}(\iota)+k} 
- \cdots. 
\end{eqnarray*}
In view of these expressions, the coefficient of $\delta^{-l}$ is either 
$-s_{\bullet}$ or $0$. Now assume the contrary that 
$v_{\iota'}(\delta)/(\delta-1)$ and 
$\delta^k\cdot v_{\iota}(\delta)/(\delta-1)$ 
have different coefficients. Let $l_1$ and $l_2$ be the minimal integers such 
that $v_{\iota'}(\delta)/(\delta-1)$ and 
$\delta^k\cdot v_{\iota}(\delta)/(\delta-1)$ 
have the coefficient $-s_{m_1}$  of $\delta^{-l_1}$ and the coefficient 
$-s_{m_2}$ of $\delta^{-l_2}$ which are different from the coefficient of $\delta^{-l_1}$ in 
$\delta^k\cdot v_{\iota}(\delta)/(\delta-1)$ 
and the coefficient of 
$\delta^{-l_2}$ in $v_{\iota'}(\delta)/(\delta-1)$ for some 
$1 \le m_1 \le n$ and $1 \le m_2 \le n$ 
respectively. Note that $s_1 > s_2 > \cdots > s_n$ and 
$\varepsilon_{m+1}(\iota'')- \varepsilon_m(\iota'') \ge 3$ 
for any $m \ge 1$ and $\iota'' \in \mathcal{K}(n)$. 
Thus, $v_{\iota'}(\delta)/(\delta-1) - 
\delta^k\cdot v_{\iota}(\delta)/(\delta-1)=0$ 
satisfies the estimates
\[
s_{m_1} \delta^{-l_1}-s_{m_2} \delta^{-l_2} - s_{1} 
\frac{\delta^{-l_2}}{\delta^3-1} < 
\frac{v_{\iota'}(\delta)}{\delta-1} - 
\frac{\delta^k \cdot v_{\iota}(\delta)}{\delta-1}
< s_{m_1} \delta^{-l_1}-s_{m_2} \delta^{-l_2} + s_{1} 
\frac{\delta^{-l_1}}{\delta^3-1}. 
\]
If $l_1 > l_2$, then it follows that 
\[
0 < s_{m_1} \delta^{-l_1}-s_{m_2} \delta^{-l_2} + s_{1}  
\frac{\delta^{-l_1}}{\delta^3-1} < 
s_{1} \delta^{-l_1}-s_{1} z_1(n)^{n-1} \delta^{-l_1+1} + s_{1} 
\frac{\delta^{-l_1}}{\delta^3-1} < s_{1} \delta^{-l_1} g_1(n),  
\]
which contradicts Lemma \ref{lem:g}. On the other hand, if $l_1 = l_2$ and 
$m_1 > m_2$, then we have 
\[
0 < s_{m_1} \delta^{-l_1}-s_{m_2} \delta^{-l_1} + s_{1} 
\frac{\delta^{-l_1}}{\delta^3-1} < 
s_{1} z_2(n)^{m_1-1} \delta^{-l_1}-s_{1} z_1(n)^{m_1-2} 
\delta^{-l_1} + s_{1} \frac{\delta^{-l_1}}{\delta^3-1} < 
- s_{1} \delta^{-l_1} g_2(n),  
\]
where the last inequality is a consequence of the fact that 
$z_2(n)^{m_1-1}-z_1(n)^{m_1-2}= - z_2(n)^{m_1-2}$
$\bigl((\frac{z_1(n)}{z_2(n)})^{m_1-2}- z_2(n) \bigr)$ 
is monotone increasing with respect to $m_1$ since 
$0 < z_2(n),\frac{z_1(n)}{z_2(n)} < 1$ and 
$(\frac{z_1(n)}{z_2(n)})^{m_1-2}- z_2(n) > \frac{g_2(n)}{z_2(n)^{m_1-2}} >0$. 
This contradicts Lemma \ref{lem:g}. In a similar manner, if $l_1 < l_2$, 
then it follows that 
\[
0 > s_{m_1} \delta^{-l_1}-s_{m_2} \delta^{-l_2} - s_{1} 
\frac{\delta^{-l_2}}{\delta^3-1} > 
s_{1} z_1(n)^{n-1} \delta^{-l_2+1}-s_{1} \delta^{-l_2} - s_{1} 
\frac{\delta^{-l_2}}{\delta^3-1} > -s_{1} \delta^{-l_2} g_1(n),  
\]
which is a contradiction. On the other hand, if $l_1 = l_2$ and 
$m_1 < m_2$, then we have 
\[
0 > s_{m_1} \delta^{-l_2}-s_{m_2} \delta^{-l_2} - s_{1} 
\frac{\delta^{-l_2}}{\delta^3-1} > 
s_{1} z_1(n)^{m_2-2} \delta^{-l_2}-s_{1} z_2(n)^{m_2-1} 
\delta^{-l_2} - s_{1} 
\frac{\delta^{-l_2}}{\delta^3-1} > s_{1} \delta^{-l_2} g_2(n),  
\]
which is also a contradiction. Thus, 
$v_{\iota'}(\delta)/(\delta-1)$ and 
$\delta^k \cdot v_{\iota}(\delta)/(\delta-1)$ 
have the same coefficients. In particular, we have 
$i_m'=i_m$ and 
$\varepsilon_m(\iota')=\varepsilon_m(\iota)-k$ 
for any $m \ge 1$. 
The relations $\varepsilon_m(\iota')=\varepsilon_m(\iota)-k$ 
yield $\kappa(\iota')=\varepsilon_1(\iota')=
\varepsilon_1(\iota)-k=\kappa(\iota)-k$ and 
$\kappa(\sigma^m(\iota'))=\kappa(\sigma^m(\iota))$ for any $m \ge 1$. 
Now, $L \ge 1$ is chosen so that $\sigma^L=\mathrm{id}$. Then we have 
$i'=i_L'=i_L=i$ and 
$\kappa(\iota')=\kappa(\sigma^L(\iota'))=\kappa(\sigma^L(\iota))=\kappa(\iota)$, 
which shows that $k=0$. 
Therefore, the proposition is established. 
\hfill $\Box$ \\[4mm] 
{\it Proof of Proposition \ref{pro:cri}}. 
From Proposition \ref{pro:coeffi}, if the relation 
$\delta^k \cdot v_{\overline{\iota}}(\delta)=v_{\overline{\iota}'}(\delta)$ 
holds, then one has $k=0$, $i_m=i_m'$ and 
$\kappa(\sigma^m(\iota))=\kappa(\sigma^m(\iota'))$ for 
any $m \ge 0$. Conversely, it is easily seen that if 
$i_m=i_m'$ and 
$\kappa(\sigma^m(\iota))=\kappa(\sigma^m(\iota'))$ for any $m \ge 0$ 
then $v_{\iota}(\delta)=v_{\iota'}(\delta)$ holds. 
In particular, it follows from Lemma \ref{lem:indA} that 
$\overline{\Gamma}_2(\tau) \cap P(\tau) = \{ 
\alpha_{\iota,\iota'}^0 \, | \, 
i_m=i_m', \, 
\kappa(\sigma^m(\iota))=\kappa(\sigma^m(\iota')), \, m \ge 0 \}$ and hence 
$(\overline{\Gamma}_2(\tau) \cap P(\tau)) \cap \Gamma_2^1(\tau) = \emptyset$. 
Moreover, if $\tau$ satisfies condition $(3)$ in 
Theorem \ref{thm:main3}, then any element 
$\alpha_{\iota,\iota'}^0 \in \overline{\Gamma}_2(\tau) \cap P(\tau)$ 
does not belong to $\Gamma_2^2(\tau)$, which shows that 
$\Gamma_2(\tau) \cap P(\tau) = \emptyset$. 
Therefore we establish Proposition \ref{pro:cri}. 
\hfill $\Box$ \par\medskip\noindent
{\bf Acknowledgment}. I would like to thank Professor Katsunori Iwasaki 
for numerous comments and kind encouragement. 
I am grateful to Professor Eric Bedford, Professor Julie Deserti, Professor Jeffrey Diller, Professor Charles Favre and 
Professor Curtis T. McMullen for useful advices, and 
to Professor Yutaka Ishii, Professor Eiichi Sato and 
Professor Mitsuhiro Shishikura for their encouragement and interest in this work. 

\end{document}